%% file: main.tex
\documentclass[12pt,leqno,oneside]{amsart}

\usepackage{amssymb,amsmath,amsthm,nccmath}
\usepackage{graphicx,color,textcomp}
\usepackage[dvipsnames]{xcolor}
\usepackage{stackrel}
\usepackage{mathtools}
\usepackage[toc,page]{appendix}
\usepackage{comment}
\usepackage[T1]{fontenc} \usepackage[utf8]{inputenc} 	
\usepackage[english]{babel} 							
\usepackage[left=2cm,top=2.5cm,bottom=2.5cm,right=2cm]{geometry}
\usepackage[colorlinks,citecolor=blue]{hyperref}
\everymath{\displaystyle}

\usepackage{tikz-cd} 
\usepackage{tikz}
\usepackage{mathdots}
\usepackage{yhmath}
\usepackage{cancel}
\usepackage{color}
\usepackage{siunitx}
\usepackage{array}
\usepackage{multirow}
\usepackage{amssymb}
\usepackage{tabularx}
\usepackage{extarrows}
\usepackage{booktabs}
\usetikzlibrary{fadings}
\usetikzlibrary{patterns}
\usetikzlibrary{shadows.blur}
\usetikzlibrary{shapes}
\usepackage{overpic}

\def \R {{\mathbb R}}

\def \Z {{\mathbb Z}}
\def \C {{\mathbb C}}

\newcommand{\uguale}{\stackrel{.}{=}} 

\newtheorem{theorem}{Theorem}[section]
\newtheorem{lemma}[theorem]{Lemma}

\newtheorem{definition}[theorem]{Definition}
\newtheorem{proposition}[theorem]{Proposition}

\newtheorem{remark}[theorem]{Remark}

\newtheorem{corollary}[theorem]{Corollary}

\newtheorem{conjecture}{Conjecture}
\newtheorem*{0theorem}{Theorem}
\newtheorem*{question}{Question}

\title[On the Birkhoff conjecture for Kepler billiards]{On the Birkhoff conjecture for Kepler billiards}
\author{S. Baranzini, V. Barutello, I. De Blasi, S. Terracini}

\address{University of Turin, Department of Mathematics, via Carlo Alberto 10, Turin, Italy}

\date{\today} 
\thanks{}


\makeatletter
\@namedef{subjclassname@2020}{\textup{2020} Mathematics Subject Classification}
\makeatother


\begin{document}

\begin{abstract}
We investigate the integrability of Kepler billiards—mechanical billiard systems in which a particle moves under the influence of a Keplerian potential and reflects elastically at the boundary of a strictly convex planar domain. Our main result establishes that, except possibly for one location of the gravitational center, analytic integrability at high energies occurs only when the domain is an ellipse and the center is placed at one of its foci. This provides a partial affirmative answer to a Keplerian analogue of the classical Birkhoff-Poritsky Conjecture.

Our approach is based on the construction of symbolic dynamics arising from chaotic subsystems that emerge in the high-energy regime. Depending on the geometric configuration of the boundary and the location of the attraction center, we construct three types of symbolic dynamics by shadowing chains of punctured Birkhoff-type trajectories. These constructions yield subsystems conjugated to Bernoulli shifts, implying positive topological entropy and precluding analytic integrability.

We further analyze the notion of \emph{focal points of the second kind}, showing that in real-analytic, non-elliptic domains there can be at most one such point, while ellipses are the only domains admitting two—coinciding with their classical foci. Finally, we demonstrate the existence of an infinite-dimensional family of non-elliptic domains possessing a focal point of the second kind, and conclude with numerical simulations illustrating chaotic behavior in such cases.
\end{abstract}
    	\maketitle
    \tableofcontents

\section{Introduction}
    \input{intro}
\section{Preliminaries and asymptotic estimates} 
    \input{preliminaries}

\section{Existence of a symbolic dynamics} 
\label{sec:dyn_sym}
\input{symbolic_dynamics}

\section{A rigidity result} 
\label{sec:rigidity}
\input{rigidity}

\appendix
\section{Proof of Lemma \ref{lem:orario_antiorario}}
\label{app:proof_Lemma}
\input{appendix}
\section{Proofs of Section \ref{subsec:critical_points}}
\input{appendix2}


\end{document}

%% file: intro.tex
In a \textit{Birkhoff billiard}, \cite{birkhoff1927acta}, a point particle moves freely (i.e., with constant velocity) within a planar domain and undergoes elastic reflections upon encountering the boundary. At each reflection point, the direction of motion changes according to the law of specular reflection: the angle of incidence equals the angle of reflection. Typically, the term \textit{Birkhoff billiard} refers to billiards in strictly convex planar domains with smooth boundary. The associated dynamics raises fundamental questions concerning integrability, existence, and stability of periodic orbits, ergodic properties, and the structure of the corresponding billiard map, which encodes the sequence of impact points and directions.

Let \( \Omega \subset \mathbb{R}^2 \) be a bounded domain with smooth boundary \( \partial\Omega \). At a point of impact \( p \in \partial\Omega \), the particle has a unit inward-pointing velocity vector \( v \in \mathbb{S}^1 \). 
Then, it travels along the straight line determined by \( (p, v) \) until it reaches another point \( p' \in \partial\Omega \). Here, it reflects its velocity about the tangent line to $\partial \Omega$ at $p'$, obtaining a new inward-pointing vector \( v' \). 
This defines the \textit{billiard map}
\(
T: \mathcal C \to \mathcal C, \quad T(p, v) = (p',v').
\)
To describe it in coordinates, we parametrize the boundary \( \partial\Omega \) by arclength \( \xi\in \mathbb{S}^1 \simeq \mathbb{R}/L\mathbb{Z} \), where \( L = \text{length}(\partial\Omega) \), and the inward directions by the angle \( \alpha \in (0, \pi) \) they make with the positively oriented tangent to \( \partial\Omega \) at \( \gamma(\xi) \). The phase space of the billiard map becomes the cylinder \( \mathbb{S}^1 \times (0, \pi) \). In these coordinates, the billiard map \( T \) is an \textit{area-preserving twist map}.

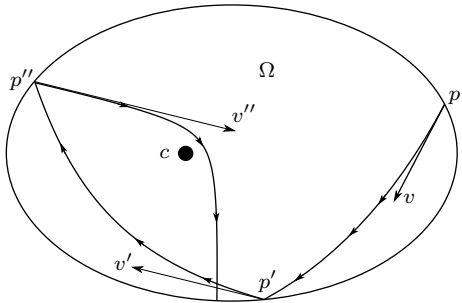
\begin{figure}[t!]
	\begin{overpic}[width=0.35\linewidth]{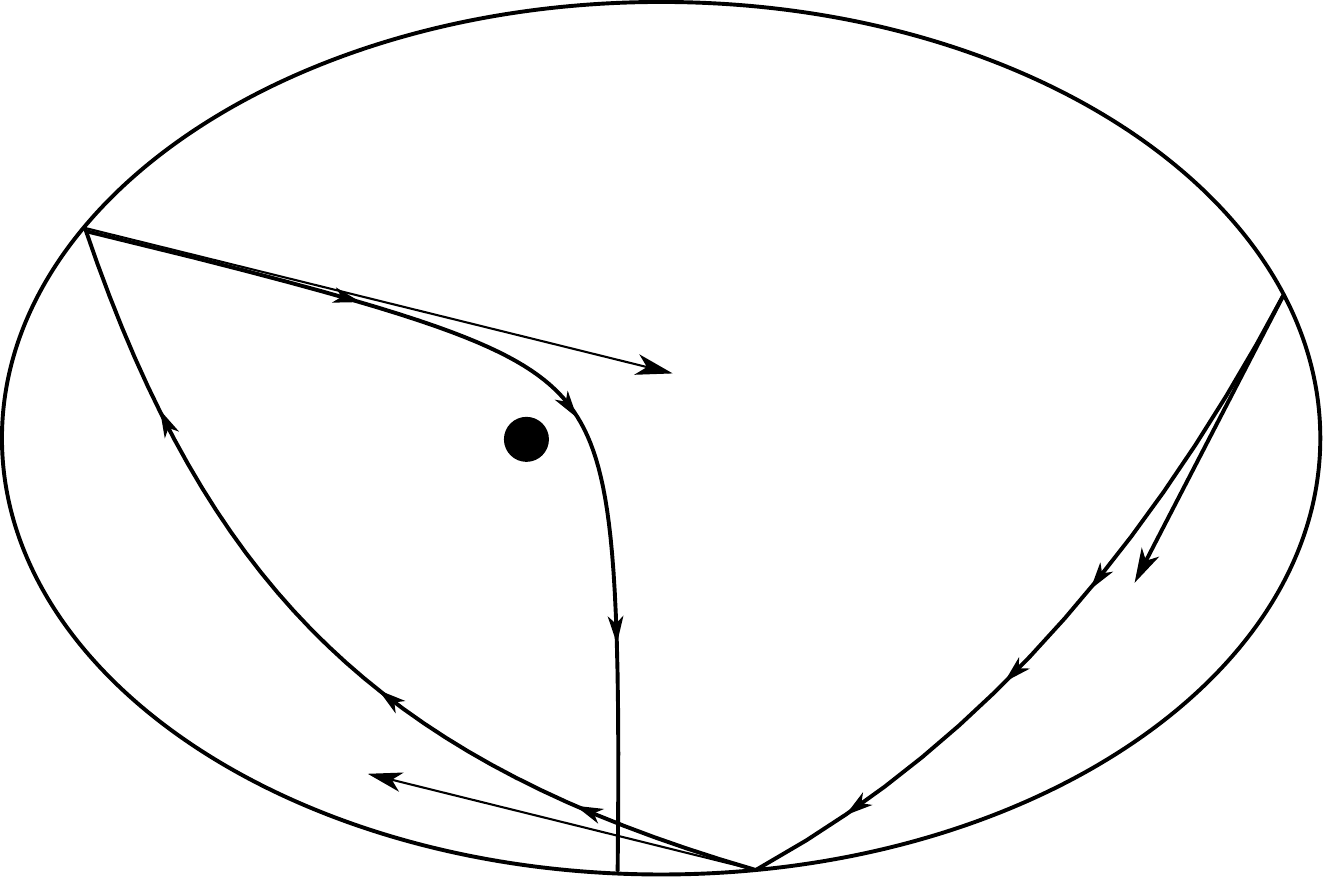}
		\put (98, 45) {\tiny$p$}
		\put (88, 22) {\tiny$v$}
		\put (56, 3) {\tiny$p'$}
		\put (24, 7) {\tiny$v'$}
        \put (1, 48) {\tiny$p''$}
		\put (50,40) {\tiny$v''$}
        \put (34,32) {\tiny$c$}
        \put (56, 50) {\tiny$\Omega$}	
	\end{overpic}
    \hspace{2cm}
	\caption{Kepler billiard with gravitational center at $c$ and positive energy.}
	\label{fig:riflessione1}
\end{figure}

This framework extends naturally to \textit{mechanical billiards}, in which the particle moves under the influence of a conservative force field \( F = \nabla V \) within the domain \( \Omega \), and still undergoes elastic reflections at the boundary (cfr. \cite{fedorov2001ellipsoidal,pustovoitov2019topological,pustovoitov2021topological}). A notable example of such billiard traces back to Boltzmann (see \cite{boltzmann1868}). A particle moves under the action of a radially symmetric potential {of the form $1/r + \beta/r^2$ for some $\beta\ge0$} and bounces elastically against a flat table not passing through  the origin. {Boltzmann introduced this example as a simple toy model to illustrate his so-called \textit{ergodic hypothesis}, claiming that the reflection with the boundary would completely destroy the integrability of the system. In recent years, Gallavotti and Jauslin (see \cite{gallavotti2020theorem}) showed that for $\beta=0$ the system is actually integrable and Felder proved that, for $\beta$ small, generic trajectories are quasi-periodic and so the system cannot be ergodic (see \cite{Felder2021}). 
Later on, other works (see \cite{zhaoBoltzmann, Gasiorek2024111}) analyzed the model further, providing a different approach to integrability based on projective dynamics and a detailed description of the phase space.}\\

In this paper, we consider \textit{Kepler billiards}, where the potential energy is of the form 
\begin{equation}
\label{eq:kepler_potential}
V(z) = \frac{\mu}{|z - c|}, \quad z \in {\mathbb R}^2 \setminus\{c\},
\end{equation}
for some fixed center \( c \in \Omega \) and \( \mu > 0 \). Between two consecutive bounces, the particle follows a solution of the Kepler problem, namely, an arc of conic section. At each impact with \( \partial\Omega \), the trajectory reflects according to the standard elastic reflection rule and total energy is preserved (as illustrated in Figure \ref{fig:riflessione1}). 
In contrast with standard Birkhoff billiards, in Kepler billiards the presence  of a center of attraction and {the energy level} significantly influence the geometry of the particle trajectories. As a result, two natural additional parameters emerge in the system: the position of the center and the value of the total mechanical energy.
{Moreover, the system is singular as} the potential becomes infinite at the point $c$.
We will take advantage of Levi-Civita regularization to obtain a complete analytic flow (see \cite{Levi-Civita}). This completion, however, is not sufficient for the corresponding Kepler billiard map to be {globally well-defined on $\mathcal{C}$ at all energies. It gives a well defined area preserving map (see Lemma \ref{lem:well_defined}) provided that the energy is high enough and $\partial \Omega$ is strictly convex.}

The dynamics of Birkhoff billiards are particularly well understood in the case of circular and elliptical tables, owing to their remarkable optical and geometric properties. In the circular case, the billiard map reduces to a rigid rotation in the phase space: the angle of reflection remains constant, resulting in a globally integrable system with all trajectories lying on invariant circles. Elliptic billiards also exhibit complete integrability: the phase space is globally foliated by smooth invariant curves, except for a 2-periodic orbit and its homoclinic connections. In both cases, the existence of a family of invariant curves (and associated caustics) underpins the integrable structure.

A central open question in this context is the \textit{Birkhoff-Poritsky Conjecture}    \cite{birkhoff1927acta,zbMATH03056709}, which posits that
\begin{conjecture}[Birkhoff-Poritsky]
Among strictly convex, planar, smooth, bounded domains, the only integrable Birkhoff billiards are those bounded by ellipses.
\end{conjecture}

Here, integrability entails the existence of a full measure set of invariant curves in the phase space, associated with a foliation by caustics to which billiard trajectories remain tangent under  reflection. In elliptical domains, this structure arises naturally from the classical optical property of confocal conics: any trajectory initially tangent to a given confocal conic remains tangent to it after each reflection. This foliation guarantees the quasi-periodic nature of most trajectories and serves as the hallmark of integrability in this setting. The conjecture is still open in general, but several partial results and local rigidity theorems have been proved notably by Bialy \cite{zbMATH00561292}, Avila, De Simoi and Kaloshin \cite{avila2016integrable}, Kaloshin and Sorrentino \cite{MR3868417,MR3815464}, Guan, Kaloshin and Sorrentino \cite{MR3788206}, Bialy and Mironov \cite{BiaMir22} and the references therein for a more comprehensive treatment of this subject. Near ellipses, local versions of the conjecture have been affirmatively proven: small smooth perturbations of ellipses that remain integrable must be ellipses. For centrally symmetric tables, even stronger results are known (see \cite{BiaMir22}).

In the context of Kepler billiards, a result by Panov in \cite{panov1994elliptical} (see also \cite{KOZLOV19951,zhaoBoltzmann,takeuchi2021conformal}) establishes a fundamental instance of integrability

\begin{0theorem}[Panov, \cite{panov1994elliptical}]
Elliptic Kepler billiards are integrable at all energy levels when the center of attraction is located at one of the foci of the ellipse.
\end{0theorem}

The result proved by Panov, later generalized by Takeuchi and Zhao (see \cite{takeuchi2021conformal, takeuchi2022projective}) to combination of confocal conics,  is based on the conjugation between Hooke and Kepler systems via the complex square mapping. Taking advantage of the conformality of the latter, integrals of an elastic billiard can be translated into integrals of Kepler ones, provided the boundaries are transformed accordingly. Notably, Boltzmann's example for $\beta=0$ falls into the more general class considered there.

Panov's result naturally leads to an analogue of the classical Birkhoff-Poritsky Conjecture in the Keplerian setting. We are thus prompted to ask

\begin{question}
Are there other examples of smooth integrable Kepler billiards, beyond  {ellipses with the center of attraction at a focus}?
\end{question}

The aim of this work is to contribute answering to this question, thereby partially establishing a Keplerian analogue of the Birkhoff-Poritsky Conjecture for domains with analytic boundary.  In this regard, we will demonstrate the following result.

\begin{theorem}\label{thm:integrable=ellipse}
Let \( \Omega \subset \mathbb{R}^2 \) be a strictly convex, bounded domain with real-analytic boundary.  Then, if $\partial\Omega$ is not an ellipse and except for possibly one position of the center of attraction, the Kepler billiard in $\Omega$ is not analytically integrable at high energies. If $\partial\Omega$ is an ellipse, these exceptional positions are two and coincide with its foci.
\end{theorem}

A point $c$ is termed \textit{focal} for a domain $\Omega$ if any Birkhoff trajectory passing through this point returns back to it after two bounces (see Definition \ref{def:focal}). This implies that there is an invariant {curve} for the associated Birkhoff billiard map determined by rays through $c$ and their reflections.   
Of course, for an ellipse, the two foci are the only focal points. Examples of non elliptical domains with one focal point are discussed in Section \ref{subsect:string}; their existence implies the optimality of Theorem  \ref{thm:integrable=ellipse}.
As a simple by-product of the Theorem \ref{thm:integrable=ellipse}, we have the following corollary.
\begin{corollary}\label{thm:integrable=ellipsesymmetry}
Under the same assumptions of Theorem \ref{thm:integrable=ellipse}, assume moreover that $\Omega$ possesses a central or a cyclic rotational symmetry. Then, {any associated Kepler billiard in $\Omega$} is not analytically integrable at high energies unless $\partial\Omega$ is an ellipse with the center of attraction located at one of its foci.
\end{corollary}

It is worth emphasizing that our objective is not to identify instability regions through perturbations of the associated Birkhoff billiard away from the keplerian center. Rather, our approach to demonstrating non-integrability is based on the identification of a subsystem that exhibits symbolic dynamics at high energy levels, thereby precluding analytic integrability. Our result extends beyond the mere failure of integrability away from the singularity, by isolating a chaotic subsystem whose trajectories undergo a sequence of close encounters with the gravitational center, which effectively acts as a scatterer.

Symbols will describe different types of close encounters with the Kepler center. Our approach is based on a \textit{shadowing} technique reminiscent of the one adopted by Bolotin and Bolotin-MacKay in \cite{Bol2017, bolotin2000periodic}.  
To better explain our result, we need to introduce the prototype of a chaotic system: the Bernoulli shift. Let $\mathcal S$  be a discrete set endowed with the discrete metric and consider the set of bi-infinite sequences
    	\[
    	\mathcal S^\Z\uguale\{(s_k)_{k\in\Z}:\,s_k\in {\mathcal S},\ \forall\,k\},
    	\]
endowed with the distance:
    	\[
    	d((s_k)_k,(t_k)_k)\uguale\sum\limits_{k\in\Z}\frac{1}{4^{|k|}} d(s_k,t_k),\quad\mbox{for}\ (s_k)_k,(t_k)_k\in {\mathcal S}^\Z.
    	\]
   \textit{the Bernoulli right shift}  is the map:
    \[
    \begin{split}
    \sigma_r\colon \mathcal S^\Z\longrightarrow \mathcal S^\Z \;:&\;(s_k)\mapsto\,T_r((s_k)_k)\uguale(s_{k+1})_k\\
    (\dots,s_{-k\phantom{+1}},\dots,s_{-1},s_0,s_1,\dots,s_{k\phantom{+1}}\dots)&\mapsto (\dots,s_{-k+1},\dots,s_{-0},s_1,s_2,\dots,s_{k+1}\dots)
    \end{split}    	
    \]
\begin{definition}\label{def:sym_dyn}
We say that a dynamical system $(X,F)$ has a symbolic dynamics with set of symbols $S$ if there exists continuous surjective map $\pi: X\to \mathcal S$ such that the diagram 
\[
\begin{tikzcd}
	X \arrow{r}{F} \arrow{d}{\pi} & X \arrow{d}{\pi} \\
\mathcal{S}^\Z \arrow{r}{\sigma_r}	& \mathcal{S}^\Z
\end{tikzcd}
\]
commutes. 
	\end{definition}

\begin{theorem}\label{thm:symbolic}
Let \( \Omega \subset \mathbb{R}^2 \) be a strictly convex, bounded domain with real-analytic boundary.  Then, if $\partial\Omega$ is not an ellipse
and except for possibly one position of the center of attraction, there is $\overline{h}>0$ such that, for any energy $h\geq\overline{h}$ the  associated Kepler billiard map has a subsystem displaying a symbolic dynamics. If $\partial\Omega$ is an ellipse, these exceptional positions are two and coincide with its foci.
\end{theorem}

Again, we have the following corollary for domains with cyclic symmetries.

\begin{corollary}\label{thm:symbdym=ellipsesymmetry}
Under the same assumptions of Theorem \ref{thm:symbolic}, assume moreover that $\Omega$ possesses a central or a cyclic rotational symmetry. Then, if $\partial\Omega$ is not an ellipse, for any choice of the gravitational center there is $\overline{h}>0$ such that, for any anergy $h\geq\overline{h}$ the  associated Kepler billiard map has a subsystem displaying a symbolic dynamics. \end{corollary}

The presence of a subsystem which is semi-conjugated with the Bernoulli shift with two or more symbols has profound implications on the dynamics. As stated in Corollary \ref{cor:katok}, it yields positive topological entropy and non integrability by analytic first integrals.
Let us remark that our results hold for sufficiently large energies. We refer to the numerical simulations in \cite{BCDeB_Nonlin} to have some examples on elliptic domains of the energy necessary to see chaotic phenomena. In the light of our further numerical simulation here (see Figure \ref{fig:primo_ritorno}) we formulate the following conjecture.

\begin{conjecture}[Birkhoff-Poritsky conjecture for Kepler Billiards]\label{conj:main}
Among all domains \( \Omega \subset \mathbb{R}^2 \) strictly convex, bounded  and with real-analytic boundary, the only Kepler billiards which are integrable at all energy levels are ellipses with the gravitational center located at one of its foci.
\end{conjecture}

\subsection*{Outline of the proofs}
Our construction naturally begins with a proper definition of the first return map for a Kepler billiard, denoted by $T_K$, in strictly convex domains. Analogously to the classical setting, the map $T_K$ is described via suitable generating functions, which encode the Jacobi length of Keplerian arcs connecting two points (see Definition \ref{def:generatrici}). In order to handle arcs colliding with the gravitational center we make use of Levi-Civita regularization. The existence of two non-homotopically equivalent arcs joining a given pair of points motivates the use of a multi-valued generating function: distinct branches of this function define $T_K$ in different regions of the phase space. 

For large energy levels ($h \gg 1$), the asymptotic length associated with this multi-valued generating function converges to the Euclidean length of the segments joining the endpoints, possibly passing through the center of attraction (see Eqs. \eqref{eq:asintotiche} and \eqref{eq:Lzazc}). This asymptotic behavior motivates the introduction—considered as a limiting case—of a variation of the classical Birkhoff billiard, which we refer to as the \emph{punctured Birkhoff billiard}, explicitly accounting for the position of the attraction center.

\begin{figure}[t!]
\begin{overpic}[width=0.2\linewidth]{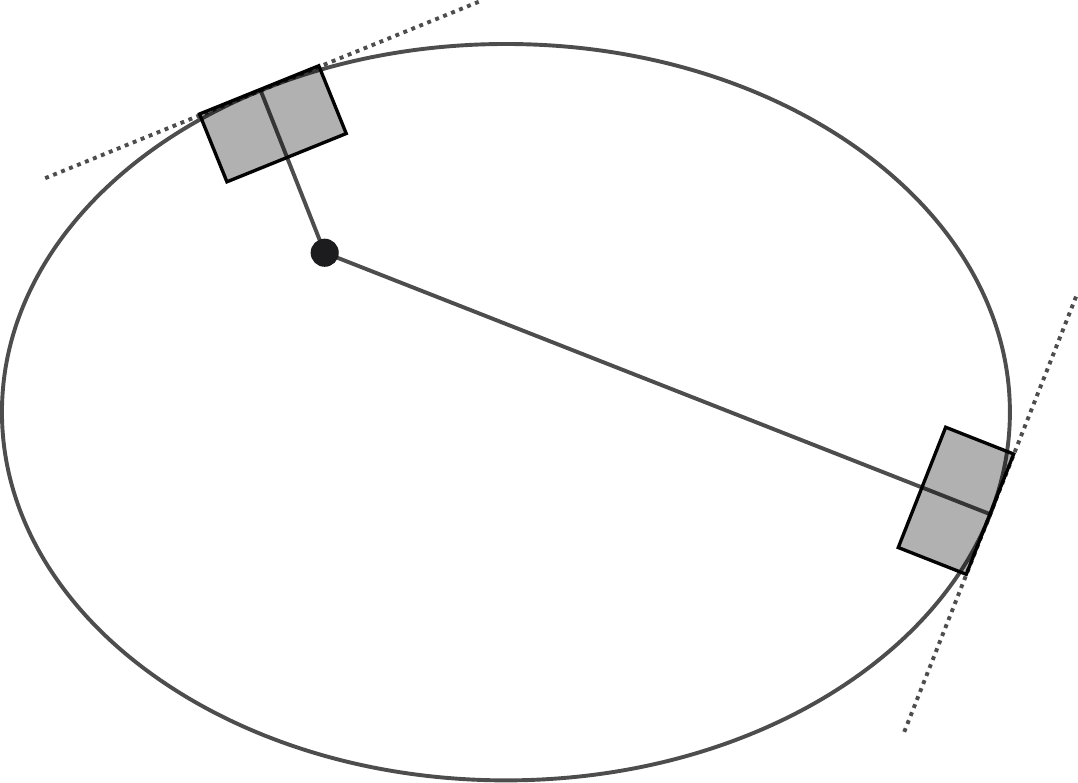}
	\end{overpic}
    \hspace{1cm}
    \begin{overpic}[width=0.25\linewidth]{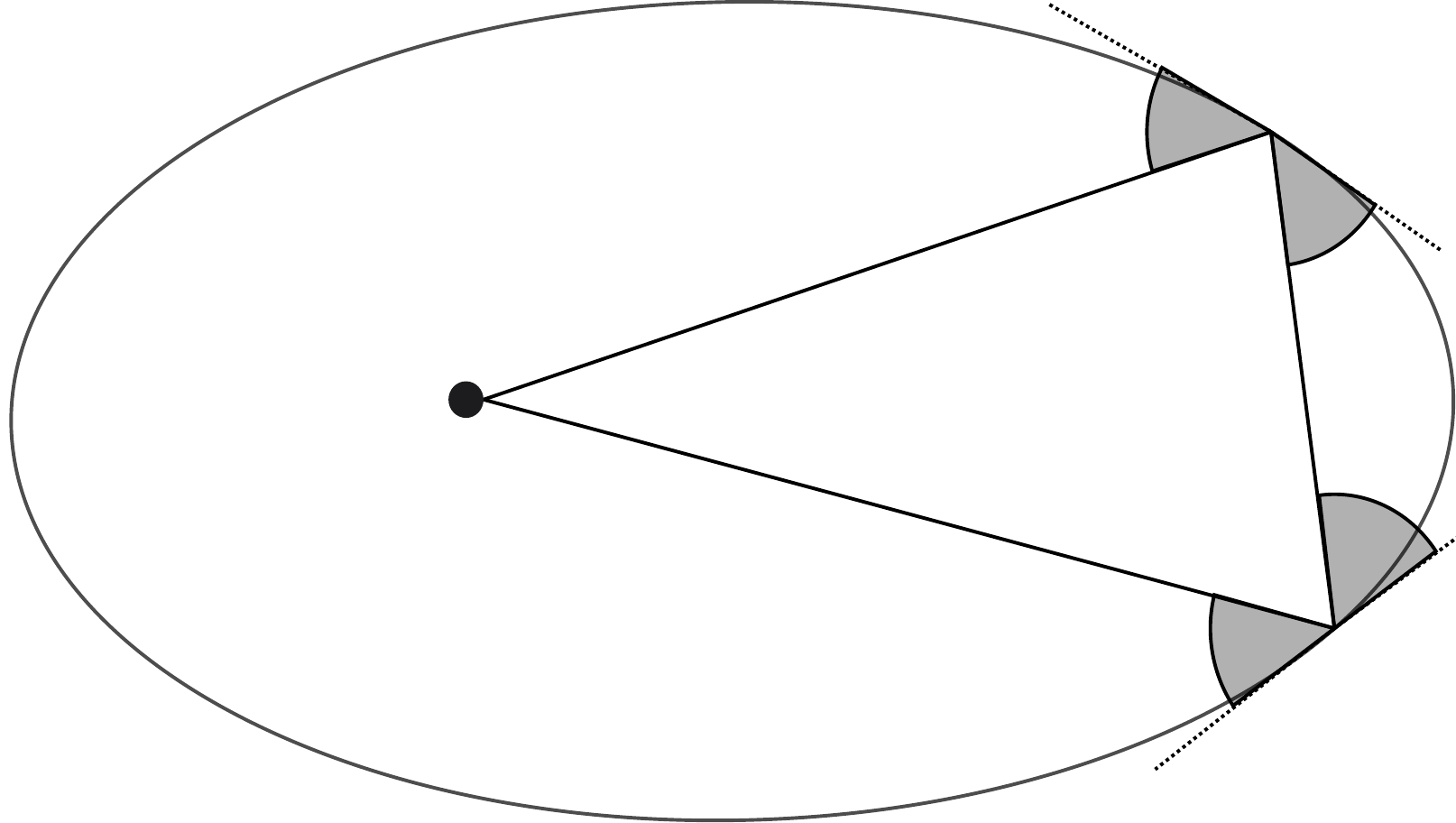}
	\end{overpic}
    \hspace{1cm}
    \begin{overpic}[width=0.3\linewidth]{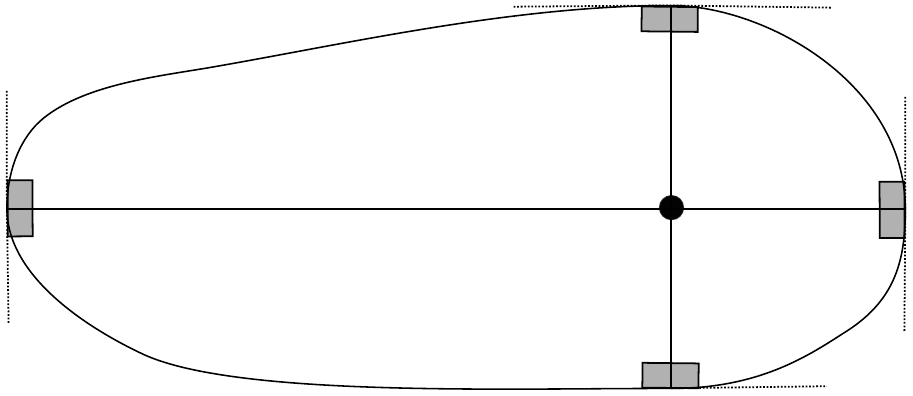}
	\end{overpic}
\caption{Elements to shadow in the three different constructions of symbolic dynamics.}	\label{fig:3dyn_symb}
\end{figure}

In light of the aforementioned asymptotic behavior, in Section~\ref{sec:dyn_sym} we investigate topologically stable trajectories of the punctured Birkhoff billiard, featuring one or two impact points, with the aim of constructing three distinct types of symbolic dynamics, each determined by the geometry of the boundary. These constructions differ in the nature of the symbolic elements that compose the chains to be shadowed (see Figure~\ref{fig:3dyn_symb}).

The first case concerns pairs of non-antipodal rays emanating from the center of attraction and striking the boundary orthogonally. The second involves non-degenerate punctured Birkhoff triangles, with one vertex located at $c$ and the remaining two on the boundary, where elastic reflections occur. Finally, the third construction shadows triangles of the same type, but with vanishing area. 

We note that the first construction has already been carried out in~\cite{IreneSusViNEW}. Furthermore, we stress that each of the three variants involves trajectories that interact with the center of attraction, resulting in a dynamics that differs fundamentally from those observed in the classical setting of Birkhoff billiards.
 
Since criticality implies elastic reflection at the boundary, the shadowed chains of geometric elements correspond to critical points of certain length functionals: in the first case, the distance from the point $c$ to points on $\partial \Omega$; in the latter two, the perimeter of triangles with one vertex at $c$ and the other two on the boundary. 

Accordingly, Section~\ref{subsec:critical_points} is devoted to the analysis of these perimeter functions, denoted by $\Psi$ (see Eqs. \eqref{eq:Psi}, \eqref{eq:Psi_a} and \eqref{eq:Psi_c}). Special variants are considered in order to deal with zero area triangles, taking into account the multivalued nature of our generating functions; therefore, we establish the existence of topologically stable critical points of $\Psi$ and its variants, and provide a characterization of their structure (see Proposition \ref{prop:critici_Psi} and Lemma \ref{lemma:grado_psi_A}). The property of topological robustness is expressed in terms of the topological degree, and the analysis is fundamentally  carried within the perturbative regime associated with high-energy limits.

Once all limiting elements have been introduced and examined, we proceed to construct the symbolic dynamics as in Definition~\ref{def:sym_dyn}. The first step entails the introduction of an appropriate set of symbols, which, for simplicity, is assumed to have cardinality two. We then construct periodic trajectories corresponding to periodic symbolic sequences. In Sections \ref{subsub:nonzero} and \ref{subsub:zero}, the existence of such trajectories is reduced to the problem of locating critical points of a suitable length functional, whose existence is ensured by the product and homotopy invariance properties of the topological degree, applied to the critical points of the perimeter functions.
The second step involves extending the construction to nonperiodic trajectories, which is carried out via standard diagonal arguments in the proof of Theorem~\ref{thm:dyn_sym}.

This result provides geometric conditions on the boundary $\partial \Omega$ and on the position of the center under which symbolic dynamics can be established. To this end, we introduce the notion of \emph{focal point} and classify points in $\Omega$ as either \emph{of the first kind} or \emph{of the second kind} (see Definitions \ref{def:focal} and \ref{def:firstkind}). Points of the first kind correspond to locations of the attracting center for which symbolic dynamics can be constructed following the approach in~\cite{IreneSusViNEW}. Conversely, non-focal points of the second kind allow the construction of symbolic dynamics via topologically stable punctured Birkhoff triangles described earlier. Thus, Theorem~\ref{thm:dyn_sym} is inconclusive only in the case of focal points of the second kind.

In Section~\ref{sec:rigidity} we analyze the possible focal points of the second kind in a strictly convex analytic domain. We show that if the domain is not an ellipse, there can be at most one such point (see Theorem~\ref{thm:at_most_one_point}). Finally, we obtain that non-circular ellipses are the only domains that admit exactly two focal points of the second kind, namely, their classical foci.

In the final section of the paper, Section~\ref{subsect:string}, we show that the set of domains admitting a focal point of the second kind is non-empty and contains an infinite dimensional family. We conclude with numerical simulations illustrating the emergence of chaotic behavior even in this exceptional case.

%% file: preliminaries.tex
{Let $\Omega \subset  \R^2$ be a bounded and strictly convex open domain and let us consider the Kepler potential defined in \eqref{eq:kepler_potential}. Without loss of generality, we can fix the origin at $c$. We consider a Kepler billiard inside $\Omega$ at positive energies, $h>0$. A particle moves along segments of Keplerian hyperbola until it reaches the boundary where it is reflected. 
As customary in billiard theory, the dynamics is described by an area preserving map $T_K:(p,v) \mapsto (p',v')$ (see Figure \ref{fig:riflessione1}). Our aim is to describe, at least locally, such a map using generating functions.\\
A natural question concerns the good definition of the map $T_K$.  Recall that, any two points $p_0,p_1 \in\partial\Omega$ can be connected by two solutions, possibly regularized, of the problem 
\begin{equation}\label{eq:Bolza}
	\begin{cases}
		z''(t) = \nabla V(z(t)), \quad &t \in [0,T],\\
		\frac12 |z'(t)|^2 - V(z(t)) = h, \quad & t \in [0,T],\\
		z(0) = p_0, \; z(T)= p_1,
	\end{cases}	
\end{equation}
for a suitable $T>0$.
{  While in the classical Birkhoff case segments connecting points of  $\partial \Omega$ are contained in $\Omega$ thanks to its strict convexity, this is not the case when we deal with non constant potentials. For this,} we need that no Kepler hyperbola starting from $\partial \Omega$ hits the boundary tangentially. { This condition is ensured by the lemma below, provided $h$ is sufficiently large.}
\begin{lemma}\label{lem:well_defined}
If $\Omega$ is strictly convex with a smooth boundary then the Kepler billiard map is well defined, provided the energy $h$ is large enough.
\end{lemma}
\begin{proof}
It is enough to prove that Keplerian arcs starting at a boundary point will not hit tangentially the boundary again. To this aim, a good criterium is a strict ordering of the curvatures, namely that the maximal possible curvature of a Keplerian hyperbola should smaller than the minimal curvature of the boundary ($k_{M}(\textrm{hyperbola},h) < k_m(\partial \Omega)$). This inequality can be easily achieved for large values of $h$ by considering that, outside a neighborhood of the center, Keplerian arcs with high energy are practically segments and their curvature can be made as small as desired {(see also Lemmas \ref{lem: asintotiche_sus} and \ref{lem:orario_antiorario})}.
\end{proof}

Next we turn to the generating functions of a Kepler billiard with positive energy. Solutions of \eqref{eq:Bolza} can be characterized as geodesics of the Jacobi metric, defined as $g_{i,j}(z)=(h + V(z))\delta_{i,j}$. The Jacobi length of a curve $z$ defined on $[0,T]$ is given by
\begin{equation}\label{eq:Jacobi}
L(z) = \int_0^T |z'(t)| \sqrt{h+V(z(t))}\, dt.
\end{equation}
Lemmas \ref{lem: asintotiche_sus} and \ref{lem:orario_antiorario} below provide asymptotic estimates for the Jacobi length as $h\to \infty$. Such formulas hold uniformly in $h$ as long as the endpoints vary in suitable compact sets of $\R^2$. 
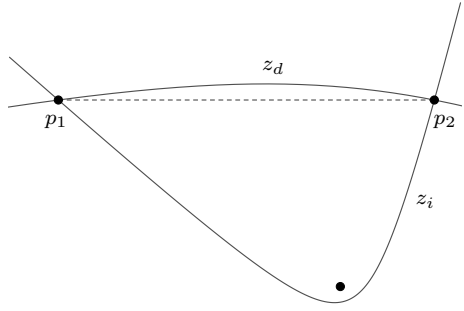
\begin{figure}[t!]
    \begin{overpic}[width=0.35\linewidth]{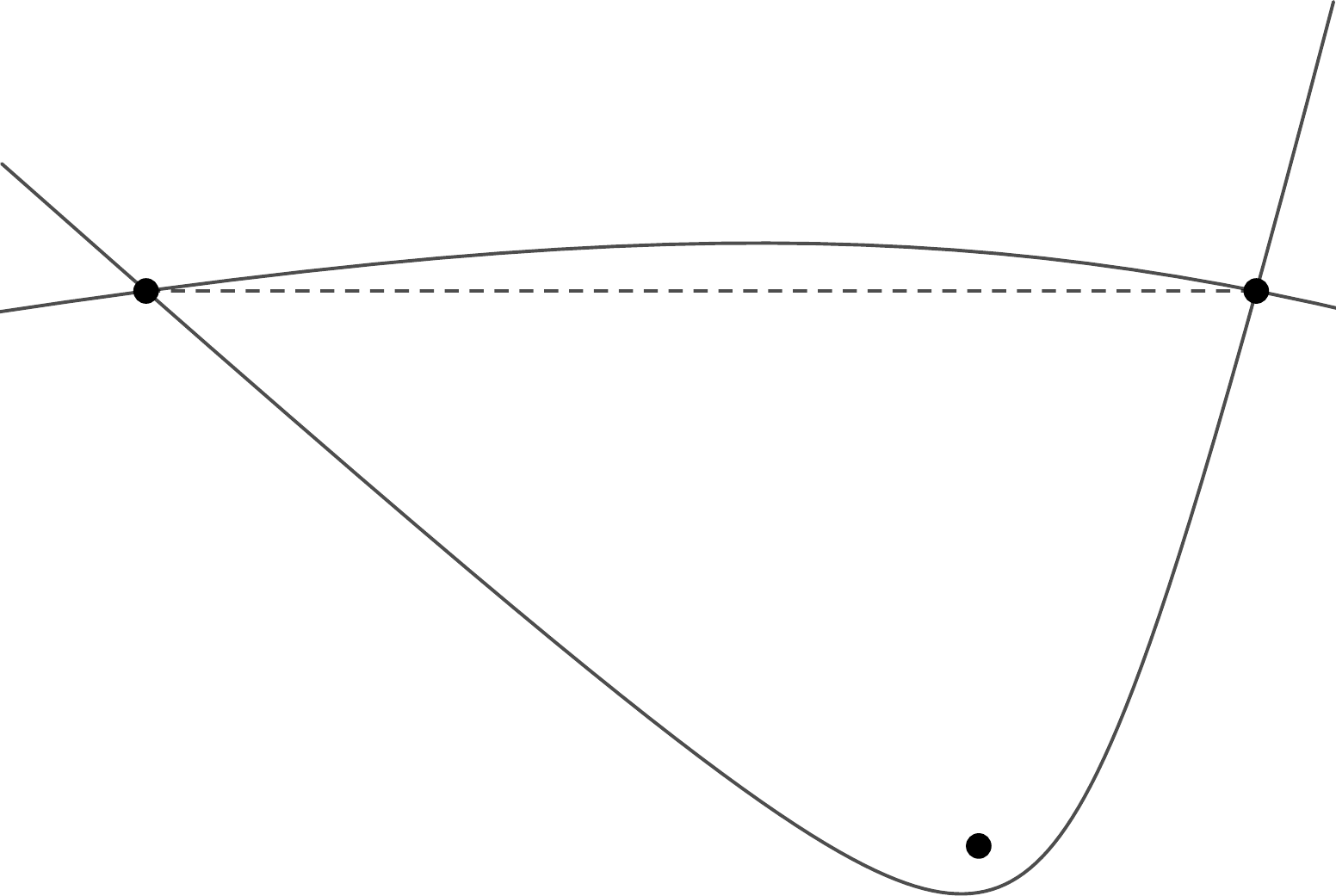}
		\put (94, 40) {\tiny$p_2$}
        \put (8, 40) {\tiny$p_1$}
		\put (90, 22) {\tiny$z_i$}
		\put (56, 52) {\tiny$z_d$}	
	\end{overpic}
	\caption{Direct and indirect arcs connecting two points in $\R^2 \setminus\{0\}$.}
	\label{fig:riflessione}
\end{figure}
\begin{lemma} \label{lem: asintotiche_sus}
	For every $\delta \in (0,1)$, let
	\[
	K_\delta = \left\{ (p_0,p_1) \in \R^2\times \R^2 \;|\; |\arg(p_0)-\arg(p_1)| \leq \pi(1-\delta), \, \delta \leq |p_i|\leq \frac1\delta,\, i=0,1 \right\}.
	\]
	Then, for every $h >0$ and $p_0,p_1 \in K_\delta$, there are two distinct solutions of \eqref{eq:Bolza}.
	One of them, $z_d(\cdot ;p_0,p_1;h)$, is homotopic to the segment $p_0p_1$ in $\R^2 \setminus \{0\}$; while the second, $z_i(\cdot ;p_0,p_1;h)$, is not (see Figure \ref{fig:riflessione}). Moreover 
	\begin{equation}\label{eq:asintotiche}
		\begin{aligned}
			L(z_d) &= \sqrt{h}|p_0-p_1| +\frac{\mu}{\sqrt{h}}g_d(p_0,p_1; h), \\
			L(z_i) &= \sqrt{h}\left(|p_0|+|p_1|\right) + \frac{\mu}{\sqrt{h}}\left(g_i(p_0,p_1; h) -\log\left(\dfrac{\mu}{2h}\right)\right),
		\end{aligned}
	\end{equation} 
	where $g_d$ and $g_i$ are $C^2$-bounded uniformly with respect to $p_0,p_1 \in K_\delta$ as $h \to +\infty$. \\
    {
    In particular, for every $\delta \in (0,1)$ and for every $\epsilon \in \left(0,\mathop{diam}(\Omega)\right)$, there exists $\bar h = \bar h ( \delta,\epsilon )$ such that for every {$(p_0,p_1) \in (\partial \Omega \times \partial \Omega) \cap K_{\delta}$}, $z_i$ is completely contained in $\Omega$. Assuming further that $|p_0-p_1| \geq \epsilon$, then also $z_d$ is completely contained in $\Omega$.}
\end{lemma}

The proof of the existence of $z_d$ and $z_i$ is a classical result by Jacobi (\cite{jacobi1837theorie}); the first asymptotic estimate follows from the definition of $L$ and the boundedness from the origin of $z_d$; the second one has been already proved in \cite{IreneSusViNEW} (see also \cite{Bol2017}). 
In accordance to the literature on this subject (see \cite{WintBook,MR4058157}), we will call $z_d$ the {\em direct arc} and $z_i$ the {\em indirect} one.
{The last part of the statement follows from uniform convergence of $z_i$ and $z_d$ to the respective segments (see \cite[p. 274]{battin}).}
Lemma \ref{lem: asintotiche_sus} does not apply to pairs of points belonging to neighborhoods of antipodal points. For this reason we state the next result, which will be proved in Appendix \ref{app:proof_Lemma}.

\begin{lemma}\label{lem:orario_antiorario}
Let $\delta\in(0,1),$ and define 
\[
\tilde K_\delta = \left\{ (p_0,p_1) \in \R^2\times \R^2 \;|\; |\arg(p_0)-\arg(p_1)| \in [\pi(1-\delta), \pi(1+\delta)], \, \delta \leq |p_i|\leq \frac1\delta,\, i=0,1 \right\}.
\]
Then, for every $h>0$, and every $p_0, p_1\in\tilde K_\delta$ there are two distinct solutions of \eqref{eq:Bolza}, called $z_a(\cdot; p_0, p_1; h)$ and $z_c(\cdot; p_0, p_1; h)$, such that 
\begin{itemize}
    \item $z_a$ has positive angular momentum, and connects $p_0$ to $p_1$ anticlockwise; 
    \item $z_c$ has negative angular momentum, and connects $p_0$ to $p_1$ clockwise. 
\end{itemize}
Moreover, once defined $f_a, f_c:\tilde K_\delta\to \R$ as 
\begin{equation}
\label{eq:def_fa_fc}
f_a(p_0, p_1) =\begin{cases}
    |p_0-p_1| \quad &\text{if }\det (p_0, p_1)\geq 0\\
    |p_0|+|p_1| \quad &\text{if }\det (p_0, p_1) < 0
\end{cases}
\quad f_c(p_0, p_1) =\begin{cases}
    |p_0|+|p_1| \quad &\text{if }\det (p_0, p_1) \geq 0\\
    |p_0-p_1|  \quad &\text{if }\det (p_0, p_1) < 0
\end{cases}, 
\end{equation}
(see Figure \ref{fig:antiorario_orario}) it holds 
\begin{equation}\label{eq:Lzazc}
L(z_{a/c})=\sqrt{h}f_{a/c}(p_0, p_1)+g_{a/c}(p_0, p_1; h), 
\end{equation}
where $g_{a/c} $ are functions of class $C^1$ in $\tilde K_\delta\times \R^+$ uniformly bounded in $p_0,p_1$. \\
In particular, if $p_0,p_1$ are in $\tilde K_\delta\cap\left( \partial \Omega\times \partial \Omega \right)$ and $h$ is sufficiently large, then both $z_a$ and $z_c$ lie inside $\Omega.$
\end{lemma}

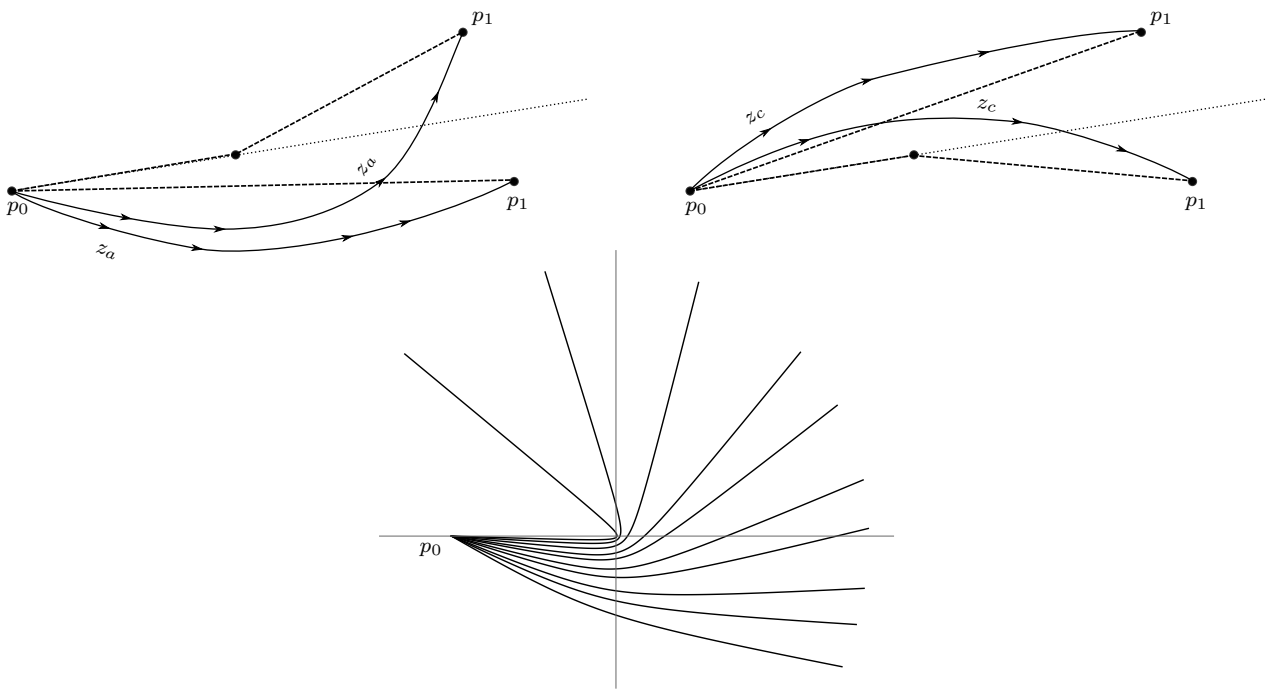
\begin{figure}[t]
	\begin{overpic}[width=0.45\linewidth]{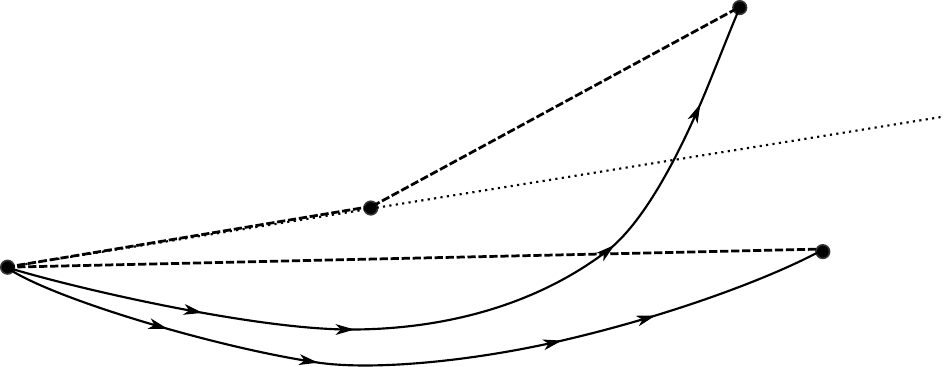}
		\put (80, 40) {\tiny$p_1$}
		\put (86, 8) {\tiny$p_1$}
		\put (0, 7) {\tiny$p_0$}
       \put (15, 0) {\rotatebox{-15}{\tiny$z_a$}}
		\put (60, 13) {\rotatebox{60}{\tiny$z_a$}}
	\end{overpic}
    \hspace{1cm}
    \begin{overpic}[width=0.45\linewidth]{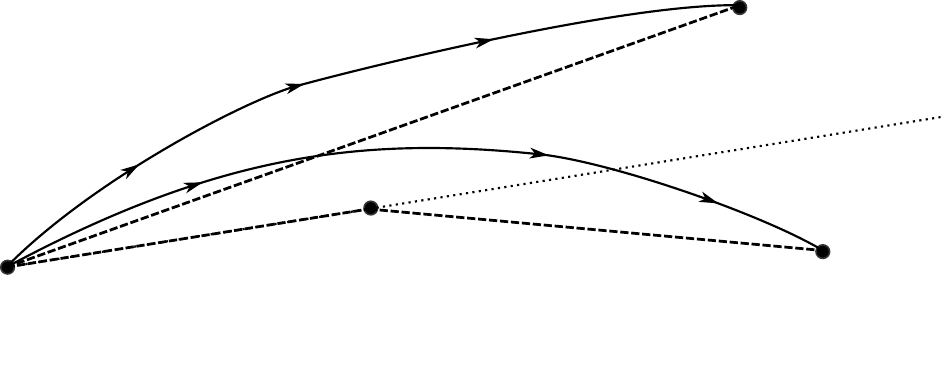}
		\put (0,7) {\tiny$p_0$}
        \put (86,8) {\tiny$p_1$}
		\put (80, 40) {\tiny$p_1$}
		\put (10,22) {\rotatebox{40}{\tiny$z_c$}}	
        \put (50,25) {\rotatebox{-5}{\tiny$z_c$}}	
	\end{overpic}
    \begin{overpic}[width=0.4\linewidth]{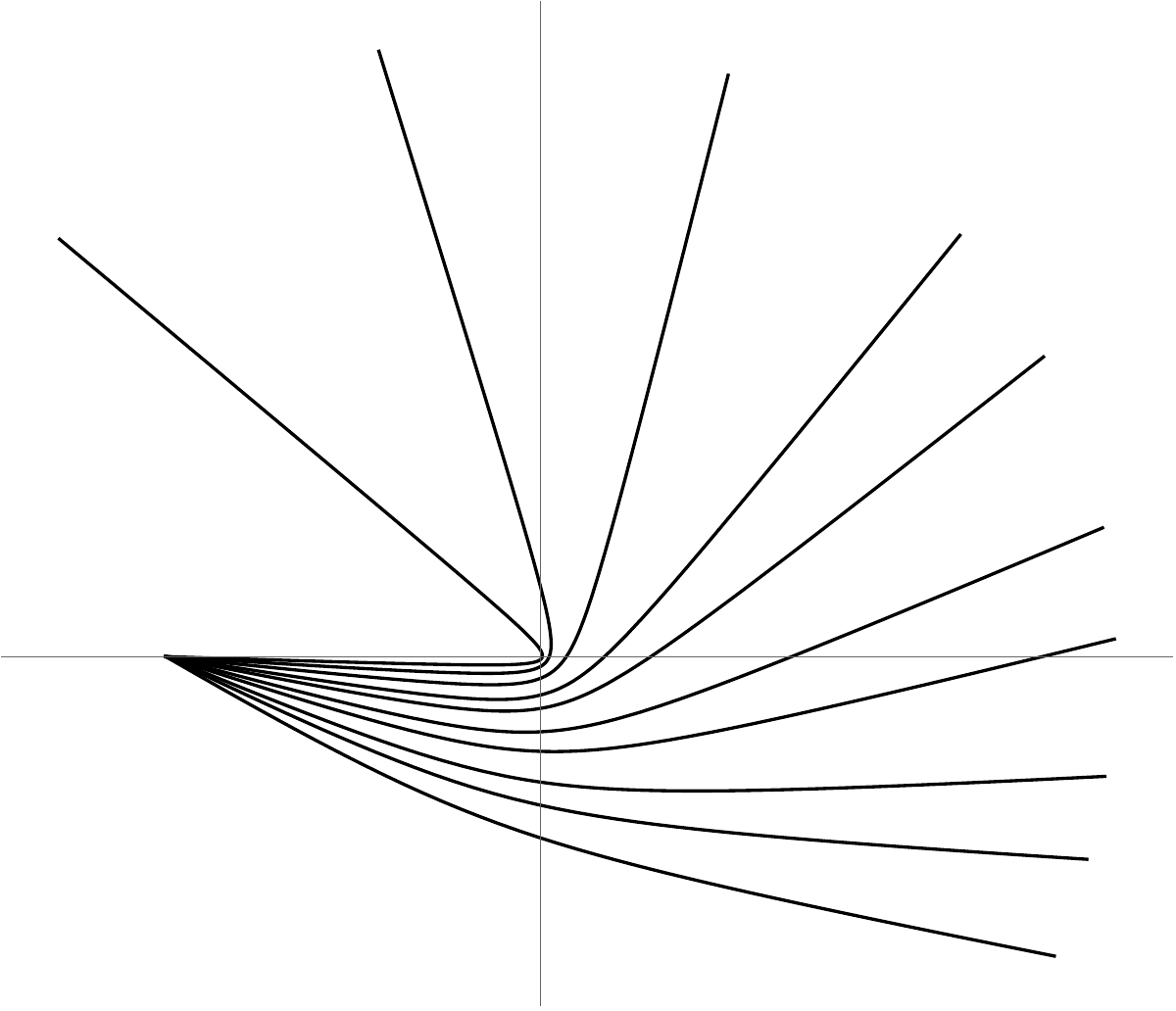}
		\put (8,27) {\tiny$p_0$}
        \end{overpic}
	\caption{In the first line, the anticlockwise, $z_a$, and clockwise, $z_c$, arcs connecting two (possibly antipodal) points. Depending on the difference in the arguments between $p_0$ and $p_1$, they correspond either to  direct or  indirect arcs. In the second line, a family of anticlockwise arcs from $p_0$.}
	\label{fig:antiorario_orario}
\end{figure}

{Let $\gamma:\mathbb{S}^1 \to \partial\Omega$ be a unit speed parametrization of $\partial \Omega$.} The presence of two different Keplerian arcs connecting any suitable pair of points on $\partial \Omega$   suggests that the generating function of the Kepler billiard should be multi-valued. For this reason, we now define the following four determinations.
\begin{definition}\label{def:generatrici}
Given $\delta \in (0,1)$, $\epsilon \in (0,1)$ we consider the sets
\begin{align*}
\mathcal{S}_\delta & = \left\{ (\xi,\eta) \in \mathbb{S}^1 \times \mathbb{S}^1 \;|\; (\gamma(\xi),\gamma(\eta)) \in K_{\delta} \right\}, \\
\mathcal{S}_{\delta,\varepsilon}  & = \left\{ (\xi,\eta) \in \mathcal{S}_\delta  \;|\; |\xi - \eta| \geq \epsilon\right\}, \\
\tilde{\mathcal{S}}_\delta & = {\left\{ (\xi,\eta) \in \mathbb{S}^1 \times \mathbb{S}^1 \;|\; (\gamma(\xi),\gamma(\eta)) \in \tilde K_{\delta} \right\}}.
\end{align*}
For any $h \geq \bar h$ (as in Lemma  \ref{lem: asintotiche_sus}), we define the {\em generating functions} 
\begin{equation*}\label{eq:gen_functions}
\begin{aligned}
	&S_i : \mathcal{S}_\delta \to \R, \quad S_d : \mathcal{S}_{\delta,\epsilon} \to \R, \quad
	S_{i/d} (\xi, \eta;h)=L\left(z_{i/d}(\cdot; \gamma(\xi), \gamma(\eta); h)\right),\\
    &{S_{a/ c}: \tilde{\mathcal S}_\delta\to \R}, \quad S_{a/c} (\xi, \eta;h)=L\left(z_{a/c}(\cdot; \gamma(\xi), \gamma(\eta); h)\right).
\end{aligned}
\end{equation*}
\end{definition}
We stress that, in view of Lemma \ref{lem: asintotiche_sus} and of the regularity of $\gamma$, the generating functions $S_i$ and $S_d$ follow the asymptotics in \eqref{eq:asintotiche}, while $S_{a}$ and $S_c$ follow the estimates in \eqref{eq:Lzazc}. Furthermore, it is well known (see \cite[Appendix A]{IreneSusViNEW}) that, for any index $k \in \{i,d,a,c\}$,
\begin{align*}
 \partial_1 S_{k} (\xi, \eta;h) &= 
 -\frac{\sqrt2}{2} \langle z'_{k} (0;\gamma(\xi), \gamma(\eta);h),
 \dot\gamma(\xi)\rangle \\
 \partial_2 S_{k} (\xi, \eta;h) &= 
 \frac{\sqrt2}{2} \langle z'_{k} (T;\gamma(\xi), \gamma(\eta);h), 
 \dot\gamma(\eta)\rangle.
\end{align*}
where $z'_{k} (0;\gamma(\xi), \gamma(\eta);h)$ and $z'_{k} (T;\gamma(\xi), \gamma(\eta);h)$ denote the initial/final velocities of the Keplerian arcs, either direct or indirect, clockwise or anticlockwise, connecting $\gamma(\xi)$ to $\gamma(\eta)$.
As a consequence, two of them reflect in $\gamma(\eta)$ if and only if
\begin{equation}\label{eq:riflessione}
 \partial_2 S_{k} (\xi, \eta;h) + 
 \partial_1 S_{k'} (\eta, \xi';h) = 0, 
 \quad \forall \xi,\xi' \in \mathbb{S}^1, 
 \quad \forall k,k' \in \{i,d,a,c\}.
\end{equation}

%% file: symbolic_dynamics.tex
The discussion of the previous section shows that every Keplerian orbit connecting two points in $\partial \Omega$ , at a positive energy $h$, converges either to a segment or to a broken line through the center of attraction $c$, as $h \to \infty$.
Moreover, since at every bounce on $\partial \Omega$ trajectories of the Kepler billiard satisfy the reflection law, we will consider pieces of trajectories passing through $c$ for a classical Birkhoff billiard. 
In the high energy regime,  we can \emph{shadow} concatenations of segments and portions of Birkhoff billiard trajectories passing through $c$.  Using this principle, we will consider three different types of symbolic dynamics for Kepler billiards, depending on the geometry of $\partial \Omega$. The first has been constructed in \cite{IreneSusViNEW}, while the second and the third will be introduced in Sections \ref{subsub:nonzero} and \ref{subsub:zero} respectively. 

The first type of symbolic dynamics shadows segments that merge from $c$ and meet $\partial \Omega$ orthogonally. The selected points on $\partial \Omega$ correspond to strict local maxima or minima for the function $|\gamma(\cdot)|$ and are not antipodal to each others. This motivates the following definition.
\begin{definition}\label{def:firstkind}
A point $c\in\Omega$ is said to be {\em of the first kind} if the distance function $|\gamma(\cdot)-c|$ has at least one pair of strict local extremal (minima or maxima) corresponding to non antipodal points of the boundary. Other points will be said of the {\em second kind}.
\end{definition}
For points of the first kind, the following result ensures the existence of a symbolic dynamics.
\begin{lemma}[\cite{IreneSusViNEW}, {Theorem 1.7}] \label{thm:dyn_sym_firstkind}
Suppose that $\partial \Omega$ is analytic and $c\in\Omega$ is a point of the first kind for $\Omega$. Then, the Kepler billiard in $\Omega$ with attraction center at $c$ admits a symbolic dynamics for $h$ sufficiently large. 
\end{lemma}
Thus, in the following, we will restrict our attention to points of the second kind.
It should be noted that, in \cite{BCDeB_Nonlin}, numerical evidence of chaotic behavior for points of the second kind in an elliptic domain is given. Thus, being a point of the first kind does not seem to be a necessary condition for non-integrability. For this reason here we propose two new types of symbolic dynamics obtained by shadowing \emph{Birkhoff triangles}, i.e. classical billiard trajectories coming back to $c$ after two bounces (see Figure \ref{fig:shadowing}).
The notion of focal point will play a central role in our construction.
\begin{figure}
    \begin{overpic}[width=0.4\linewidth]{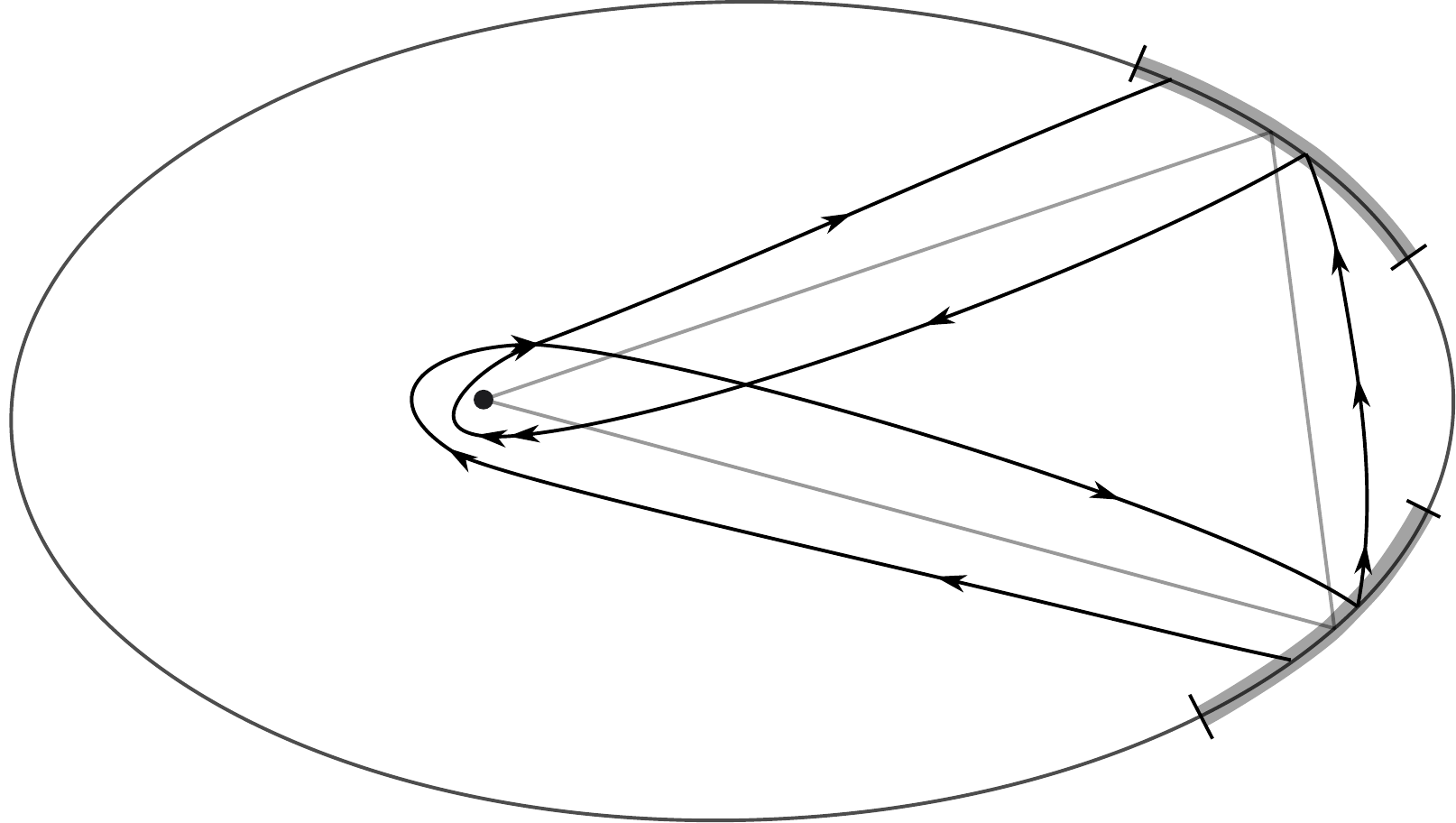}
    \put(24,27){\tiny$c$}
    \put (86, 50) {\tiny$\gamma(\hat\xi)$}
\put (97, 40) {\tiny$\gamma(I)$}
    \put (93, 10) {\tiny$\gamma(\hat\xi')$}
	\put (100, 20) {\tiny$\gamma(I')$}
\end{overpic}
    \caption{Shadowing of a Birkhoff triangle (in light grey): a vertex is at $c$, the other two on $\partial \Omega$ and the bounces are elastic.}
    \label{fig:shadowing}
\end{figure}
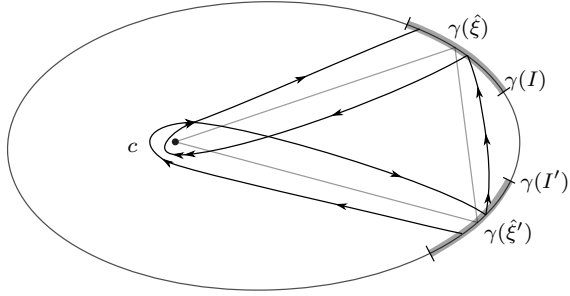
\begin{definition}\label{def:focal}
Let $c\in\Omega$, and let $\varphi_c:{\mathbb S}^1 \to \R$ be defined as follows 
\[
\varphi_c(\xi) = |\gamma(\xi)-c| + |\gamma(\xi)-\gamma(\xi')|+|\gamma(\xi')-c|,
\]
where the segment $c\gamma(\xi)$ elastically reflects into $\gamma(\xi)\gamma(\xi')$.
We say that $c$ is a \emph{focal point} for $\Omega$ if $\varphi_c$ is constant as a function of $\xi$.
\end{definition}
Section \ref{subsec:critical_points} is devoted to the search of non degenerate (in a suitable sense) critical points of $\varphi_c$ which correspond to Birkhoff triangles. They are a key ingredient in our construction, being the stable limit objects to shadow. Section \ref{subsec:per_traj} provides a complete description of the shadowing procedure for periodic trajectories. Finally, in Section \ref{subsec:Miranda} we prove the following theorem and discuss some of its consequences.
\begin{theorem}\label{thm:dyn_sym}
Suppose that $\partial \Omega$ is analytic and $c\in\Omega$ is either of the first kind or is not a focal point for $\Omega$.
Then the Kepler billiard in $\Omega$ with center of attraction at $c$ admits a symbolic dynamics for $h$ sufficiently large. 
\end{theorem}
{Without loss of generality and to simplify the notation we will fix the origin at $c$ and denote the function $\varphi_c$ appearing in Definition \ref{def:focal} by $\varphi$.}

\input{critical_points_new}

\input{periodic}

\input{miranda}

%% file: critical_points_new.tex
\subsection{Existence of punctured Birkhoff triangles}
\label{subsec:critical_points}
In this section we discuss the existence of {two distinct critical points of $\varphi$; all the proofs are given in Appendix \ref{app:proof_critical_points}.\\
Keeping in mind the asymptotic expansions of the Jacobi length given in Lemma \ref{lem: asintotiche_sus} and Lemma \ref{lem:orario_antiorario}, instead of working directly with the function $\varphi$ we will consider
\begin{equation}\label{eq:Psi}
	\Psi: \mathbb{S}^1\times \mathbb{S}^1 \to \mathbb{R}, \quad \Psi(\xi,\eta) = \vert\gamma(\xi)\vert+\vert \gamma (\eta)\vert +\vert \gamma(\eta)-\gamma(\xi)\vert.
\end{equation} 
The relation between the two functions is clear: for any $\xi \in \mathbb{S}^1$, pick the line in the direction $\gamma(\xi)$, reflect it about the normal at $\gamma(\xi)$ and denote by $\xi'=\xi'(\xi)$ the parameter in $\mathbb{S}^1$ defining the intersection point between this line and $\partial \Omega$ (compare with Definition \ref{def:focal}). Then
\begin{equation}
\label{eq:def_phi}
\varphi: \mathbb{S}^1 \to \mathbb{R}, \quad \varphi(\xi) = \Psi(\xi,\xi').
\end{equation}

Our aim is now to show that there are at least two critical points for $\Psi$, which are non degenerate in the sense that they have non-zero topological index (see \cite{Krawcewicz}).
Let us observe that $\Psi$ } is continuous and symmetric on the torus. Furthermore, $\Psi$ is differentiable on the open set 
$$
{\mathcal U} := \mathbb{S}^1\times \mathbb{S}^1 \setminus \Delta, \quad \Delta = \{(\xi,\xi) : \xi \in \mathbb{S}^1\} 
$$
and, by triangular inequality, attains its maximal value
\[
M = \max_{\mathbb{S}^1 \times \mathbb{S}^1} \Psi
\]
outside $\Delta$. Hence $\Psi$ admits, at least, a critical point in $\mathcal{U}$ at level $M$. On the other hand, absolute minima are achieved in $\Delta$, where $\Psi$ is not differentiable, and correspond to points of minimal distance from the origin. \\
The next result analyzes the link between critical points of $\Psi$ and $\varphi$.
\begin{lemma}
\label{lemma:critical_points_psi_phi}
 Let $(\hat\xi,\hat\eta)\in{\mathcal U}$. Then 
\begin{itemize}
\item[(i)] if $(\hat\xi,\hat\eta)$ is a critical point for $\Psi$ then both $\hat \xi$ and $\hat \eta$ are critical points of $\varphi$;
\item[(ii)] if $\varphi$ is constant then all critical points of $\Psi$ are at level $M$ and $\Psi^{-1}(M)$ is connected.
\end{itemize}
\end{lemma} 
{Since we are interested in critical points of $\Psi$ which correspond to distinct trajectories in the phase space, we introduce the following definition.}
\begin{definition}
Two distinct critical points $(\xi,\eta)$ and $(\xi',\eta')$ for the function $\Psi$ are said {\em geometrically distinct} if $(\xi,\eta) \neq (\eta',\xi')$.
\end{definition}

{Finally, the following result describes the structure of the set of critical points of $\Psi$.} 

\begin{proposition}\label{prop:critici_Psi}
   Let us assume that $\partial\Omega$ is analytic. Then, only one of the following holds
   \begin{itemize}
       \item[(i)] the map $\varphi$ is constant,
       \item[(ii)] there exists two isolated geometrically distinct critical points of $\Psi$ both with non-zero index. A maximum point, with index $+1$, and a point with negative index.
   \end{itemize}
\end{proposition}

{Let us remark that Proposition \ref{prop:critici_Psi} might provide \emph{degenerate} triangles, having zero area. In this case Lemma \ref{lem: asintotiche_sus} can not be applied and we have to resort to the generating functions given in Lemma \ref{lem:orario_antiorario}. This motivates the following definitions
\begin{equation}
\label{eq:Psi_a}
\Psi_a: \mathbb S^1 \times \mathbb S^1\to\R \quad \Psi_a(\xi, \eta)= |\gamma(\xi)| + f_a(\gamma(\xi), \gamma(\eta)) + |\gamma(\eta)|. 
\end{equation}
\begin{equation}
\label{eq:Psi_c}
\Psi_c: \mathbb S^1 \times \mathbb S^1\to\R \quad \Psi_c(\xi, \eta)= |\gamma(\xi)| + f_c(\gamma(\xi), \gamma(\eta)) + |\gamma(\eta)|. 
\end{equation}
where $f_a$ and $f_c$ are the ones given in \eqref{eq:def_fa_fc}.

We conclude this paragraph showing that, under some suitable assumptions, non-zero index critical points of $\Psi$ corresponding to zero-area triangles are non-zero index critical points of at least one between $\Psi_a$ and $\Psi_c$.
}
\begin{lemma}
\label{lemma:grado_psi_A}
Let us assume that $(\hat{\xi},\hat{\eta})$ is an isolated critical point of $\Psi$ having non-zero index and that $\gamma(\hat \xi)$ and $\gamma(\hat \eta)$ are antipodal points. Let us further assume that $(\hat{\xi},\hat{\eta})$ is an isolated critical point of $\Psi^*(\xi,\eta) = \vert \gamma(\xi)\vert+\vert \gamma(\eta)\vert$ having zero index. Then, $(\hat{\xi},\hat{\eta})$ is an isolated critical point of $\Psi_a$ and $\Psi_c$ and has non-zero index for at least one of them. 
\end{lemma}

%% file: periodic.tex
\subsection{Existence of periodic Kepler billiard trajectories}\label{subsec:per_traj}

To construct a symbolic dynamics, we need different shadowing strategies depending on whether the critical points determined in the previous section correspond to triangles either with positive or zero area. 
In the first case, we use the generating functions $S_i$ and $S_d$ introduced in Lemma \ref{lem: asintotiche_sus}. In the second, we shall rely on $S_{a}$ and $S_{c}$ defined in Lemma \ref{lem:orario_antiorario}. {In Section \ref{subsec:critical_points} we have shown that, if $\varphi$ is not constant, there are always at least two distinct punctured Birkhoff triangles with non-zero index (Proposition \ref{prop:critici_Psi}). We are thus brought to consider two cases: either there exists a non-zero area Birkhoff triangle or there are at least two zero area ones and all of them are non-zero index critical points of $\Psi$.}

\subsubsection{Periodic trajectories involving a non-zero area triangle}\label{subsub:nonzero}
\begin{figure}
\begin{overpic}[width=0.4\linewidth]{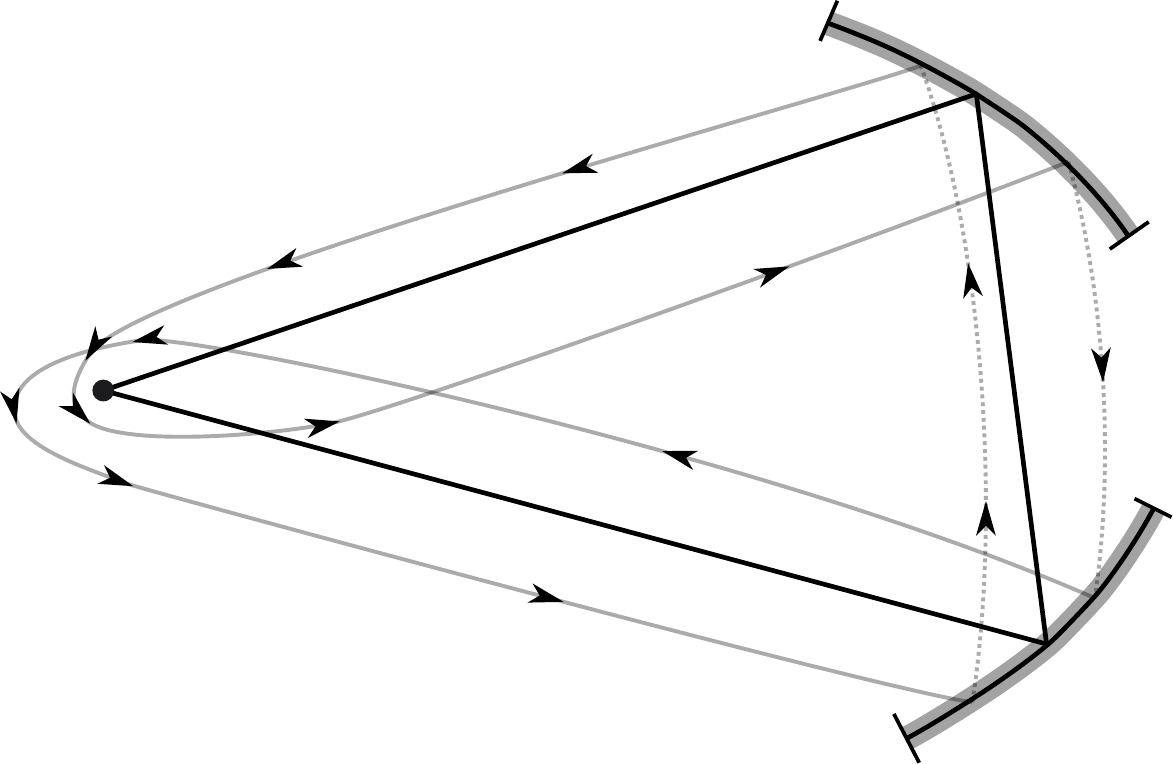} 
\put (87,60) {\tiny$\gamma(I)$}
\put (87,3) {\tiny$\gamma(I')$}
\end{overpic}
\hspace{1cm}
\begin{overpic}[width=0.4\linewidth]{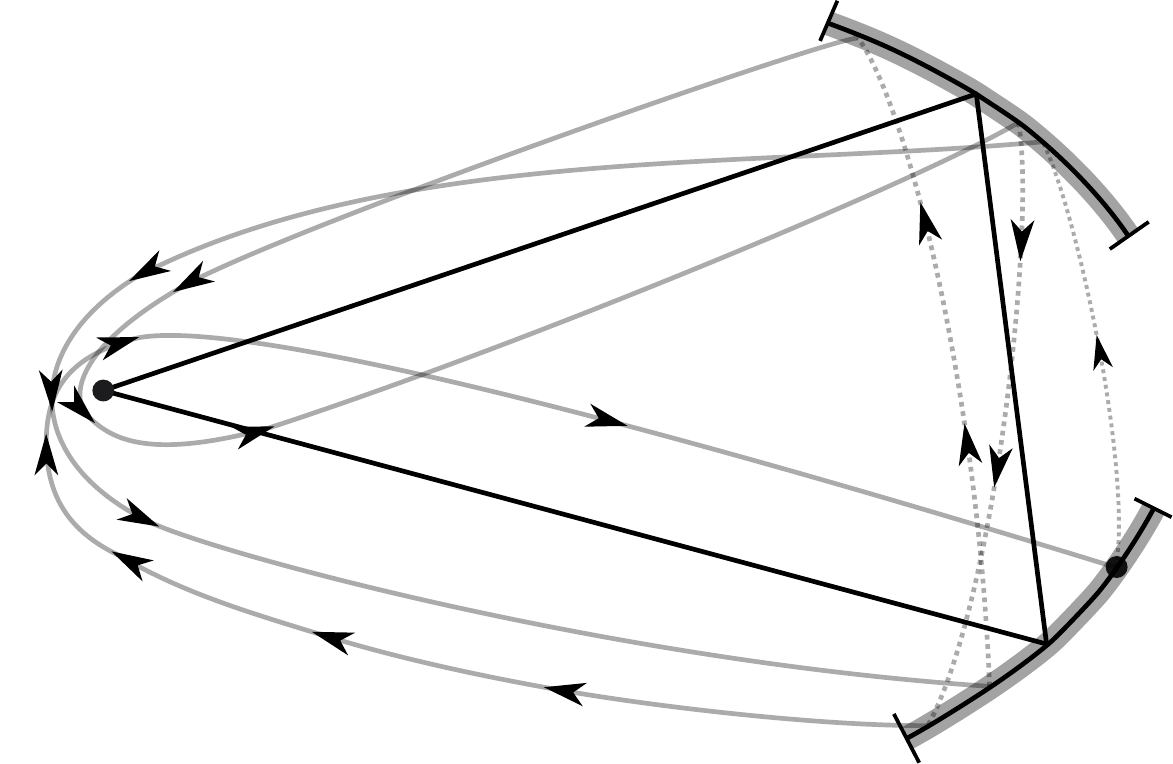}
\put (87,60) {\tiny$\gamma(I)$}
\put (87,3) {\tiny$\gamma(I')$}
\end{overpic}
\vspace{1cm}
\caption{ Possible concatenations of triangular trajectories. Left: a concatenation of arcs realizing  a 2-periodic word with periodicity modulus $TT'$. Right: a concatenation of arcs realizing a periodicity modulus of a 3-periodic world with $s_0=s_1=T'$ and $s_2=T$. As one can observe, every edge of the triangle is shadowed $N$ times, $N$ being the period of $\underline s$.}
\label{fig:4tipi}
\end{figure}

Let us assume that $\partial\Omega$ admits a non-zero area punctured Birkhoff triangle, $\mathcal{T}$, associated with a non-zero index critical point of $\Psi$. {In this section we build periodic Kepler billiard trajectories shadowing $\mathcal{T}$ and its opposite  $\mathcal{T}'$ (i.e. the same triangle but with opposite orientation). Let $\hat \xi, \hat \xi' \in \mathbb{S}^1$ correspond to the vertices of $\mathcal{T}$. They satisfy}  
\begin{equation*}
\langle \gamma(\hat \xi), \dot \gamma (\hat \xi)\rangle = \left\langle\gamma(\hat \xi')-\gamma(\hat \xi), \dot \gamma (\hat \xi)\right\rangle 
\; \text{ and } \;
\langle \gamma(\hat \xi') , \dot \gamma (\hat \xi')\rangle = \left\langle\gamma(\hat \xi)-\gamma(\hat \xi'), \dot \gamma (\hat \xi')\right\rangle.
\end{equation*}
Let then $I$, $I' \subset \mathbb{S}^1$ be disjoint and closed neighborhoods of $\hat \xi$ and $\hat \xi'$ respectively (see Figure \ref{fig:shadowing}) such that no pairs in $\gamma(I)\times\gamma(I')$ are antipodal to each others. Equivalently, recalling Lemma \ref{lem: asintotiche_sus}, 
\begin{equation}\label{eq:intorni}
\left(\gamma(\xi),\gamma(\eta)\right) \in K_\delta \text{ for some $\delta>0$  and for every $(\xi,\eta) \in I\times I'$}.
\end{equation}
According to Lemma \ref{lem: asintotiche_sus} and provided that the energy is sufficiently large, there exist indirect arcs between any pair of points in $\gamma(I\cup I')$ and direct arcs connecting $\gamma(I)$ and $\gamma(I')$, which are completely contained in $\Omega$.\\
Let $\mathcal{S}^\Z = \{T,T'\}^\Z$ be the set of bi-infinite words with letters $T$ and $T'$ and let us start considering its periodic elements. Define \(J_T = I\times I'\) and \(J_{T'} = I'\times I\), choose $\underline{s} \in \mathcal{S}^\Z$ a $N$-periodic word with periodicity modulus $(s_0,\ldots,s_{N-1})$, and consider 
\[
\mathbb{U}_{\underline{s}} = \prod_{i=0}^{N-1} J_{s_i} \text{ whose elements are } u=(u_{0,1},u_{0,2},\ldots,u_{N-1,1},u_{N-1,2}).
\]
Given $u \in \mathbb{U}_{\underline{s}}$, we connect the points parametrized by the elements of $u$ as follows (all indices are intended $\mod N$)
\begin{itemize}
	\item we connect $\gamma(u_{i,1})$ with $\gamma(u_{i,2})$ via a direct arc, $z_d$; 
	\item we connect $\gamma(u_{i,2})$ with $\gamma(u_{i+1,1})$ via an indirect arc, $z_i$.
\end{itemize} 
Thus, given a periodic word in $\underline{s}\in\mathcal{S}^\Z$ and $u \in \mathbb{U}_{\underline{s}}$ we are able to produce a unique concatenation of arcs touching the boundary in the order prescribed by $\underline{s}$. In this case, we say that the concatenation \emph{realizes} the word $\underline{s}$ (see Figure \ref{fig:4tipi}) and we can consider its associated total length
\begin{equation}\label{eq:W_s}
W_{\underline{s}} : \mathbb{U}_{\underline{s}}  \to \R, \quad W_{\underline{s}}(u;h) = \sum_{i=0}^{N-1} \left[ S_{d}(u_{i,1},u_{i,2};h) + S_{i}(u_{i,2},u_{i+1,1};h) \right].
\end{equation}
Nevertheless, this concatenation is generally not a billiard trajectory since the bounces are not elastic. 
In view of \eqref{eq:riflessione} and Lemma \ref{lem: asintotiche_sus} we have the following result.
\begin{lemma}\label{lem:punto_critico}
{Let us assume that $0$ is not a focal point and that there exists at least a non-zero area punctured Birkhoff triangle $\mathcal{T}$. Then, there exists $\bar{h}>0$ such that, for any $h \ge \bar{h}$ and any periodic word $\underline{s} \in \mathcal{S}^\mathbb{Z}$ there exists $\hat{u}\in \mathbb{U}_{\underline{s}}$ such that 
\begin{equation*}
\nabla W_{\underline{s}} (\hat u;h) = 0.
\end{equation*}
In particular, $\hat{u}$ determines a periodic billiard trajectory at energy $h$ and,} for $k=1,2$, $i=0, \dots, N-1$, 
	\[
	\frac{\partial}{\partial u_{i,k}} W_{\underline{s}} (u;h) = \sqrt{h} \frac{\partial}{\partial u_{i,k}} 
\Psi(u_{i,1},u_{i,2}) 
+O\left(\frac{\mu}{\sqrt{h}}\right).
\]
Where the remainder depends only on $u_{i-1, j}, u_{i, j}, u_{i+1, j}$, $j=1,2$. 
\end{lemma}
\begin{proof}  
Let $\hat \xi, \hat \xi' \in \mathbb{S}^1$ denote the parameters of  the vertices of $\mathcal{T}$. 
From the above discussion $(\hat\xi,\hat\xi')$ and $(\hat\xi',\hat\xi)$ are non-zero index critical points of $\Psi$.
Let us consider the function
\[
\Psi_{\underline{s}}: \mathbb{U}_{\underline{s}} \to \mathbb R, \qquad
\Psi_{\underline{s}}(u) = \sum_{i=0}^{N-1}\Psi(u_{i,1},u_{i,2}).
\]
By construction and thanks to Lemma \ref{lem: asintotiche_sus} we have that $
\nabla W_{\underline{s}} (u;h) = \sqrt{h} \nabla  \Psi_{\underline{s}}(u) + O\left({\mu}/{\sqrt{h}}\right)$.
Furthermore, the point $u^* \in \mathbb{U}_{\underline{s}}$, defined as $(u_{i,1},u_{i,2}) = (\hat \xi,\hat{\xi}')$ if $s_i = T$ and  $(u_{i,1},u_{i,2}) = (\hat \xi',\hat{\xi})$ if $s_i = T'$, is the only critical point of $\Psi_{\underline{s}}$ in $\mathbb{U}_{\underline{s}}$ (possibly reducing the neighborhoods $I$ and $I'$).
The product property of the index (see \cite{Krawcewicz}) ensures that the index of $\nabla  \Psi_{\underline{s}}$ at $u^*$ is different from zero.
By definition, since there are no other critical points in $\mathbb{U}_{\underline{s}}$ a part from $u^*$, the topological degree of $\nabla \Psi_{\underline{s}}\vert_{\partial \mathbb{U}_{\underline{s}}}$ coincides with the index of $\nabla  \Psi_{\underline{s}}$ at $u^*$. Furthermore, the homotopy property of the topological degree ensures the existence of a nonzero index critical point  of $W_{\underline s}$ in $\mathbb{U}_{\underline{s}}$ provided $h\ge \bar{h}$ for some large $\bar{h}$. Finally, let us observe that this bound on the energy is uniform in $s \in \mathcal{S}^{\mathbb{Z}}$ since the partial derivatives of $\partial_{ u_{i,k} }W_{\underline{s}} (u;h)$ depend only on $u_{i-1, j}, u_{i, j}, u_{i+1, j}$, for $j=1,2$ and are close to the ones of $W_{\underline{s}}$, uniformly in $\underline{s}$, in the $C^0$ topology.
\end{proof}

\subsubsection{Periodic trajectories involving zero area triangles}\label{subsub:zero}
\begin{figure}
\vspace{1cm}
\begin{overpic}[width=0.7\linewidth]{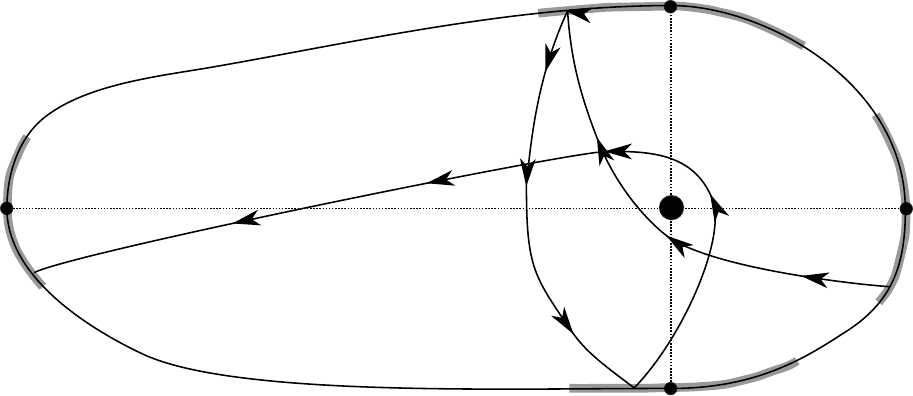}
\put(-10, 20){\tiny{$\gamma(\xi_M)$}}
\put(101, 20){\tiny{$\gamma(\xi_m)$}}
\put(70, 45){\tiny{$\gamma(\xi_P)$}}
\put(70, -2){\tiny{$\gamma(\xi_Q)$}}
\end{overpic}
    \caption{Degenerate case: the minimum and maximum of $|\gamma(\cdot)|$ correspond to antipodal points, and there is a second pair of antipodal homothetic directions corresponding to inflection points for $|\gamma(\cdot)|$. The piecewise curve with arrows is an example of transfer concatenation from $I_m$ to $I_M$ through an anticlockwise arc. }
    \label{fig:caso_degenere}
\end{figure}
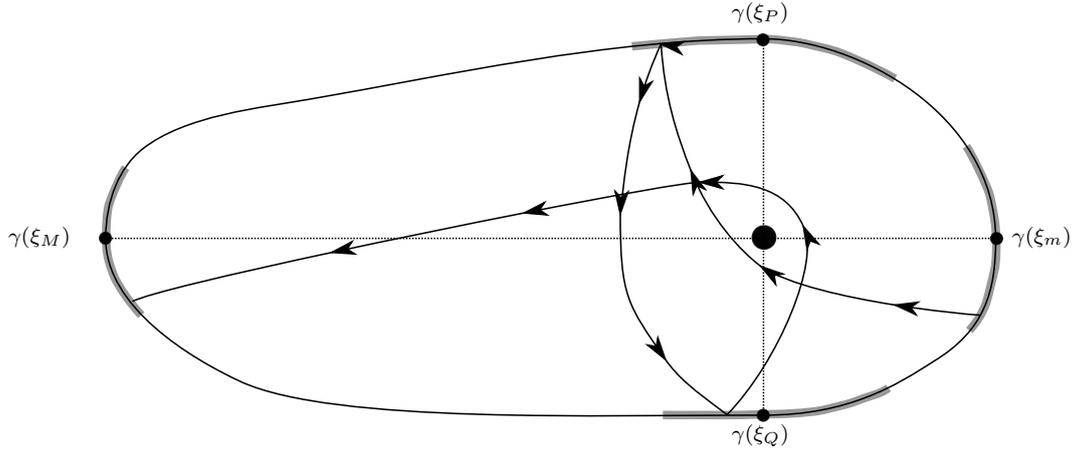

Let us now suppose that any non-zero index critical point of $\Psi$ corresponds to a zero area punctured Birkhoff triangle. These \emph{degenerate} triangles are two-periodic Birkhoff billiard trajectories passing through the origin.\\
In light of Lemma \ref{thm:dyn_sym_firstkind}, the only case left to be checked is when the local minimum $\xi_m$ and maximum $\xi_M$ are unique and antipodal. 
Since we are assuming that all (non-zero index) punctured  Birkhoff triangles have zero area, the function $|\gamma(\cdot)|$ admits a second pair of critical points $\xi_P, \xi_Q$, which are inflection points and correspond to antipodal directions (see Figure \ref{fig:caso_degenere}). Thanks to Lemma \ref{lemma:grado_psi_A} we know then that the chord $\gamma(\xi_P)\gamma(\xi_Q)$ determines a non-zero index critical point for at least one between $\Psi_a$ and $\Psi_c$. Without loss of generality we can assume that the index is non-zero for $\Psi_a$.

We build now periodic trajectories shadowing the segments $0\gamma(\xi_m)$ and $0\gamma(\xi_M)$, using the zero-area triangle determined by $\gamma(\xi_P)\gamma(\xi_Q)$ as a transfer orbit.
Let us then consider four disjoint compact neighborhoods, namely, $I_m, I_M, I_P, I_Q$, respectively of $\xi_m, \xi_M, \xi_P, \xi_Q$, such that 
\[
\text{if } (\xi, \eta)\in (I_{m/M}\times I_P)\cup(I_Q\times I_{m/M}) \text{ then } (\gamma(\xi),\gamma(\eta))\in  K_\delta \ \text{for some }\delta\in(0,1). 
\]
According to Lemma \ref{lem:orario_antiorario}, for every pair of points  in $\gamma(I_P)\times \gamma(I_Q)$ there are exactly two arcs connecting them: the first, $z_a$, rotating in anticlockwise sense and the second, $z_c$, in the opposite direction. 
Let us consider now $\widetilde {\mathcal S}^\Z = \{m, M\}^\Z$, let $\underline s\in \widetilde{\mathcal S}^\Z$ having periodicity modulus $\underline s=(s_0, \dots, s_{N-1})$ and define then 
\begin{equation*}
\begin{aligned}
 J_{m/M} &= I_{m/M}\times I_P\times I_Q, \quad \mathbb U_{\underline s}=\prod_{i=0}^{N-1} J_{s_i}\\
    \widetilde W_{\underline s}: \mathbb U_{\underline s}\to \R, \quad \widetilde 
 W_{\underline s}(u; h) &= \sum_{i=0}^{N-1}\left[S_i(u_{i, 1}, u_{i,2}; h)+ S_a(u_{i,2}, u_{i,3}; h)+ S_i(u_{i,3}, u_{i+1,1}; h)\right].
    \end{aligned}
\end{equation*}
We are connecting  points of $\gamma(I_{m/M})$ with points of $\gamma(I_{P/Q})$ trough  indirect arcs, and points of $\gamma(I_P)$ with points of $\gamma(I_Q)$ with anticlockwise ones (see Figure \ref{fig:caso_degenere}). Of course, one can consider an analogous construction by taking $S_c$, namely, connecting $\gamma(I_P)$ to $\gamma(I_Q)$ with $z_c$.\\ {With the same reasoning, we can prove a result analogous to Lemma \ref{lem:punto_critico}.}
\begin{lemma}\label{lem:punto_critico_degenerate}
Let us assume that $0$ is not a focal point of the second kind and that any non-zero index critical point of $\Psi$ corresponds to a zero area punctured Birkhoff triangle.
Then, there exists $\bar{h}>0$ such that, for any $h \ge \bar{h}$ and any periodic word $\underline{s} \in \widetilde{\mathcal{S}}^\mathbb{Z}$ there exists $\hat{u}\in \mathbb{U}_{\underline{s}}$ such that 
\begin{equation}\label{eq:critico_rif}
	\nabla \widetilde{W}_{\underline{s}} (\hat u;h) = 0.
\end{equation}
In particular, $\hat{u}$ determines a periodic billiard trajectory at energy $h$ and, for $k=1,2$, $i=0, \dots, N-1$, 
\begin{equation}\label{eq:gradWTilda}
\frac{\partial}{\partial u_{i, k}}\widetilde  W_{\underline s}(u) = \sqrt{h}\frac{\partial}{   \partial u_{i, k}} \big(2|\gamma(u_{i, 1})|+ 
\Psi_a(u_{i, 2},u_{i, 3})
\big) + O\left(1\right), 
\end{equation}
Where the remainder depends only on $u_{i-1, j}, u_{i, j}, u_{i+1, j}$, $j=1,2$. 
\end{lemma}

%% file: miranda.tex
\subsection{The symbolic dynamics}\label{subsec:Miranda}

We are now in the position to prove Theorem \ref{thm:dyn_sym}. 
\begin{proof}[Proof of Theorem \ref{thm:dyn_sym}]
When $0$ is a point of the first kind then  Lemma \ref{thm:dyn_sym_firstkind} applies, concluding the proof.
Let us then suppose that $0$ is a point of the second kind and not focal. \\
In view of Proposition \ref{prop:critici_Psi}, the function $\Psi$ admits at least two critical points with non-zero index. If one of them corresponds to a non-zero area triangle, we are in the case covered by Lemma \ref{lem:punto_critico}. If all of them correspond to zero area triangles, we are in the case of Lemma \ref{lem:punto_critico_degenerate}.
Thus, in the hypotheses of the theorem, we can realize every periodic sequence in $\mathcal{S}^\mathbb Z$ or  $\widetilde{\mathcal{S}}^\mathbb Z$. Using the notation introduced in Definition \ref{def:sym_dyn}, we have proved that the projection map $\pi$ is well defined on initial conditions of periodic trajectories realizing (periodic) sequences in  $\mathcal{S}^\mathbb Z$ or   $\widetilde{\mathcal{S}}^\mathbb Z$.

The map $\pi$ is then extended to a map which is onto on the whole $\mathcal{S}^\mathbb Z$ (resp.  $\widetilde{\mathcal{S}}^\mathbb Z$)
following the limit procedure detailed in  \cite{IreneSusViNEW}. 
For the reader's convenience, we give here a sketch. First of all, we define the set $X$ of initial conditions generating trajectories with bounces inside the intervals considered in the proofs of Lemma \ref{lem:punto_critico} or \ref{lem:punto_critico_degenerate}, for all forward and backward iterations of the billiard map. This set is non-empty since it contains at least the initial conditions of all periodic trajectories built so far. Furthermore, the set $X$ is invariant and the map $\pi: X \to \mathcal S^\Z$ (resp. $\widetilde{\mathcal S}^\Z$) is well defined. 
    Then, one shows that this map is  surjective: this is based on a diagonalization process, taking suitable limits of periodic billiard trajectories. Finally, continuity follows from continuous dependence on initial conditions and uniform boundedness of the time needed to connect any two couple of points with a chosen type of arc. 
\end{proof}

The existence of a symbolic dynamics, although not sufficient to guarantee that our system admits a topologically chaotic subsystem in the sense of Devaney (see \cite{Dev_book}), implies some kind of \emph{complexity}, which can be expressed in terms of non-existence of analytical first integral of the motion, as well as of the positivity of the topological entropy associated to a subsystem. The proof of the next Corollary follows from classical results in  \cite{MR1995704} and, for the second statement, adapting a classical argument by Katok (see \cite{koz1983} and \cite[Theorem 4.5]{IreneSusViNEW}). 

\begin{corollary}\label{cor:katok}
Suppose that $\partial \Omega$ is analytic, and that $0$ is not a focal point of the second kind for $\Omega$.
Then the following fact hold: 
\begin{itemize}
    \item the Kepler billiard in $\Omega$ with attraction center at $0$ has positive topological entropy; 
    \item the first return map associated to such Kepler billiard in $\Omega$ does not admit analytic first integrals of the motion.  
\end{itemize}
\end{corollary}

Theorem \ref{thm:dyn_sym} gives a geometric condition that guarantees the presence of a symbolic dynamics for the  Kepler billiard inside an analytic domain if the energy is high enough. We aim now to characterize the boundaries that {\it do not} satisfy this condition, namely, boundaries that admit focal points of the second kind.

%% file: rigidity.tex
In this section, we focus on counting the number of possible focal points of the second kind (according to Definitions \ref{def:firstkind} and  \ref{def:focal}) in a strictly convex analytic domain. 
We recall that saying that $c$ is focal means that the function $\varphi_c$ is constant; for example, when $\partial\Omega$ is an ellipse and $c$ is located at one of the two foci,  the value of $\varphi_c$ equals the double of the major axis. 

The main result of this section is the following.
\begin{theorem}
 \label{thm:at_most_one_point}
For any $\partial \Omega$ real-analytic there are at most two focal points of the second kind. If there are exactly two then $\partial \Omega$ is a non-circular ellipse and the focal points are its foci.
\end{theorem}

As a consequence, under some additional symmetry constraints on $\partial \Omega$, in the non-elliptic case there are no focal points of the second type. 

\begin{corollary}
\label{cor:focal_points_with_symmetry}
    Assume that $\partial \Omega$ is not a circle and that it has a cyclic symmetry of order $p$, that is if it is invariant by a rotation of angle $2\pi/p$. Then it has no focal points of the second kind unless $p=2$ and $\partial\Omega$ an ellipse.
\end{corollary}

In this section, we are going to deal specifically with some properties of the Birkhoff billiard associated with $\partial\Omega$, which is understood to be real-analytic, when not specified. No Kepler potential will be involved. The terms \textit{trajectory, billiard, billiard map} etc. should be interpreted in this sense. In order to prove Theorem \ref{thm:at_most_one_point} we need to make a preliminary remark that follows immediately from Definition \ref{def:firstkind} and to prove a series of lemmata.

\begin{figure}[t]
\begin{overpic}[width=0.40\linewidth]{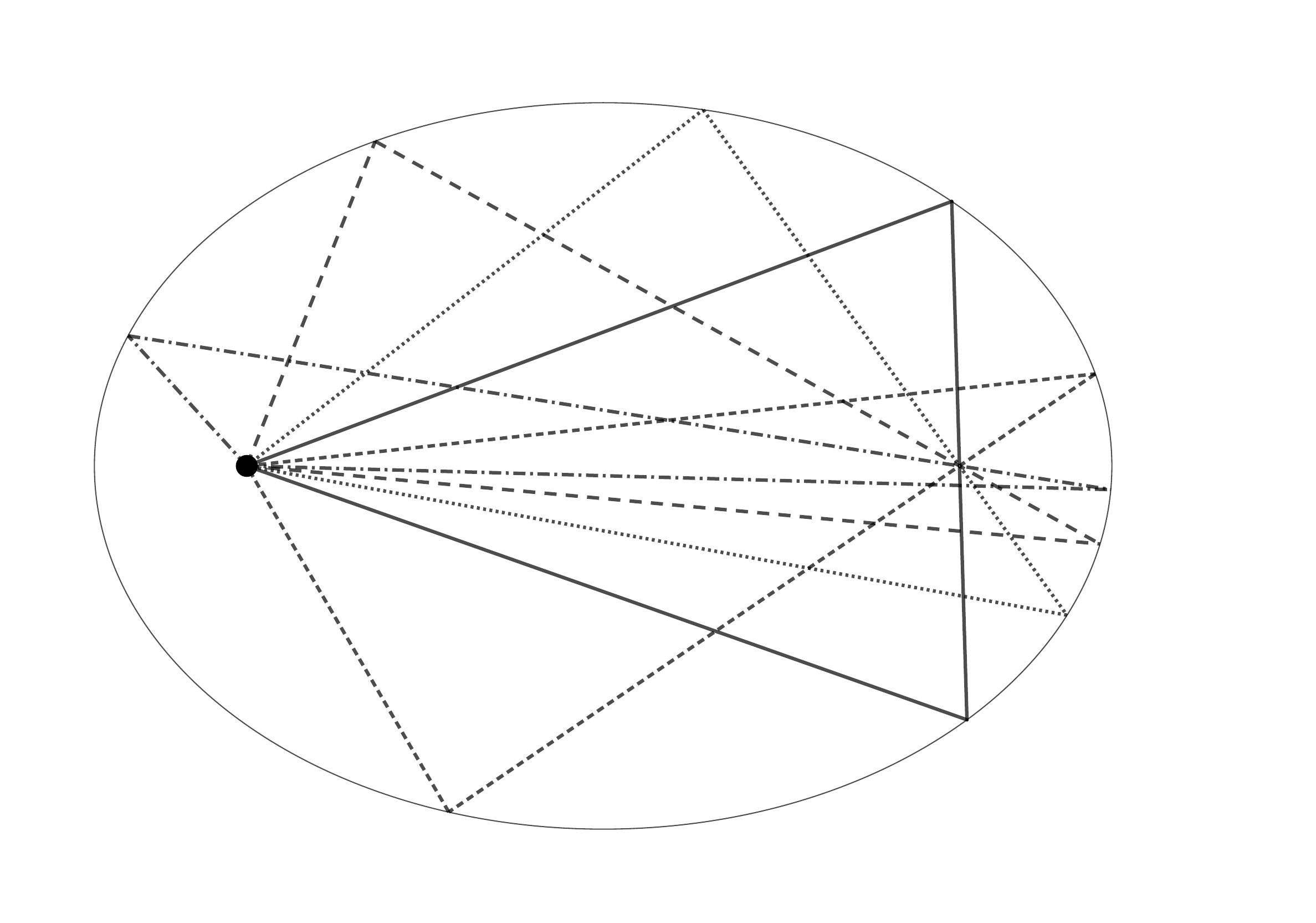}
\end{overpic}
\caption{If $0$ is focal, then all triangles with a vertex at the origin and one reflection on the boundary satisfy the elastic reflection condition also at the second bouncing point. These triangles have all the same perimeter, the value of $\varphi$.}
\label{fig:triangoli}
\end{figure}

\begin{remark}\label{rem:secondkind}
    When $0$ is a point of the second kind and $\partial\Omega$ is not a circle, then the distance function \(|\gamma(\cdot)|\) has exactly two strict local and global extremals (minimum and maximum) corresponding to antipodal points. 
    There might be other inflection-type critical points. Pairs of antipodal critical points give rise to two-periodic trajectories. 
\end{remark}

\begin{lemma}\label{lem:geometric}
If $0$ is a focal point for $\Omega$, then all rays from the origin are pieces of billiard trajectories that return to the origin after exactly two bounces (see Figure \ref{fig:triangoli}) (except for the chords).
Moreover, if $\Omega$ admits two distinct focal points then they must lie on the same chord orthogonal to the boundary and the corresponding values of the functions $\varphi$ are equal.    
\end{lemma}
\begin{proof}
Being $0$ a focal point, $\varphi$ is constant and Lemma \ref{lemma:critical_points_psi_phi} guarantees that all triangles with a vertex at the origin and one reflection on $\partial \Omega$ are actually critical points of $\Psi$. Thus, they satisfy the elastic reflection condition also at the second bouncing point, and they are pieces of billiard trajectories. As for the second statement, without loss of generality, let us suppose that the two focal points are 0 and $c \neq 0$. Since the functions $\varphi$ and $\varphi_c$ are constant, the billiard trajectory that passes through 0 and $c$ after two elastic bounces must go back to 0 and $c$. Being $c \neq 0$, this implies that the first bounce, and hence the second, is orthogonal.
We deduce that 0 and $c$ lie on the same orthogonal chord and $\varphi$ and $\varphi_c$ are constantly equal to the double of the length of such chord.
\end{proof}

\begin{lemma}
\label{lemma:at_most_2_bad _points}
Let $0$ and $c\neq 0$ be two focal points of the second kind, and let $p_1=\gamma(\xi_1)$ and $p_2=\gamma(\xi_2)$ be the end-points of their shared orthogonal chord.
Then $\xi_1$ and $\xi_2$ are non-degenerate critical points for both 
distance functions from 0 or $c$ to $\partial \Omega$.
Moreover, $\Omega$ can not contain more than two focal points of the second kind.
\end{lemma}

\begin{proof}
For $i=1,2$, let $k_i$ be the curvatures at $p_i$, $d_i$ the distance between $0$ and $p_i$ and $d = d_1+d_2$ the length of the chord. 
Since $\varphi$ is constant, by Lemma \ref{lemma:critical_points_psi_phi} (ii), the Hessian of $\Psi$ is degenerate at $(\xi_1,\xi_2)$.
A straightforward computation shows that the Hessian of $\Psi$ and its determinant at $(\xi_1,\xi_2)$ are
\begin{equation}
\label{eq:det_hess_Psi}
\begin{aligned}
\mathop{Hess} &= \begin{pmatrix}
\frac{1}{d}+\frac{1}{d_1}-2 k_1 & \frac{1}{d}\\
\frac{1}{d} &\frac{1}{d}+\frac{1}{d_2}-2 k_2
\end{pmatrix}\\ 
\det\mathop{Hess} & = \left( \frac{1}{d}+\frac{1}{d_1}-2 k_1 \right)\left( \frac{1}{d}+\frac{1}{d_2}-2 k_2 \right)-\frac{1}{d^2}.
\end{aligned}
\end{equation}
Since also $\varphi_c$ is constant,
the Hessian of the corresponding function $\Psi_c$
\[
\Psi_c(\xi,\eta)=\vert \gamma(\xi)-c\vert+ \vert \gamma(\eta)-c\vert+\vert \gamma(\eta)-\gamma(\xi)\vert
\]
is degenerate as well, and has essentially the same form as $\mathop{Hess}$, except for the fact that the distances $d_i$ are computed with respect to $c$. Thus, all focal points on $p_1p_2$ must satisfy the equation $\det \mathop{Hess}=0$, that is
\begin{equation}
\label{eq:equation_bad_points}
\left( \frac{1}{d}+\frac{1}{t}-2 k_1 \right)\left( \frac{1}{d}+\frac{1}{d-t}-2 k_2 \right)=\frac{1}{d^2}, \text{ for some } t \in (0,d).
\end{equation}
Without loss of generality, since $0$ is a point of the second kind, we can assume that $p_1$ is either the minimum of the distance from $0$ or an inflection point and $p_2$ is a maximum or an inflection point (Remark \ref{rem:secondkind}) and thus
\[
\frac{d^2}{d \xi^2}\Big\vert_{\xi=\xi_1}  \vert \gamma(\xi) \vert = \frac{1}{d_1}-k_1 \ge 0, \quad \frac{d^2}{d \xi^2}\Big\vert_{\xi=\xi_2} \vert \gamma(\xi) \vert= \frac{1}{d_2}-k_2 \le 0. 
\] 
Since $d>d_2$ we immediately obtain 
\[
 \frac{1}{d}+\frac{1}{d_2}-2 k_2 <0,
\] 
and, comparing with \eqref{eq:equation_bad_points}, we have
\[
\frac{1}{d}+\frac{1}{d_1}-2 k_1 <0.
\] 
and in particular $\frac{1}{d}<k_1$ and $\frac{1}{d}<k_2$. Thus, a straightforward computation shows that solutions of \eqref{eq:equation_bad_points} correspond to the roots in $(0,d)$ of the following polynomial
\[
P(t) = 1 - 2 k_1 t + \frac{k_1 + k_2 - 2 d k_1 k_2}{d (1 - d k_2)} t^2.
\]
Since $P(0) = 1$, $P(d)= \frac{1-d k_1}{1- d k_2}>0$ and $\frac{k_1 + k_2 - 2 d k_1 k_2}{d (1 - d k_2)}>0$,
$P$ can have either two solutions in $(0,d)$ or one, when its discriminant is zero. 
The discriminant of $P$ reads
\[
\Delta = (-1 + d k_2) (-k_2 - k_1 + d k_1 k_2),
\]
and it can be zero if and only if $d = (k_1+k_2)/k_1k_2$. 
Since we have assumed the existence of two distinct roots of $P$, then  $d > 1/k_1+1/k_2$; as a consequence $\xi_1$ and $\xi_2$ can not be both inflection points, hence (by Remark \ref{rem:secondkind}) they are necessarily global extremals of $|\gamma(\cdot)|$ and $|\gamma(\cdot)-c|$.\\
Next, we show that $\xi_1$ and $\xi_2$ are non-degenerate critical points. Assume by contradiction that $\xi_1$ is degenerate, i.e. that $d_1 = 1/k_1$. Expanding $P(1/k_1) = 0$ we find that
\[
(-1 + d k_1) (-k_2 - k_1 + d k_1k_2) =0
\] 
and thus also $d_2 = 1/k_2$. Thus $\Delta =0$, contradicting the existence of two distinct focal points.\\
Assume now that $c' \neq 0$ is a third focal point of second kind. Then the end-points of the chord containing 0 and $c'$ are the global extremals for $|\gamma(\cdot)|$. Then necessarily $c'=c$ since this chord contains at most two focal points of second kind.
\end{proof}

To simplify the notation, when we deal with two distinct focal points we will assume they are $0$ and $c\neq 0$. Furthermore, the endpoints of their shared orthogonal chord will be called $p_1=\gamma(\xi_1)$ and $p_2=\gamma(\xi_2)$.
Now we consider the effect of the presence of a focal point on the phase portrait associated with the Birkhoff first return map or possibly with its second iterate.

\begin{lemma}\label{lem:invariant_curve}
If $0$ is a focal point for $\Omega$, then to any 2-periodic trajectory through the origin with end-points $\gamma(\xi_1)$ and $\gamma(\xi_2)$ we can associate an invariant curve in the phase-space cylinder $\mathbb{S}^1 \times (0,\pi)$ of the billiard map $T$ as follows
\begin{equation}
\label{eq:def_invariant_curves}
\delta(\xi) = 
\begin{cases}
\left(\xi,\arccos\left(-\langle\frac{\gamma(\xi)}{\vert \gamma(\xi)\vert},\dot\gamma(\xi)\rangle\right)\right), \quad \text{if }\,\xi \in (\xi_1,\xi_2], \\[1ex]
\left(\xi,\arccos\left(\langle \frac{\gamma(\xi)}{\vert \gamma(\xi)\vert},\dot\gamma(\xi)\rangle\right)\right), \quad \text{if }\,\xi \in (\xi_2,\xi_1+L].
\end{cases}
\end{equation}
\end{lemma}
\begin{proof}
Let us consider a $2-$periodic trajectory passing through the origin with bouncing points $\gamma(\xi_1)$ and $\gamma(\xi_2)$. This segment divides $\Omega$ into two regions, $\Omega^\pm$, and its boundary, respectively, into two arcs $\partial\Omega^\pm$.
We define now a piecewise analytic curve $\hat \delta:\mathbb{S}^1 \to \partial\Omega \times \mathbb{R}^2$ (see Figure \ref{fig:riccetto}, left) as
\begin{equation}
\label{eq:invariant_curve_lines}
\hat \delta(\xi) = 
\begin{cases}
\left(\gamma(\xi),-\frac{\gamma(\xi)}{\vert \gamma(\xi)\vert}\right), \quad \text{if }\,\xi \in (\xi_1,\xi_2], \\[1ex]
\left(\gamma(\xi),R_\xi\frac{\gamma(\xi)}{\vert \gamma(\xi)\vert}\right), \quad \text{if }\,\xi \in (\xi_2,\xi_1+L], 
\end{cases}
\end{equation}
where $R_\xi$ denotes the reflection about the tangent line at the point $\gamma(\xi)$. 
The invariance of $\hat \delta$ follows from the fact that no triangle of our family is contained in a single region bounded by the $2-$periodic trajectory, as we now prove. Assume by contradiction that there is a triangle that it is contained in, say, $\Omega^+$. Label the vertices on $\partial \Omega^+$ clock-wise as $y$ and $z$ as in Figure \ref{fig:riccetto}, right. Moving $y$ clockwise, by continuity, there is a $y_0$ for which the corresponding vertex $z_0$ either touches the bouncing point of the $2-$periodic trajectory or coincides with $y_0$. The first case is not possible since $y_0$ should then coincide with the other end-point of the periodic trajectory. The second is not possible either since the angle between the segment $y-x$ and the tangent line is always bounded from below. 
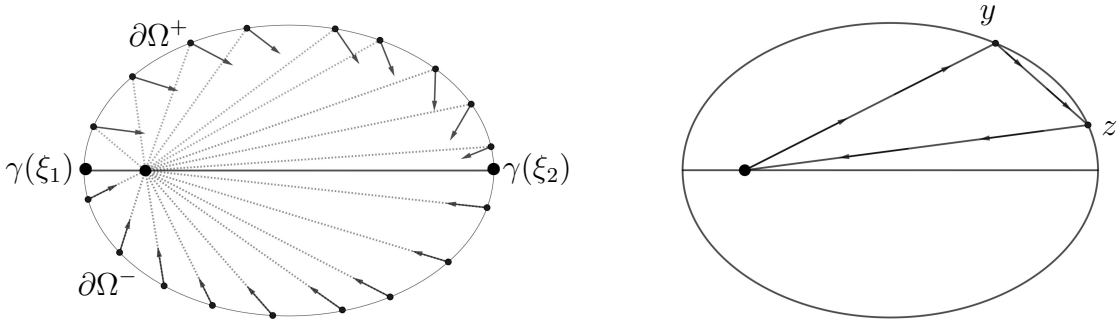
\begin{figure}[ht]
\begin{overpic}[height=4cm]{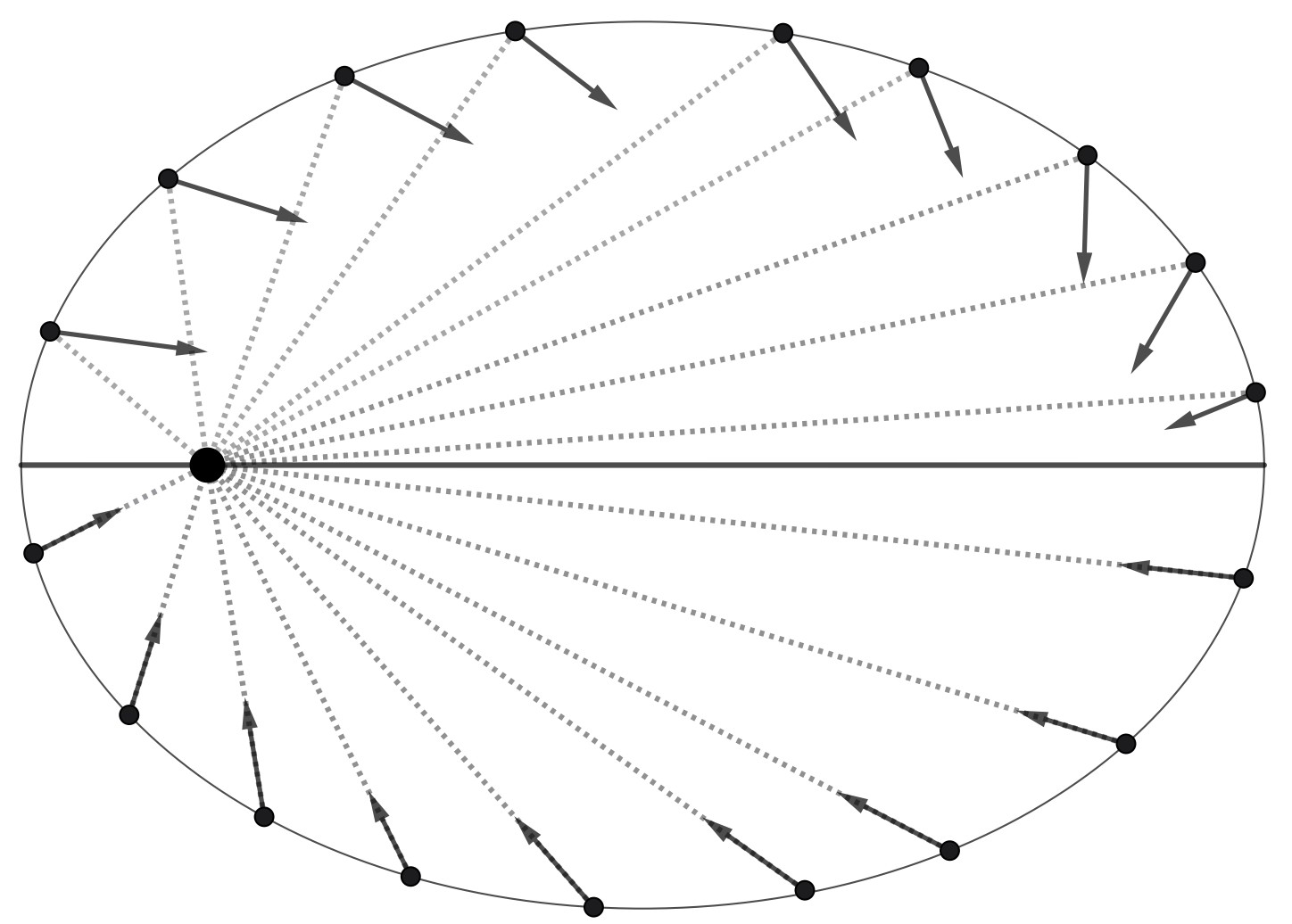}
\put(12,65){$\partial \Omega^+$}
\put(0,6){$\partial \Omega^-$}
\put(-17,34){$\gamma(\xi_1)$}
\put(100,34){$\gamma(\xi_2)$}
\put(96,34){$\bullet$}\put(0,34){$\bullet$}
\end{overpic}     
\hspace{2cm}
\begin{overpic}[height=4cm]{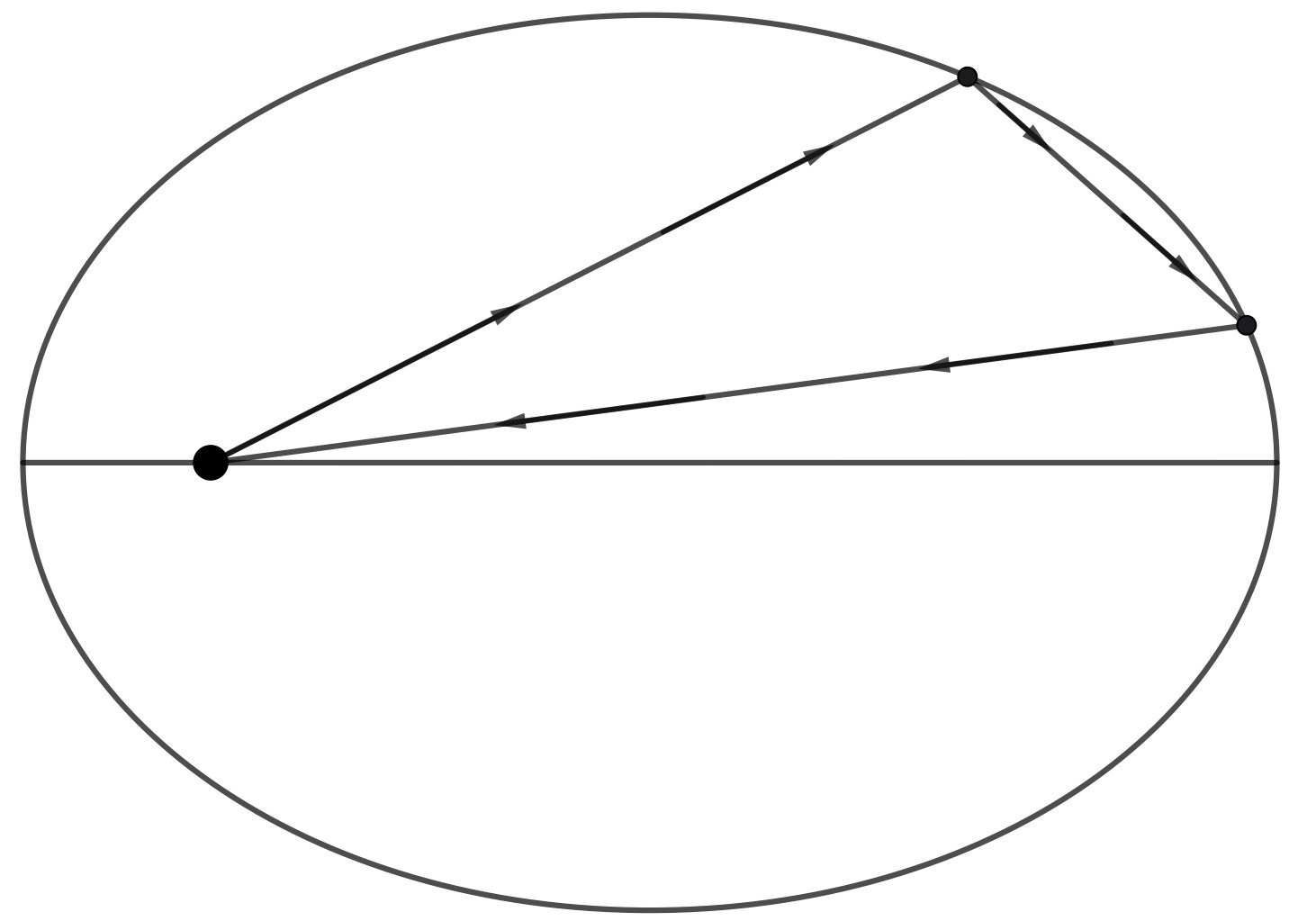}
\put(70,70){$y$}
\put(98,43){$z$}
\end{overpic} 
\caption{Left, the vector field $\hat \delta$ introduced in the proof of Lemma \ref{lem:invariant_curve}. Right, following the proof of Lemma \ref{lem:invariant_curve}, every double-bouncing triangle has a vertex on $\partial \Omega^+$ and the second on $\partial \Omega^-$. }
\label{fig:riccetto}
\end{figure}
\noindent Hence, $\hat\delta$ corresponds to the homotopically non-trivial curve in the phase cylinder, $\delta$, which is invariant and piece-wise analytic.
\end{proof}

A direct consequence of Lemma \ref{lem:invariant_curve} is the next result.

\begin{lemma}\label{lem:graphs}
If $0$ is focal, then there are two symmetric analytic invariant graphs in the phase space for $T^2$, which are given by
\begin{equation*}
\delta_\pm(\xi) = 
\left(\xi,\arccos\left(\pm\langle\frac{\gamma(\xi)}{\vert \gamma(\xi)\vert},\dot\gamma(\xi)\rangle\right)\right), \quad \xi \in \mathbb{S}^1.
\end{equation*}
These graphs contain all 2-periodic points whose chord passes through $0$. In addition, to every consecutive pair of 2-periodic points $\gamma(\xi_1)$ and $\gamma(\xi_2)$ also both the graphs of the restrictions $\delta_\pm$ to $[\xi_1,\xi_2]$ are invariant curves for $T^2$.
\end{lemma}

We now investigate the local behavior of the map $T^2$ around the fixed points $(\xi_i, \pi/2)$, $i=1,2$. 
\begin{lemma}	\label{lemma_differential_is_hyperbolic}
Assume that $\Omega$ admits two distinct focal points of the second kind. Then, both $(\xi_i,{\pi}/{2})$, $i=1,2$, are hyperbolic fixed points for the square of the billiard map.
\end{lemma}
\begin{proof}    
We compute the differential of $T^2$ at $(\xi_i,{\pi}/{2})$, $i=1,2$. Denoting with $d = |\gamma(\xi_1)-\gamma(\xi_2)|$ and $k_i$ the curvature at $\gamma(\xi_i)$, it turns out that (see \cite[Theorem 4.2,Part V]{katok2006invariant}) the Jacobian of $T$ at the same points are
\[
DT\left(\xi_1,\frac{\pi}{2}\right) =
\begin{pmatrix}
k_1d - 1               & d \\
k_1 k_2 d -(k_1 + k_2) & k_2d - 1
\end{pmatrix},
\quad 
DT\left(\xi_2,\frac{\pi}{2}\right) =
\begin{pmatrix}
k_2d - 1               & d \\
k_1 k_2 d -(k_1 + k_2) & k_1d - 1
\end{pmatrix}
\]
hence 
\[
DT^2\left(\xi_1,\frac{\pi}{2}\right) =
DT\left(\xi_2,\frac{\pi}{2}\right)
DT\left(\xi_1,\frac{\pi}{2}\right) = 
\begin{pmatrix}
A & B_1 \\
C_1 & A 
\end{pmatrix}
\]
\[
DT^2\left(\xi_2,\frac{\pi}{2}\right) =
DT\left(\xi_1,\frac{\pi}{2}\right) 
DT\left(\xi_2,\frac{\pi}{2}\right) =
\begin{pmatrix}
A & B_2 \\
C_2 & A 
\end{pmatrix},
\]
whose determinant is clearly 1 and the trace is the double of
\[
A = -1 - 2 k_1k_2 d^2 +2d (k_1+k_2).
\]
By direct computations the discriminant of the characteristic polynomial of both matrices is
\[
4d(dk_1-1)(dk_2-1)(dk_1k_2-k_1-k_2).
\]
Following the proof of Lemma \ref{lemma:at_most_2_bad _points},  we can conclude that, since we are assuming that $c \neq 0$, both matrices have two distinct real eigenvalues.
\end{proof}

\begin{lemma}\label{lemma:stable_unstable}
Assume that $\Omega$ admits  two focal points of the second kind. Then $\delta_{\pm}$ and $\delta^c_{\pm}$, introduced in Lemma \ref{lem:graphs}, locally coincide with the stable/unstable manifolds of the points $(\xi_1,\pi/2)$ and $(\xi_2,\pi/2)$. In particular, $\delta_{\pm}(\mathbb{S}^1) =\delta^c_{\mp}(\mathbb{S}^1)$.
\end{lemma}
\begin{proof}
Without loss of generality, let us suppose that $p_1=\gamma(\xi_1)$ achieves the minimal distance from the origin, while $p_2=\gamma(\xi_2)$ the maximal one. Then the function $\langle {\gamma(\cdot)}/{\vert \gamma(\cdot) \vert}, \dot{\gamma}(\cdot)\rangle$ strictly increases in a neighborhood of $\xi_1$, and strictly decreases near $\xi_2$. Hence the graphs of $\delta_{+}$ and $\delta_{-}$ intersect transversally at $(\xi_1,\pi/2)$ and $(\xi_2,\pi/2)$. \\
We now show that $p_2$ is an attracting fixed point for the restriction of $T^2$ to the graph of $\delta_-$. The argument for the other cases is completely analogous.
If they exist, let us consider the points $q_1, \dots, q_{2n} \in \partial \Omega$ corresponding to all periodic trajectories through the origin different from the chord $p_1p_2$, ordered anti-clockwise as in Figure \ref{fig:enter-label}. Every $q_i$ is fixed by $T^2$ and any arc between consecutive points is invariant. 
For any point $p$ on $\partial \Omega$ between $p_{2}$ and $q_{n+1}$, the line from $p$ through the origin intersects $\partial \Omega$ in a point $p'$ between $p_1$ and $q_1$.
The reflected segment does not intersect $p_10$ since the minimality of $p_1$ implies that the angle it spans with $\dot{\gamma}$ is smaller than $\pi/2$. 
Thus, $T^2(p)$ lies closer to $p_2$ than $p$. This proves that locally $p_2$ is an attracting fixed point on the graph of $\delta_-$. By Lemma \ref{lemma_differential_is_hyperbolic} and Hartman Grobman Theorem, the stable manifold of $(\xi_2,\pi/2)$ locally coincides with  the graph of $\delta_-$.
Analogously, the graph of $\delta_+$ locally coincides with the unstable manifold of $(\xi_2,\pi/2)$. The role of $\delta_-$ and $\delta_+$ exchanges considering $(\xi_1,\pi/2)$.

Next, we consider the analogous invariant graphs associated with the second point $c$, which are $\delta_\pm^c$. With the same reasoning they must locally coincide with the same stable/unstable manifolds of the two fixed points $(\xi_i,\pi/2)$, $i=1,2$.  By the unique continuation principle for analytic maps we can conclude that the two pairs of invariant graphs coincide and necessarily $\delta_{\pm}(\mathbb{S}^1) =\delta^c_{\mp}(\mathbb{S}^1)$ (see Figure \ref{fig:Delta}).
\begin{figure}
\centering
\begin{tikzpicture}[x=0.75pt,y=0.75pt,yscale=-1,xscale=1]

\draw [color={rgb, 255:red, 155; green, 155; blue, 155 }  ,draw opacity=1 ] [dash pattern={on 4.5pt off 4.5pt}]  (100,91) -- (360.33,91) ;
\draw [color={rgb, 255:red, 155; green, 155; blue, 155 }  ,draw opacity=1 ] [dash pattern={on 4.5pt off 4.5pt}]  (322.33,20) -- (110.33,119) ;
\draw  [line width=1.5]  (322.33,20) .. controls (347.33,36) and (359.33,61) .. (360.33,91) .. controls (361.33,121) and (347.33,159) .. (327.33,179) .. controls (307.33,199) and (121.33,137) .. (110.33,119) .. controls (99.33,101) and (100.33,102) .. (100,91) .. controls (99.67,80) and (102.33,77) .. (120.33,60) .. controls (138.33,43) and (297.33,4) .. (322.33,20) -- cycle ;
\draw [color={rgb, 255:red, 155; green, 155; blue, 155 }  ,draw opacity=1 ] [dash pattern={on 4.5pt off 4.5pt}]  (327.33,179) -- (120.33,60) ;
\draw    (105.33,112) -- (338.33,33) ;
\draw    (105.33,112) -- (355.33,62) ;
\draw [color={rgb, 255:red, 103; green, 168; blue, 243 }  ,draw opacity=1 ]   (105.33,112) -- (144.4,101.52) ;
\draw [shift={(146.33,101)}, rotate = 164.98] [color={rgb, 255:red, 103; green, 168; blue, 243 }  ,draw opacity=1 ][line width=0.75]    (10.93,-3.29) .. controls (6.95,-1.4) and (3.31,-0.3) .. (0,0) .. controls (3.31,0.3) and (6.95,1.4) .. (10.93,3.29)   ;
\draw [color={rgb, 255:red, 103; green, 168; blue, 243 }  ,draw opacity=1 ]   (105.33,112) -- (118.55,143.16) ;
\draw [shift={(119.33,145)}, rotate = 247.01] [color={rgb, 255:red, 103; green, 168; blue, 243 }  ,draw opacity=1 ][line width=0.75]    (10.93,-3.29) .. controls (6.95,-1.4) and (3.31,-0.3) .. (0,0) .. controls (3.31,0.3) and (6.95,1.4) .. (10.93,3.29)   ;
\draw  [fill={rgb, 255:red, 0; green, 0; blue, 0 }  ,fill opacity=1 ] (169,90.24) .. controls (169,88.17) and (170.68,86.48) .. (172.76,86.48) .. controls (174.83,86.48) and (176.52,88.17) .. (176.52,90.24) .. controls (176.52,92.32) and (174.83,94) .. (172.76,94) .. controls (170.68,94) and (169,92.32) .. (169,90.24) -- cycle ;

\draw (112,32.4) node [anchor=north west][inner sep=0.75pt]    {$q_{2n}$};
\draw (329.33,182.4) node [anchor=north west][inner sep=0.75pt]    {$q_{n}$};
\draw (321.33,-3.6) node [anchor=north west][inner sep=0.75pt]    {$q_{n+1}$};
\draw (96.33,120.88) node [anchor=north west][inner sep=0.75pt]    {$q_{1}$};
\draw (68.33,81.4) node [anchor=north west][inner sep=0.75pt]    {$p_{1}$};
\draw (375.33,77.38) node [anchor=north west][inner sep=0.75pt]    {$p_{2}$};
\draw (82,108.4) node [anchor=north west][inner sep=0.75pt]    {$p'$};
\draw (361,46.4) node [anchor=north west][inner sep=0.75pt]    {$T^{2}( p)$};
\draw (347,16.4) node [anchor=north west][inner sep=0.75pt]    {$p$};
\end{tikzpicture}
    \caption{An illustration of the proof of Lemma \ref{lemma_differential_is_hyperbolic}. The curves $\delta_\pm$ are contained in the stable/unstable set of $p_1$ and $p_2$.}
    \label{fig:enter-label}
    \end{figure}
\begin{figure}
\includegraphics[width=0.5\textwidth]{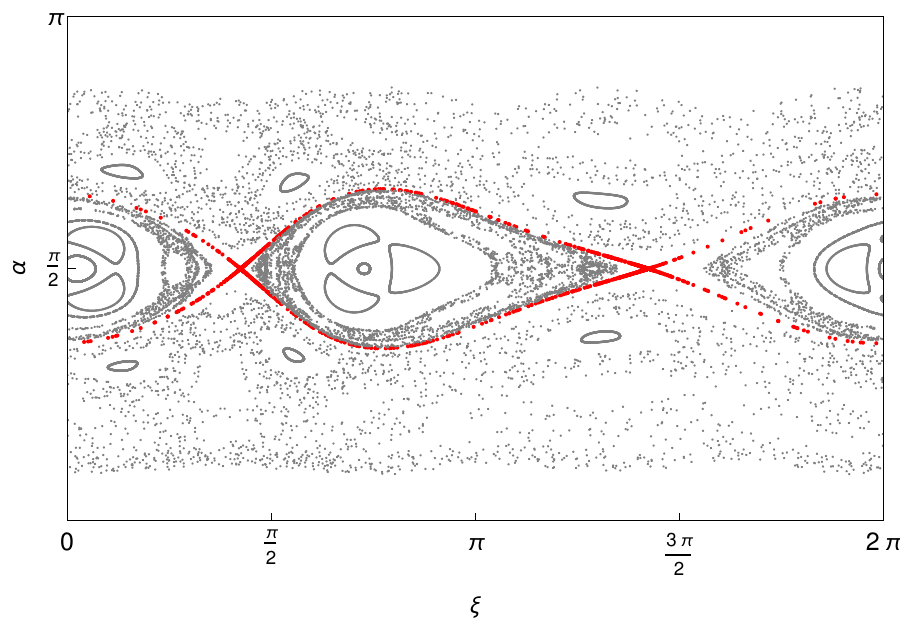}
\caption{Phase portrait of the second iterate $T^2$, choosing $\Omega$ as in Section \ref{subsect:string}. In red, the invariant curves $\delta_{\pm}(\mathbb S^1)$, see  Lemmas \ref{lem:graphs} and \ref{lemma:stable_unstable}. }
\label{fig:Delta}
\end{figure}
\end{proof}
    
\begin{proof}[Proof of Theorem \ref{thm:at_most_one_point}]
From Lemma \ref{lemma:at_most_2_bad _points} we know that $\Omega$ has at most two focal points of the second kind. If there are exactly two, say 0 and $c\neq 0$, Lemma \ref{lemma:stable_unstable} guarantees $\delta_{\pm}(\mathbb{S}^1) =\delta^c_{\mp}(\mathbb{S}^1)$, hence 
\[
\langle\frac{\gamma(\xi)}{\vert \gamma(\xi)\vert},\dot\gamma(\xi)\rangle=-\langle\frac{\gamma(\xi)-c}{\vert \gamma(\xi)-c\vert},\dot\gamma(\xi)\rangle,  \quad \xi \in \mathbb{S}^1.
\]
This implies that $\Omega$ is an ellipse with foci $0$ and $c$. Indeed, this identity implies that the family of invariants lines defined in \eqref{eq:invariant_curve_lines} consists of the pencil of lines through the origin and the one through $c$, drawing triangles having the same perimeter. So, $\partial \Omega$ is an ellipse with foci in $0$ and $c$.
\end{proof}

 \begin{proof}[Proof of Corollary \ref{cor:focal_points_with_symmetry}]
Assume that 0 is a focal point of the second kind for $\Omega$. Since $\partial \Omega$ has a cyclic symmetry of order $p$ there is a point $C$ in the interior of $\Omega$ fixed by a rotation of angle $2 \pi /p$ preserving $\Omega$. Thus, thanks to Theorem \ref{thm:at_most_one_point}, $C\neq 0$ leads to a contradiction when $p>2$ while it imposes that $\partial \Omega $ is an ellipse when $p=2$. \\
On the other hand, when $C=0$ two cases are possible: when $p>2$ the focal point can not be of the second kind. When $p=2$, points of minimal and maximal distance from the origin can not be antipodal unless $\partial \Omega$ is a circle centered at $0$.
 \end{proof}

\subsection{Tables with focal points of the second kind}\label{subsect:string}
Finally, we consider the problem of establishing the existence of non-elliptic analytic domains admitting a focal point of the second kind. We illustrate here the construction of an infinite-dimensional family of such domains via a string-type argument. When an attracting center is located at such a point, our results do not apply and  we are unable to ascertain whether symbolic dynamics is present. Moreover, our construction provides families of curves with focal points arbitrarily close to any ellipse.

\begin{figure}[ht]
\begin{overpic}[height=4.5cm]{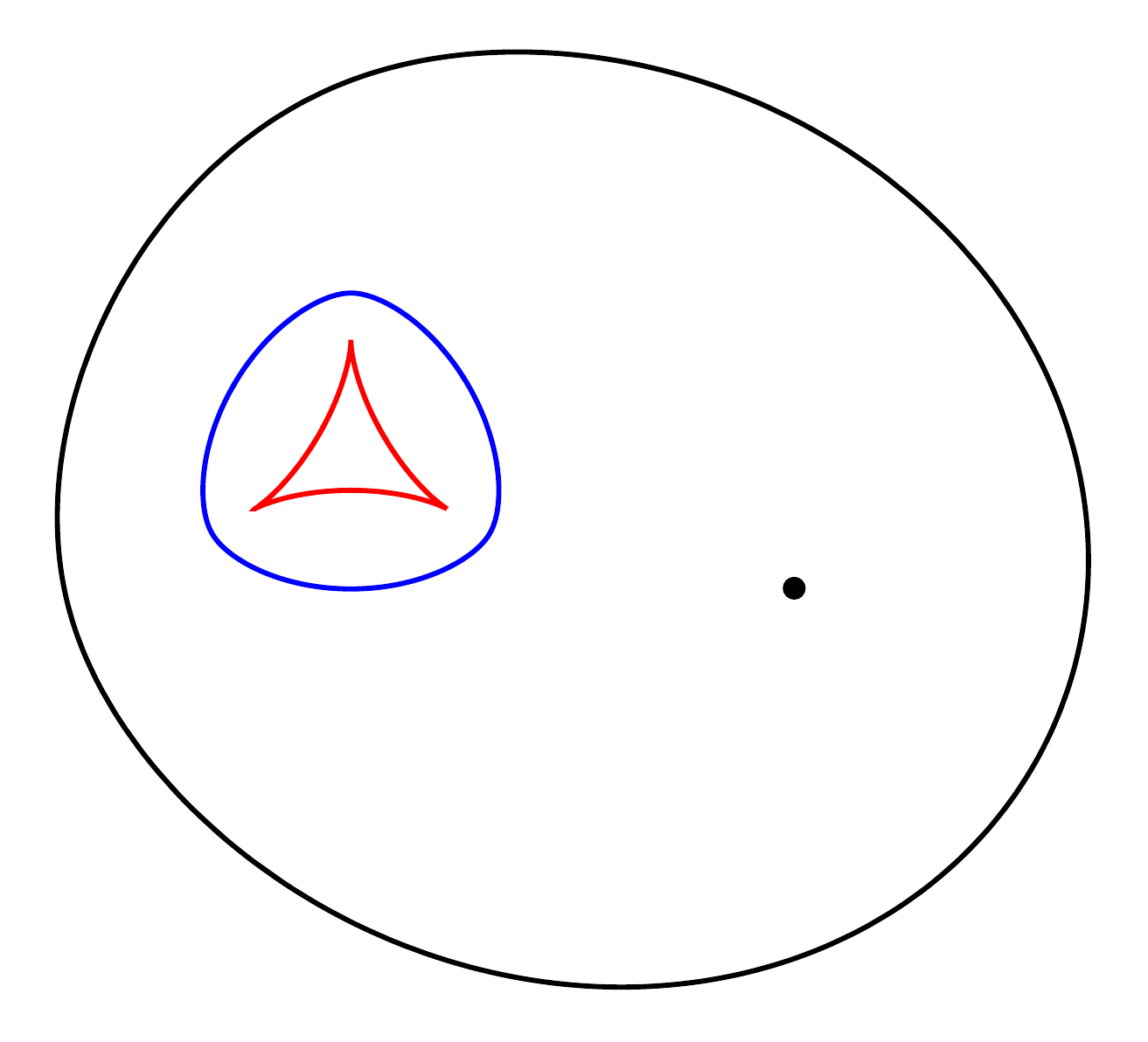}
\put(2,65){$\gamma$}
\put(40,60){$T$}
\put(28,43){$\alpha$}
\put(70,35){$c$}
\end{overpic}   
\hspace{.5cm}
\begin{overpic}[height=4.5cm]{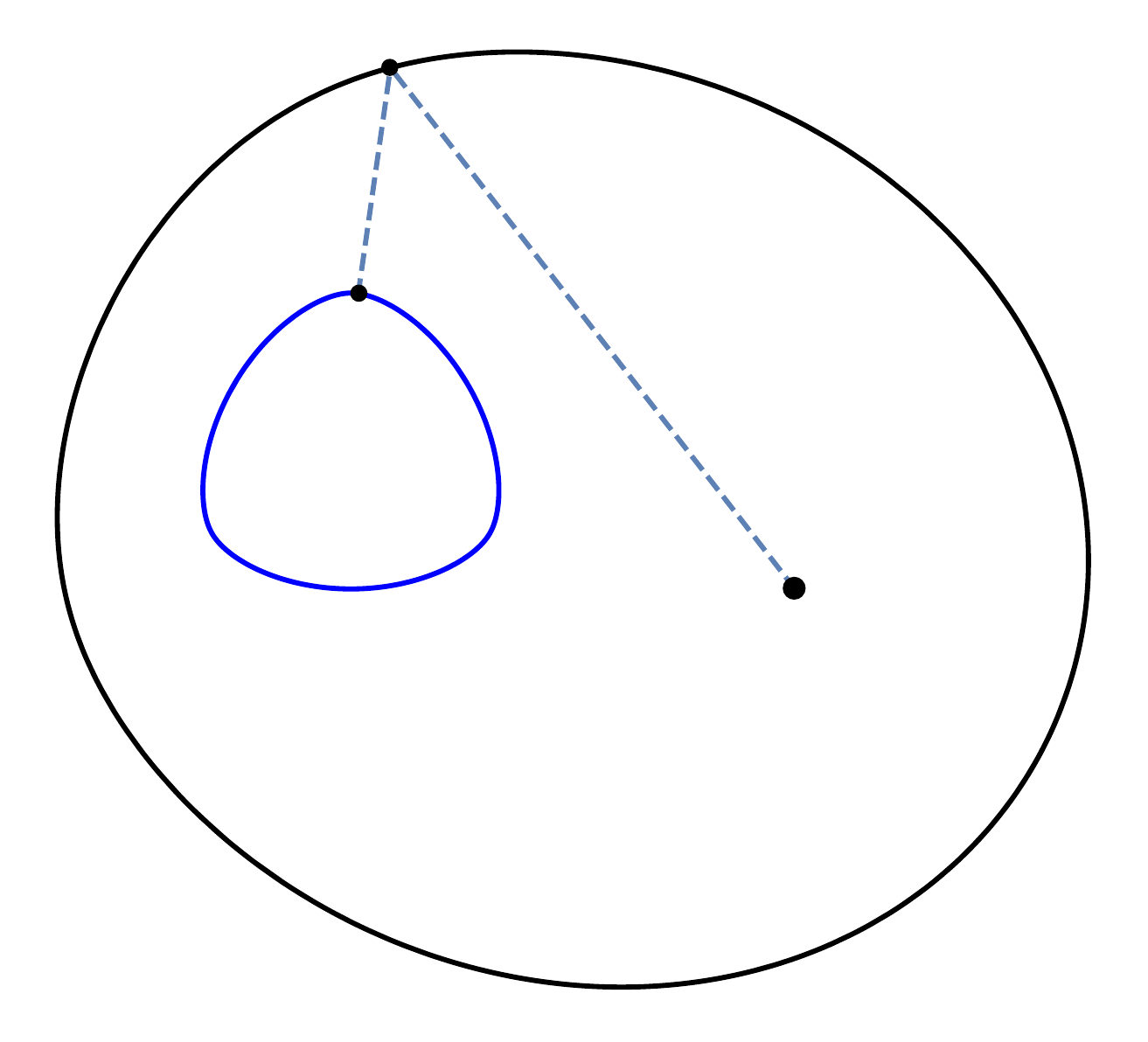}
\put(70,35){$c$}
\put(28,50){$\tau$}
\put(30,88){$x$}
\end{overpic}   
\hspace{.5cm}
\begin{overpic}[height=4.5cm]{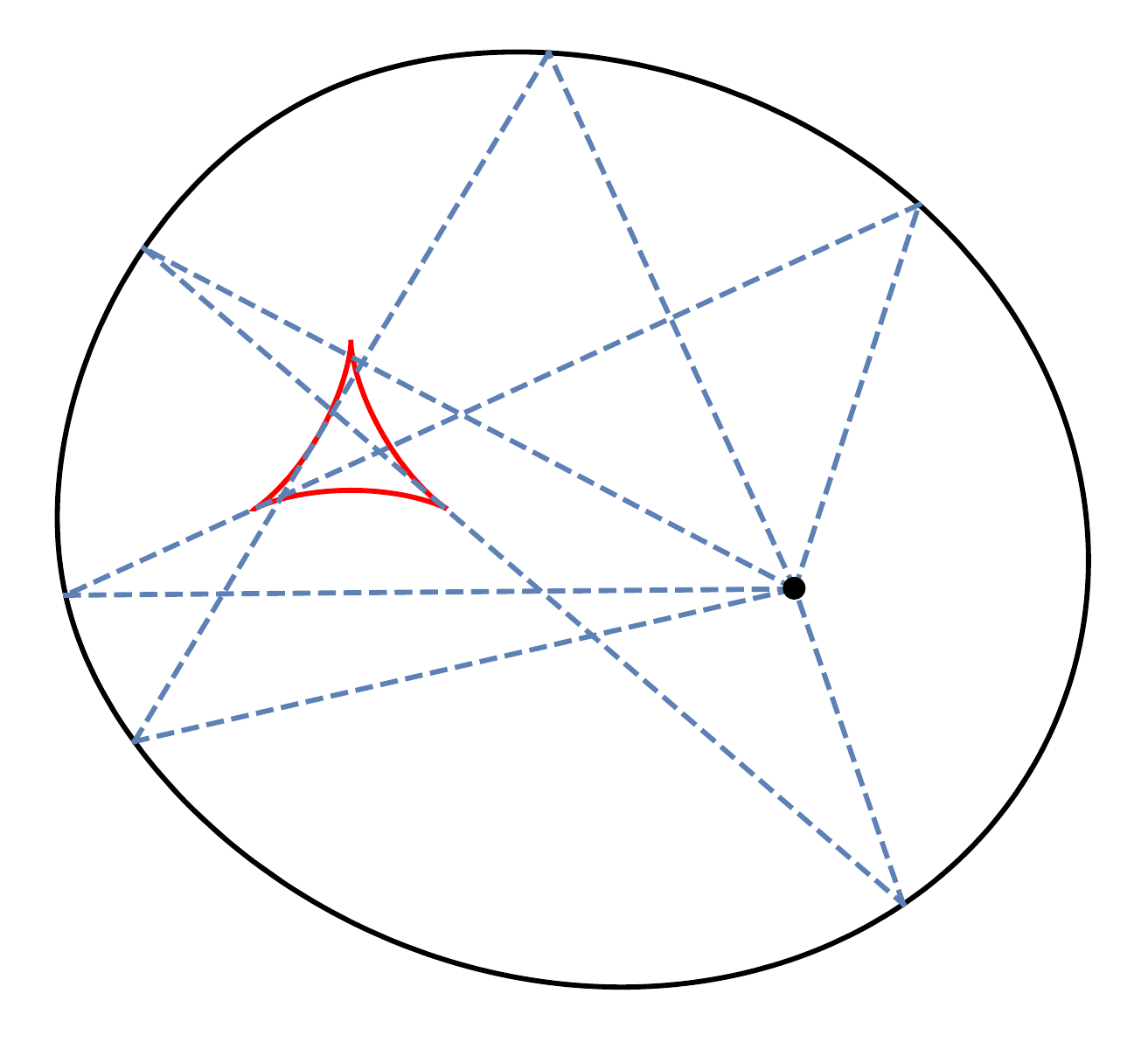}
\end{overpic} 
\caption{Left, the curves $T$, $\alpha$ and $\gamma$ as in Eqs. \eqref{eq:constan_width_caustic} and \eqref{eq:parametrization_curve_focal_point} when $a_0=1$, $a_3=1/3$, with $a_k=0$ when $k\neq0,3$. The point $c$ is located at $(3,0)$ and $\tilde \ell =6$. Center, a point $x$ belongs to $\gamma$ if and only if the sum of length of the dashed segments is constantly equal to $\tilde \ell$ (see Eq. \eqref{eq:f}). Right, Birkhoff triangles centered at $c$ are tangent to the curve $\alpha$.}
\label{fig:gamma}
\end{figure}

We will build our tables starting from \emph{constant-width} domains, following the exposition in \cite{knill_non_convex_caustics}. Let $T$ be an analytic curve and let us parametrize it with constant angular velocity. This means that $T$ can be written as
\[
T(\xi) = T(0) + \int_0^{\xi}\rho(s) e^{i s}ds,\quad  \dot{T}(\xi) = \rho(\xi) e^{i \xi}.
\]
Here the function $\rho(\xi)$ corresponds to the radius of curvature of the curve $T$ at $\xi$. A domain is termed of \emph{constant-width} if it possesses an invariant curve of $2-$periodic billiard trajectories.    
    A general formula by R. Douady (see \cite[Lemma 3.1]{knill_non_convex_caustics}) implies that the caustic $\alpha$ determined by this curve coincides with the evolute of $T$ and is thus given by the formula
    \[
    \alpha(\xi)  = T(\xi) + i \rho(\xi) e^{i \xi} = T(0) + \int_0^{\xi}\rho(s) e^{i s}ds+ i \rho(\xi) e^{i \xi}.
    \]
    Moreover, the radius of curvature must satisfy $\rho(\xi)+\rho(\xi+\pi)= \ell$ for some fixed $\ell>0$. This is geometrically clear since the orthogonal line to $T(\xi)$ and $T(\xi+\pi)$ must determine a unique tangency point on the caustic $\alpha$, which is $\alpha(\xi)$. Let us expand the radius of curvature in Fourier series on $[0,2\pi]$ as $\rho(\xi) = \sum_{k \in \mathbb{Z}} a_k e^{i k \xi}$.
    Summarizing the above discussion we require that $\rho$ is a real, positive function and so $ a_{-k} = \bar{a}_k$. We require that $a_{2k}=0 $ for all $k \ne 0$ to have the symmetry condition $\rho(\xi)+\rho(\xi+\pi)=\ell= 2 a_0$. Finally $a_1 = 0$ in order for the formula for $T$ to define a periodic function.
    With these choices the parametrization of $\alpha$ and $T$ are explicit and given by the following formulas (compare with \cite[Proposition 4.1]{knill_non_convex_caustics} and see Figure \ref{fig:gamma} for an example)
    \begin{equation}
    \label{eq:constan_width_caustic}
       \begin{aligned}
            T(\xi) &= T(0) - i \sum_{k \in \mathbb{Z}} \frac{a_k}{k+1}\left(e^{i(k+1)\xi} -1\right),\\
            \alpha(\xi) &= T(0)+ i \sum_{k \in \mathbb{Z}}\frac{a_k}{k+1}\left(k e^{i(k+1)\xi} +1\right).
       \end{aligned}
    \end{equation}

We are now ready to illustrate a \emph{generalized} string construction for boundaries with focal points.
\begin{figure}[ht]
\begin{overpic}[height=4.5cm]{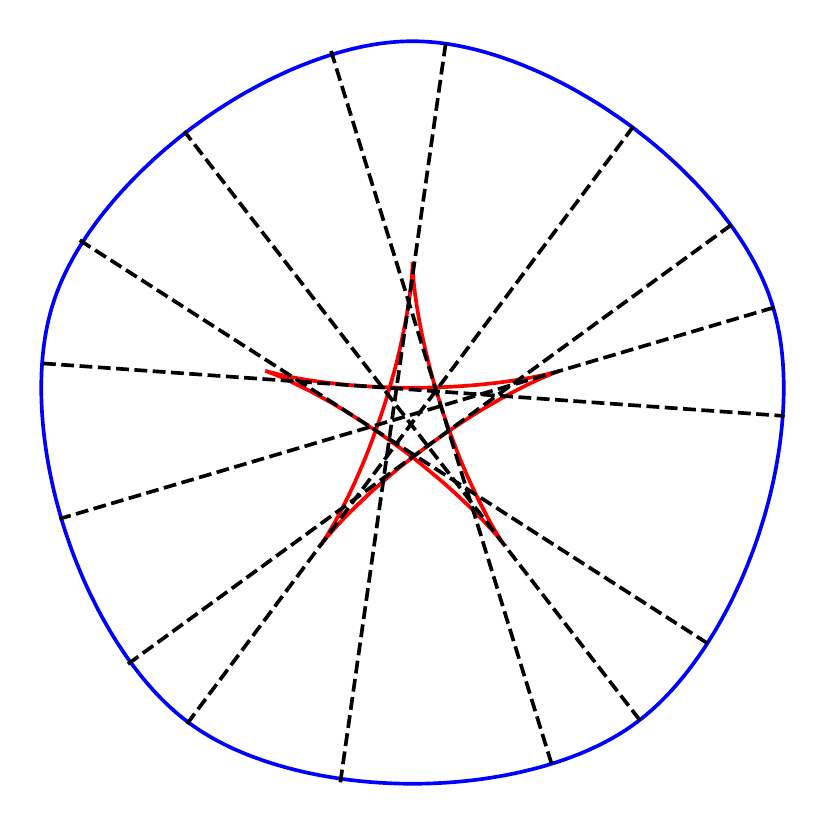}
\put(0,65){$T$}
\put(28,43){$\alpha$}
\end{overpic}   
\hspace{.5cm}
\begin{overpic}[height=4.5cm]{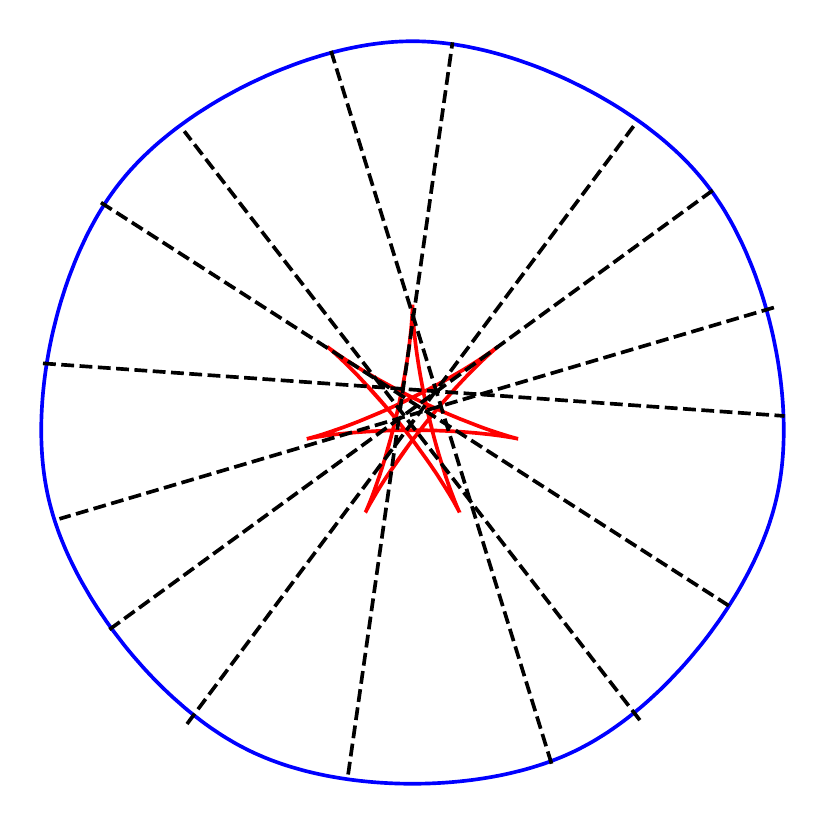}
\put(0,65){$T$}
\put(28,43){$\alpha$}
\end{overpic}   
\caption{Left, the curves $T$ and $\alpha$ as in Eq. \eqref{eq:constan_width_caustic} when $a_0=1$, $a_5=1/5$, with $a_k=0$ when $k\neq0,5$. Right, the curves $T$ and $\alpha$ as in Eq. \eqref{eq:constan_width_caustic} when $a_0=1$, $a_7=1/7$, with $a_k=0$ when $k\neq0,7$. }
\label{fig:gamma57}
\end{figure}

\begin{figure}[ht]
\includegraphics[width=0.7\textwidth]{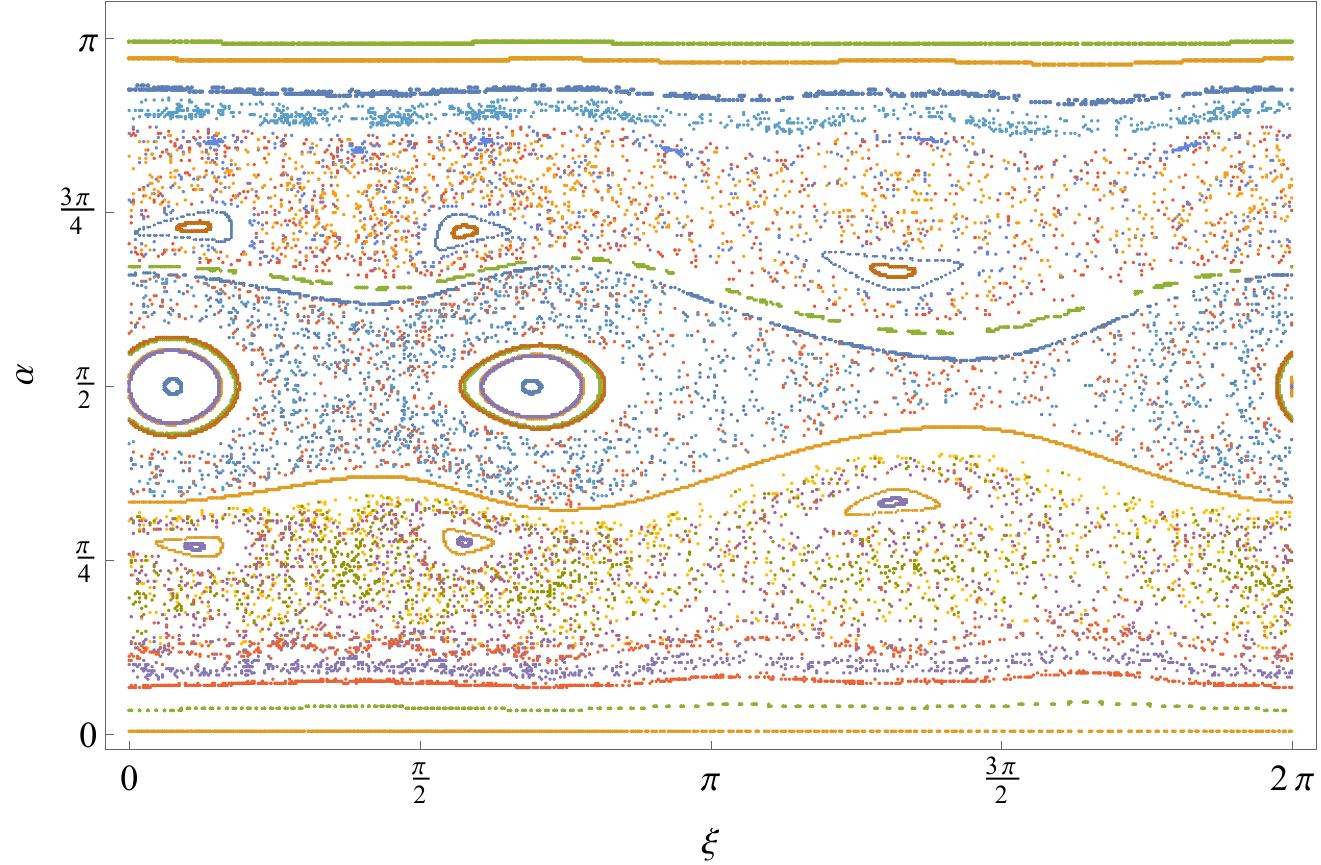}
\caption{First return map associated to the Kepler billiard when the boundary is constructed as in Figure \ref{fig:gamma}. The physical parameters associated to $V(z)$ are $h=10$ and $\mu=5.$ The simulation has been obtained with the software $Mathematica^\copyright$.}
\label{fig:primo_ritorno}
\end{figure}

\begin{proposition}
\label{prop:curve_with_focal_points}
 For any analytic curve $T$ as in \eqref{eq:constan_width_caustic} bounding a convex constant-width body $\tau$ and any point $c \in \mathbb{R}^2 \setminus \tau$ there exists a 1-parameter family of convex analytic curves having $c$ as a focal point of the second kind.
\end{proposition}
\begin{proof}
We consider the family of lines orthogonal to the boundary of $T$, they can be explicitly parametrized as
\[
v(\xi,s) = T(\xi) - s i e^{i \xi}, \quad (\xi, s) \in \mathbb{S}^1 \times \mathbb{R}.
\]
Any point $x$ in $\mathbb{R}^2\setminus \tau$ can be written uniquely as $x = T(\xi) -s ie^{i \xi}$ for some $s>0$ and $\xi \in \mathbb{S}^1$. In analogy with the classical string construction, we look for level sets of the sum of the distances from $c$ and $\tau$. This reduces to the following equation in $s$ and $\xi$ 
\[
\vert  T(\xi) - s i e^{i \xi}- c\vert + s = \tilde{\ell}, \quad  \tilde{\ell}>0.
\]
A straightforward computation shows that
\[
s(\xi) = \frac12\frac{\tilde{\ell}^2- \vert T(\xi)-c \vert^2}{\tilde{\ell}-\langle T(\xi)-c,i e^{i \xi }\rangle}
\]
and thus, for sufficiently big $\tilde{\ell}$, the curve 
\begin{equation}
\label{eq:parametrization_curve_focal_point}
    \gamma(\xi) = T(\xi) -\frac12  \frac{\tilde{\ell}^2- \vert T(\xi)-c \vert^2}{\tilde{\ell}-\langle T(\xi)-c,i e^{i \xi }\rangle} i e^{i\xi}
\end{equation}
    is well defined. Since the curve $\gamma$ parametrizes the level set of the function
    \begin{equation}\label{eq:f}
    f(x) = \vert x-c\vert+d(x,\tau)
    \end{equation}
    then it is analytic on $\mathbb{R}^2\setminus (\partial \tau \cup \{c\})$
    as soon as $\nabla f \ne 0$ and $\tilde \ell \ge \max_{y\in \tau \cup \{c\}}f(y)$. Since the gradient of the distance function at a point $x$ is the unitary vector pointing in the direction of the chord realizing the distance, $\nabla f(x)$ is zero if and only if $x$ belongs to the segment having $c$ as one endpoint and orthogonal to $\partial \tau$ at the other.
    For an analogous reason, the table $\gamma$ has a unique periodic trajectory passing through $c$, contained in the line orthogonal to $T$ and passing through $c$.  Thus, for $\tilde{\ell}$ large enough, equation \eqref{eq:parametrization_curve_focal_point} parametrizes an analytic curve.
    Clearly, the corresponding function $\varphi_c$ is constant since all chords orthogonal to $T$ have the same length.
    \end{proof}
    Figure \ref{fig:primo_ritorno} shows the orbits of the Kepler billiard map $T_K$ when $\gamma$ is of the form \eqref{eq:parametrization_curve_focal_point}, with $\alpha_0=1$, $\alpha_{\pm}=1/3$, $a_k=0$ when $k\neq0,\pm3$, $c=(3,0)$ and $\tilde \ell=6$, and physical parameters $h=10$ and $\mu=5$. It is clear how, even if Theorem \ref{thm:dyn_sym} does not apply to this case, the resulting dynamics presents strong evidences of chaotic behavior, both close and far from collisional arcs. It is noticeable the presence of invariant curves close to the boundary of $\Omega$. Such numerical results motivates the formulation of Conjecture \ref{conj:main}. 

\begin{figure}[ht]
\begin{overpic}[height=4.5cm]{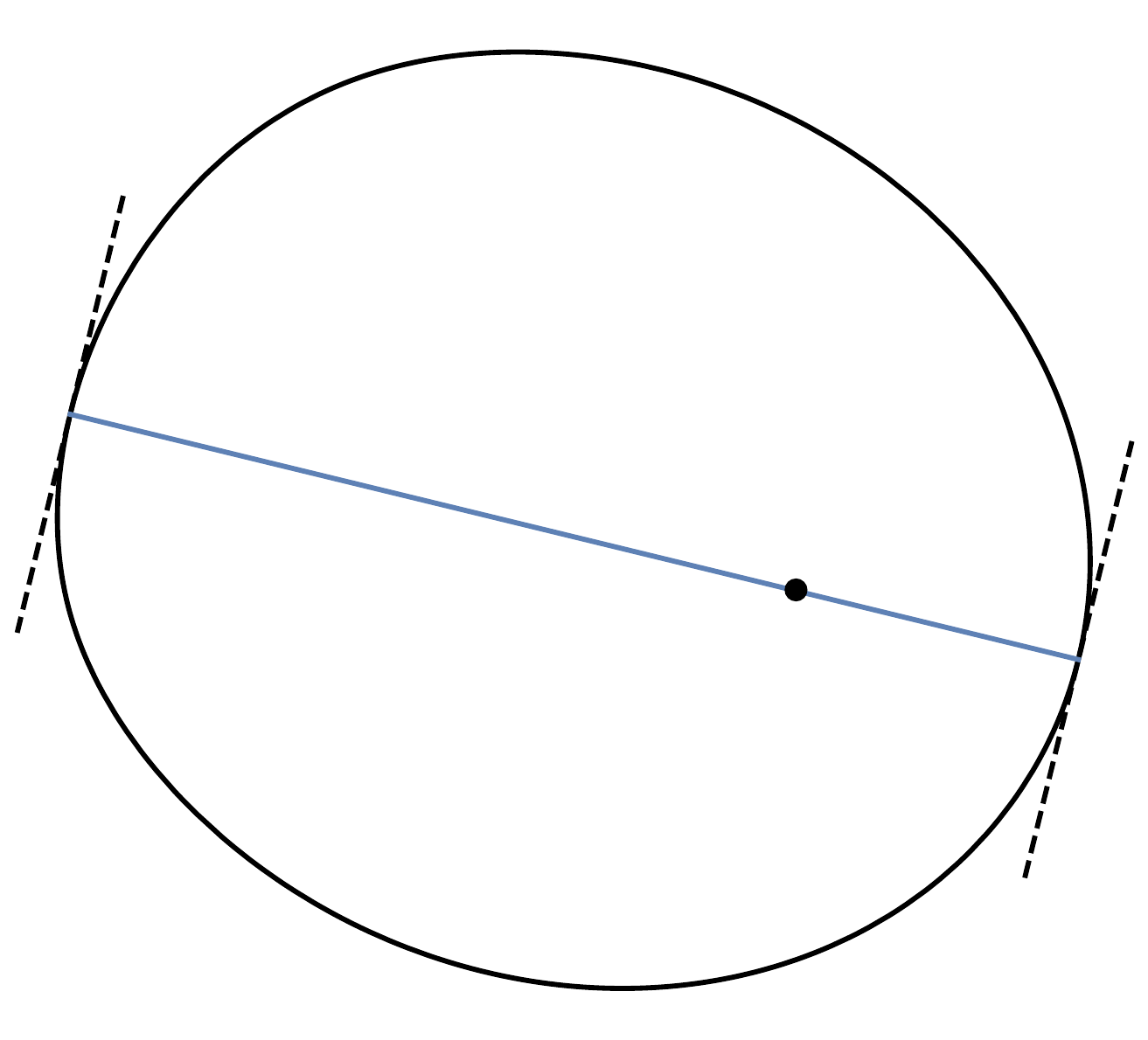}
\put(15,80){$\gamma$}
\put(70,33){$c$}
\end{overpic}   
\hspace{1.5cm}
\begin{overpic}[height=4.5cm]{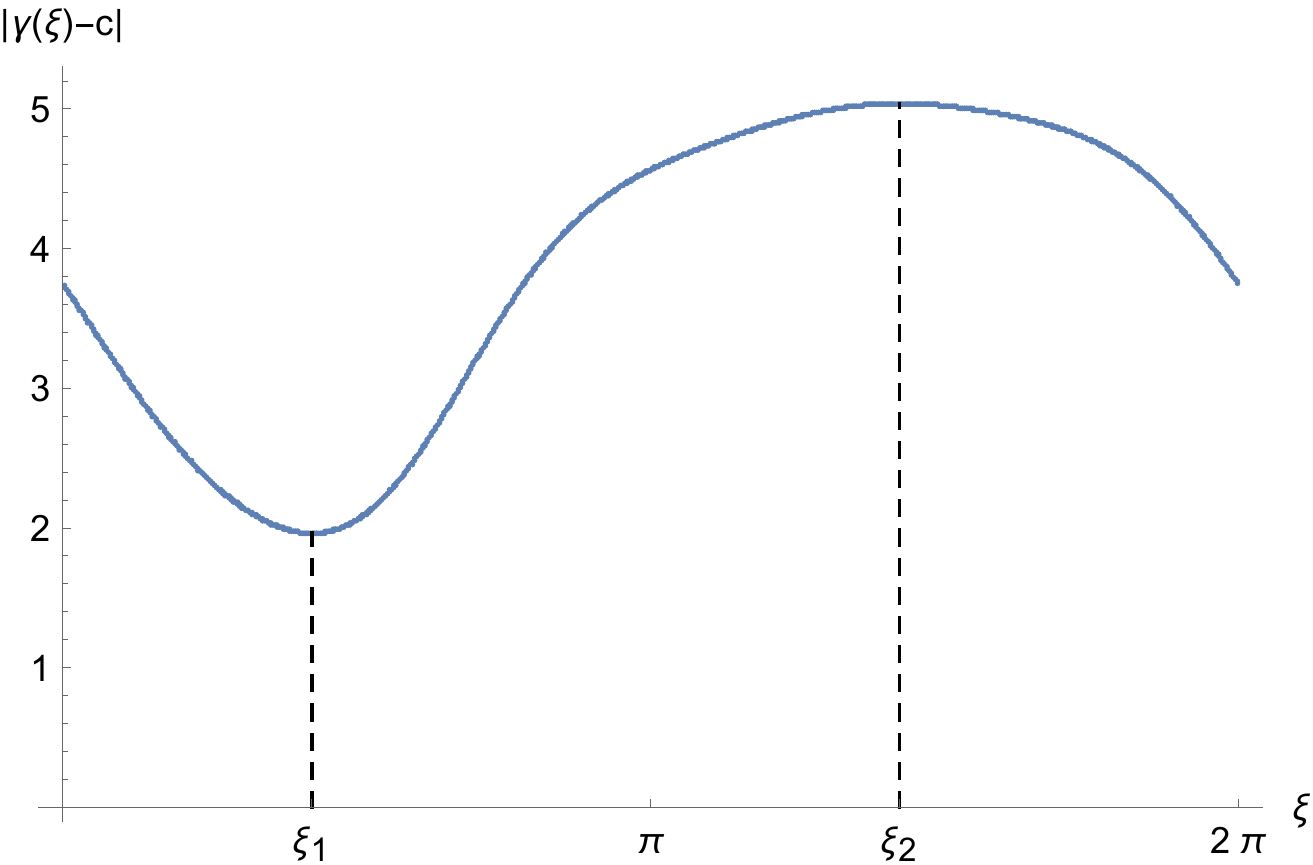}
\end{overpic}   
\caption{In the same setting of Figure \ref{fig:gamma}, $c$ is a focal point of the second kind. In particular in this case there exists a unique orthogonal chord (left)  and the function $|\gamma(\cdot)-c|$ admits exactly two critical points (right).}
\label{fig:gamma2}
\end{figure}

\section{Final remarks and comments}

\begin{enumerate}
    \item \label{first}
    The notions of integrability considered in~\cite{avila2016integrable,MR3815464,MR3788206,BiaMir22} are significantly stronger than the one adopted in the present work. In particular, in~\cite{zbMATH00561292,BiaMir22}, integrability requires that the entire phase cylinder be foliated by invariant curves, either globally or at least for rotation numbers exceeding $1/4$. In~\cite{avila2016integrable,MR3815464,MR3788206}, by contrast, the persistence of KAM tori justifies requiring the existence of invariant curves only for rational rotation numbers below a fixed threshold, under the additional assumption that all such curves consist solely of periodic points. This concept of integrability aligns with the Arnold-Liouville Theorem, which asserts that a system admitting action-angle coordinates near the boundary necessarily yields such a foliation.
    
    \item In contrast, the present paper concerns \emph{analytic integrability}, namely, the existence of analytic first integrals. Absent further assumptions on the behavior of such an integral near the boundary (see, e.g.,~\cite[Corollary 1.3]{BiaMir22}), the billiard map may admit an analytic first integral without being integrable in the stronger sense described in point~(\ref{first}).
    
    \item It is also worth noting that Theorem~\ref{thm:symbolic} yields an invariant subset indicative of a \emph{Birkhoff instability region}. However, the emergence of such a region appears unrelated to the rigidity phenomena addressed in~\cite{MR3788206}, for two main reasons. First, the presence of a Birkhoff instability region near the boundary does not directly follow from the notion of non-integrability adopted in that work. Second, the invariant chaotic set constructed here crucially depends on interactions with the singularity—a distinctive feature of Kepler billiards not shared by classical Birkhoff systems. Hence our result does not contribute to the solution of the classical Birkhoff-Poritsky conjecture.

    \item In this paper we focus on high energy regimes. It was brought to our attention by Lei Zhao in a personal communication that the zero energy case is actually conjugated with the usual Birkhoff billiard on the double covering of the domain. As such, it can be covered by the results of \cite{zbMATH00561292,BiaMir22}.
    
    \item Our proof further provides a lower bound for the topological entropy of the system. This bound may be improved depending on the geometry of $\Omega$ and the number of symbolic elements employed in the shadowing construction.
    
    \item Finally, we observe that our numerical simulations for a non elliptical Kepler billiard with the gravitational center placed at a focal point suggest that chaotic behavior persists even in this configuration, thereby lending support to our Conjecture that the only tables which are integrable at all energy levels are ellipses with the Keplerian center located at one of its foci.
\end{enumerate}

%% file: appendix.tex
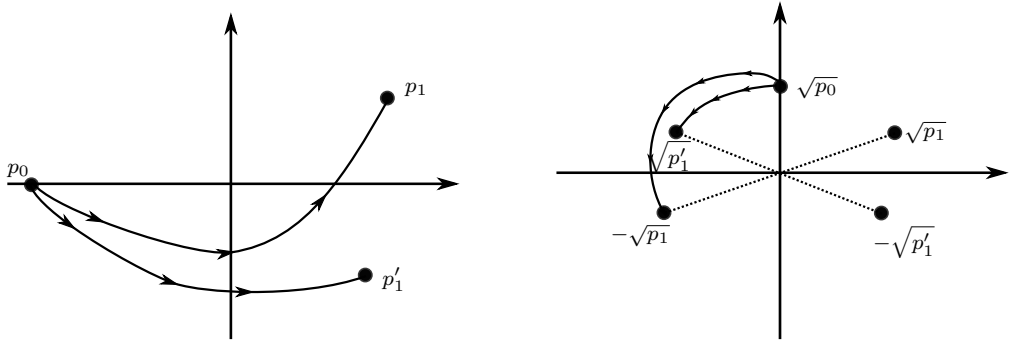
\begin{figure}[t]
	\begin{overpic}[width=0.35\linewidth]{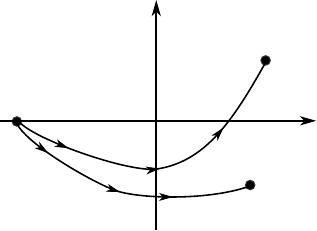}
		\put (88, 55) {\tiny$p_1$}
		\put (83, 12) {\tiny$p_1'$}
		\put (0, 37) {\tiny$p_0$}
	\end{overpic}
    \hspace{1cm}
    \begin{overpic}[width=0.35\linewidth]{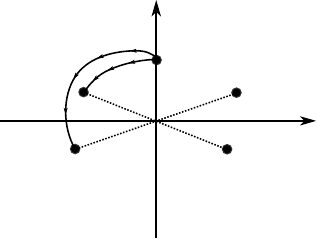}
		\put (53, 55) {\tiny$\sqrt{p_0}$}
        \put (70,20) {\tiny$-\sqrt{p_1'}$}
		\put (77, 45) {\tiny$\sqrt{p_1}$}
		\put (12,22) {\rotatebox{0}{\tiny$-\sqrt{p_1}$}}
        \put (20,39) {\rotatebox{0}{\tiny$\sqrt{p_1'}$}}	
	\end{overpic}
	\caption{The anticlockwise and clockwise arcs connecting two (possibly antipodal) points. Depending on the difference in the arguments between $p_0$ and $p_1$, they correspond either to  direct or  indirect arcs. }
	\label{fig:w01}
\end{figure}

The existence of the two arcs comes again from \cite{jacobi1837theorie}; as for Eq.\eqref{eq:Lzazc}, let us consider only the case $L(z_a)$, the reasoning for $L(z_c)$ is completely analogous. \\
Let us start by observing that, whenever $\det(p_0,p_1)>0$, the arc $z_a(\cdot; p_0,p_1; h )$ corresponds to the direct one $z_d(\cdot; p_0, p_1; h)$, while when $\det(p_0,p_1)<0$ one has $z_a(\cdot; p_0,p_1; h )=z_i(\cdot; p_0,p_1; h )$ ($z_d$ and $z_i$ have been introduced in Lemma \ref{lem: asintotiche_sus}). Furthermore, the function $f_a$ is of class $C^1$, and by Eq. \eqref{eq:asintotiche} the lengths $L(z_{i/d})$ tend to $\sqrt{h}f_a$ whenever $\det(p_0,p_1)\neq 0$: one can then argue by continuity to obtain the continuity at order zero of $L(z_a)$. To study the convergence of its derivatives, let us recall that (see \cite{IreneSusViNEW})
\begin{equation*}
    \begin{aligned}
    &\nabla_{p_0} L(z_a(\cdot; p_0, p_1; h)) = -\sqrt{h + \frac{\mu}{|p_0|}}\frac{z_a'(0;  p_0, p_1; h)}{|z_a'(0;  p_0, p_1; h)|}, \\
    &\nabla_{p_1} L(z_a(\cdot; p_0, p_1; h)) =  \sqrt{h + \frac{\mu}{|p_1|}}\frac{z_a'(T;  p_0, p_1; h)}{|z_a'(T;  p_0, p_1; h)|}; 
    \end{aligned}
\end{equation*}
without loss of generality, we will study only $\nabla_{p_1} L(z_a)$. Let us pass to the regularized problem, considering the Levi-Civita transformation (see \cite{Levi-Civita}). In complex notation
    \[
    \frac{d}{dt}=\frac{1}{\sqrt{2}|z(t(\tau))|}\frac{d}{d\tau}, \quad w^2(\tau)=z(t(\tau)), 
    \]
    obtaining the harmonic repulsor system for a suitable $\tilde T$
    \begin{equation}\label{eq:prob_LC}
    \begin{cases}
    w''(\tau)=w(\tau) & \tau \in [0, \tilde T]\\
    \dfrac12 |w'(\tau)|^2-\frac{1}{2}|w(\tau)|^2=h' \quad & \tau \in [0, \tilde T]\\
    w(0)=w_0, \ w(\tilde T)=w_1
    \end{cases}
    \end{equation}
    where 
    \[
    h'=\dfrac{\mu}{2h}, \ w_0^2=p_0, \ w_1^2=p_1.  
    \]
    If $w(\tau)$ is the regularized solution corresponding to $z_a$, one has that 
    \begin{equation} \label{eq:grad}
    \nabla_{p_1}L(z_a)=\sqrt{h}\frac{w_1}{|w_1|^2}w'(\tilde T), 
    \end{equation}
    so that we are left to study the behavior of $w'(\tilde T)$ under variations of $w_1$. To obtain the arc corresponding to $z_a$ in the Levi-Civita plane, let us fix $w_0=\sqrt{p_0}$, and choose $w_1=\sqrt{p_1}$ when $\det(p_0, p_1) \geq 0$, while $w_1=-\sqrt{p_1}$ whenever $\det(p_0, p_1)<0$ (see Figure \ref{fig:w01}). In the regularized problem, the transition from $z_i$ to $z_d$ corresponds to a sign change of $\langle w_0, w_1\rangle$ (from positive to negative).\\
    Once $w_0,w_1$ are fixed, the solution of \eqref{eq:prob_LC} is given by 
    \begin{equation*}
        w(\tau)=\dfrac{w_1-w_0e^{-\tilde T}}{e^{\tilde T}-e^{-\tilde T}} e^\tau + \dfrac{w_0e^{\tilde T}-w_1}{e^{\tilde T}-e^{-\tilde T}} e^{-\tau},  
    \end{equation*}
    so that 
    \begin{equation}\label{eq:wprimoT}
        w'(\tilde T)=\dfrac{e^{\tilde T}+e^{-\tilde T}}{e^{\tilde T}-e^{-\tilde T}}w_1-\dfrac{2}{e^{\tilde T}-e^{-\tilde T}} w_0. 
    \end{equation}
    From the energy conservation law, one has that 
    \[
    e^{\tilde T}+e^{-\tilde T}=-\dfrac{\langle w_0, w_1 \rangle}{h'}+\sqrt{\dfrac{2(|w_0|^2+|w_1|^2)}{h'}+\left(\dfrac{\langle w_0, w_1 \rangle}{h'}\right)^2+4}; 
    \]
    let us consider the quantity $y=1/(e^{\tilde T}+e^{-\tilde T})$, looking for estimates for $h'\to 0^+$ which are uniform in $w_0, w_1$ when $|\langle w_0, w_1 \rangle|\leq C$, $C>0$ being a suitable constant. \\
    If $\langle w_0, w_1 \rangle\geq 0, $ one has 
    \[
    y=\frac{\langle w_0, w_1 \rangle + \sqrt{2 (|w_0|^2+|w_1|^2)h' + (\langle w_0, w_1 \rangle)^2+4 h'^2}}{2 (|w_0|^2+|w_1|^2) + 4 h'}, 
    \]
    which converges to 
    \[
    y_0^+ = \frac{\langle w_0, w_1 \rangle}{|w_0|^2+|w_1|^2} \quad \text{when } h'\to0.
    \]
    Let us then study $(y-y_0^+)/\sqrt{h'}$, calling $\alpha=\langle w_0, w_1 \rangle/\sqrt{h'}\geq0$: 
    \[
    \frac{y-y_0^+}{\sqrt{h'}}=\frac{(|w_0|^2+|w_1|^2)\left[\sqrt{2(|w_0|^2+|w_1|^2)+ \alpha^2+4 h'}-\alpha\right]-4\sqrt{h'}\langle w_0, w_1 \rangle}{(|w_0|^2+|w_1|^2)\left[2(|w_0|^2+|w_1|^2)+4h'\right]},
    \] 
    which is uniformly bounded as $h'\to0$. We can then conclude that, whenever $\langle w_0, w_1 \rangle\geq 0, $
    \[
    y=y_0^+ +\sqrt{h'}g_1(w_0,w_1; h'), 
    \]
    the function $g_1$ being bounded in all variables. On the other hand, when $\langle w_0, w_1 \rangle<0$, we can write 
    \[
    y = \frac{h'}{-\langle w_0, w_1 \rangle + \sqrt{2 (|w_0|^2+|w_1|^2) h' + \langle w_0, w_1 \rangle^2 + 4 h'^2}} \leq \frac1{2(|w_0|^2+|w_1|^2)}\sqrt{h'}.
    \]
    Hence
    \[
    y = \sqrt{h'}g_2(w_0,w_1; h'), 
    \]
    $g_2$ being bounded in all variables. As a consequence, whenever $w_0,w_1$ belong to a compact set of $\C$, $|\langle w_0, w_1 \rangle|\leq C$, we find 
    \[
    y = y_0+\sqrt{h'}g(w_0, w_1;h'), \quad y_0=\begin{cases}
        \frac{\langle w_0, w_1 \rangle}{|w_0|^2+|w_1|^2} \quad & \text{if }\langle w_0, w_1 \rangle\geq 0\\
        0 & \text{otherwhise}
    \end{cases} 
    \]
    for some bounded function $g$. 

    Let us return to Eq. \eqref{eq:wprimoT}, and observe that, by the regularity of the involved functions (note that $y\in (0,1/2)$), 
    \begin{equation*}
    \begin{aligned}
        \frac{e^{\tilde T}+e^{-\tilde T}}{e^{\tilde T}-e^{-\tilde T}} = \frac{1}{\sqrt{1-4y^2}}=\frac{1}{\sqrt{1-4y_0^2}}+\sqrt{h'}\ a(w_0,w_1; h')\\
         \frac{1}{e^{\tilde T}-e^{-\tilde T}}=\frac{y}{\sqrt{1-4y^2}}=\frac{y_0}{\sqrt{1-4y_0^2}}+\sqrt{h'}\ b(w_0,w_1; h'), 
        \end{aligned}
    \end{equation*}
    so that one obtains
    \begin{equation*}
    \displaystyle 
        w'(\tilde T)=v_0 + \sqrt{h'}G(w_0,w_1; h'), \quad v_0=\begin{cases}
            \dfrac{(|w_0|^2+|w_1|^2)w_1 -2 \langle w_0, w_1 \rangle w_0}{\sqrt{(|w_0|^2+|w_1|^2)^2-4\langle w_0, w_1 \rangle^2}}, \quad &\text{if }\langle w_0, w_1 \rangle\geq0\\
            w_1, &\text{otherwise},
        \end{cases}
    \end{equation*}
    for some bounded function $G$. {Replacing in Eq. \eqref{eq:grad} we have 
    \begin{equation*}
    \nabla_{p_1}L(z_a(\cdot; p_0, p_1; h))=\sqrt{h}\frac{w_1}{|w_1|^2}v_0 + 
     g_a(p_0, p_1; h),
    \end{equation*}
    for some bounded function $g_a$, and to conclude we need to show that 
    $\nabla_{p_1} f_a(p_0,p_1) = \frac{w_1}{|w_1|^2}v_0$. When $\langle w_0, w_1 \rangle<0$, the identity is trivial since $w_1^2=p_1$. Assuming now $\langle w_0, w_1 \rangle\geq0$, recalling that $w_0^2=p_0$ as well, we conclude verifying that
    \[
    \frac{p_1-p_0}{|p_1-p_0|} = \frac{w_1}{|w_1|^2}\dfrac{(|w_0|^2+|w_1|^2)w_1 -2 \langle w_0, w_1 \rangle w_0}{\sqrt{(|w_0|^2+|w_1|^2)^2-4\langle w_0, w_1 \rangle^2}}
    \]
    which is equivalent to the following identity in the complex plane
    \[
    |w_0|^2 w_1^2 + |w_1|^2 w_0^2 = 2 \langle w_0,w_1\rangle w_0w_1,
    \]
    since $|p_1-p_0| = \sqrt{(|w_0|^2+|w_1|^2)^2-4\langle w_0, w_1 \rangle^2}$.
    }
    This concludes the proof.
    \qed

%% file: appendix2.tex
\label{app:proof_critical_points}
\begin{proof}[Proof of Lemma \ref{lemma:critical_points_psi_phi}]
    Let us compute the partial derivatives of $\Psi$ assuming, without loss of generality, that $\gamma$ is parameterized by arc length:
    \begin{align*}
 	\partial_\xi \Psi (\xi,\eta) &= \langle \frac{\gamma(\xi)}{\vert \gamma (\xi) \vert}-\frac{\gamma(\eta)-\gamma(\xi)}{\vert\gamma(\eta)-\gamma(\xi)\vert}, \dot{\gamma}(\xi)\rangle \\
 	\partial_\eta \Psi(\xi,\eta) &= \langle \frac{\gamma(\eta)}{\vert \gamma (\eta) \vert}-\frac{\gamma(\xi)-\gamma(\eta)}{\vert\gamma(\xi)-\gamma(\eta)\vert}, \dot{\gamma}(\eta)\rangle. 
    \end{align*}            
    Thus, $\partial_\xi \Psi=0$ if and only if the oriented line spanned by $\gamma(\xi)$ is reflected into the one spanned by $\gamma(\eta)-\gamma(\xi)$. The same conclusion holds for $\partial_\eta \Psi$.   
    
    \noindent  We prove (i). The function $B$ is defined by the equation $\partial_\xi \Psi(\xi,B(\xi)) =0$, thus we have
    \[
    \varphi'(\xi) = \partial_\xi\Psi(\xi,B(\xi)) + B'(\xi) \partial_\eta\Psi(\xi,B(\xi)) =  B'(\xi) \partial_\eta\Psi(\xi,B(\xi)). 
    \]
    Hence critical points of $\Psi$ are contained in the set $\{(\xi',B(\xi')): \varphi'(\xi')=0\}$.

     \noindent We now prove (ii) Assume that $\varphi$ is constant.  Thanks to (i), all critical points are necessarily at level $M$. Similarly, if $\varphi$ is constant, all points of the form $(\eta,B(\eta))$ are maxima of $\Psi$ and thus critical points. So, $\Psi^{-1}(M)$ is connected.  
\end{proof}

\begin{proof}[Proof of Proposition \ref{prop:critici_Psi}]       Assume $\varphi$ is constant. Thanks to Lemma \ref{lemma:critical_points_psi_phi} (ii),  $\Psi^{-1}(M)$ is connected and all critical points of $\Psi$ are at level $M$. Hence there  are no other isolated critical points  for $\Psi$.

       Assume now that $\varphi$ is not constant.
       Since the boundary $\partial\Omega$ is analytic, the function $\varphi$ is analytic as well. Thus, the set of its critical points is discrete.
       Thanks to Lemma \ref{lemma:critical_points_psi_phi} (i), also the set of critical points of $\Psi$ is discrete, in particular $\Psi$ admits, at least, a maximum point in $\mathcal{U}$.
       Consider the vector field $\nabla \Psi$, continuously defined on the open set $\mathcal{U}$. This set is diffeomorphic to an open cylinder whose boundary components correspond to two copies of the diagonal $\Delta$ with opposite orientation. We denote by $\Sigma$ the corresponding closed cylinder (see Figure \ref{fig:cylinder}).
        \begin{figure}
        \begin{overpic}[height=3cm]{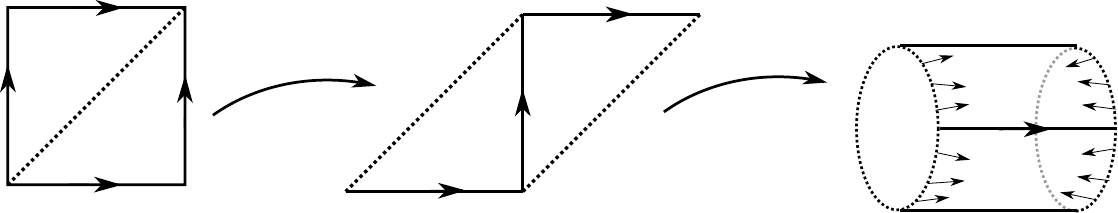}
        \put(8, 8){$\Delta$}
        \end{overpic}
        \caption{The open set $\mathcal U$ (left) is diffeomorphic to an open cylinder whose boundary components correspond to two copies of $\Delta$ (right). The vector field $\nabla \Psi$ extends continuously to the closure of such cylinder, $\Sigma$, and points inward along its boundary.} \label{fig:cylinder}
        \end{figure}
        Now we show that the vector field $\nabla \Psi$ extends continuously to $\Sigma$ and it is always inward-pointing along the boundary. Indeed, $\nabla \Psi$ extends continuously if and only if the gradient of $\vert \gamma(\xi)-\gamma(\eta)\vert $ does so. We have
        \[
        \nabla \left(\vert \gamma(\xi)-\gamma(\eta)\vert\right)  = 
            \begin{pmatrix}
            \frac{\langle\gamma(\xi)-\gamma(\eta), \dot{\gamma}(\eta)\rangle}{\vert\gamma(\xi)-\gamma(\eta)\vert}\\
            - \frac{\langle\gamma(\xi)-\gamma(\eta), \dot{\gamma}(\xi)\rangle}{\vert\gamma(\xi)-\gamma(\eta)\vert}
            \end{pmatrix},
        \]
        hence, since we have chosen the arc-length parametrization, we obtain
        \[
        \lim_{\xi\to \eta^\pm}   \frac{\gamma(\xi)-\gamma(\eta)}{\vert\gamma(\xi)-\gamma(\eta)\vert}  = \pm \dot\gamma (\eta). 
        \]
        So, we can extend $\nabla \Psi$ to the boundary of the cylinder $\Sigma$ setting
        \[
        \nabla \Psi (\eta,\eta)= \frac{\langle \gamma(\eta), \dot\gamma(\eta)\rangle}{ \vert \gamma(\eta)\vert} 
        \begin{pmatrix} 1\\1 \end{pmatrix}\pm 
        \begin{pmatrix} 1\\-1 \end{pmatrix}
        \]
        which is always pointing inward $\Sigma$.
       We conclude that $\nabla \Psi$ is a vector field on $\Sigma$, pointing inwards along its boundary with only discrete zeros.
       The Hopf index theorem \cite[Chapt. 3, Sec. 5]{guillemin_pollack} implies that the algebraic count of these zeros equals the Euler characteristic of the cylinder and thus
       \[
         \sum_{p :\nabla \Psi(p)=0}\mathrm{ index}(\nabla \Psi,p) = 0.
       \]
       Since isolated maxima have a positive index $+1$, then there exist at least a critical point with negative index.
       \end{proof}

       \begin{remark}\label{rem:mp}
    The critical point with negative index, whose existence is proved in Proposition \ref{prop:critici_Psi}, can be found with a max-min construction, as in the classical Poincaré-Birkhoff theorem \cite[Theorem 6.2]{Tabbook}.
    Let $p= (\xi,\eta)$ be an isolated maximum point of $\Psi$, and let $p' = (\eta,\xi)$ be its symmetric. 
    Consider the family of continuous curves
    \[
    \mathcal C = \{\gamma:[0,1]\to \mathcal{U}: \gamma(0) = p, \gamma(1)=p'\},
    \] 
    and define the following critical value
    \[
    	c^* = \sup_{\gamma \in \mathcal{C}}\, \min_{u \in[0,1]} \Psi\circ \gamma (u).
    \]
    Since the flow of $\nabla\Psi$ points away from the boundary of the compact cylinder $\Sigma$, the critical value $c^*$ is achieved in $\mathcal{U}$.
\end{remark}

 \begin{proof}[Proof of Lemma \ref{lemma:grado_psi_A}]
    The maps $\Psi_a$ and $\Psi_c$ are $C^1$ and can be written in terms of $\Psi$ and $\Psi^*$, indeed
\begin{equation}\label{eq:def_psiAeC}
     \Psi_a(\xi,\eta) = 
     \begin{cases}
     \Psi (\xi,\eta), &\hspace{-0.3cm}\text{if }\det (\gamma(\xi),\gamma(\eta))\geq 0,\\
    2\Psi^*(\xi,\eta), &\hspace{-0.3cm}\text{if }\det (\gamma(\xi),\gamma(\eta))< 0,
    \end{cases}
         \; \Psi_c(\xi,\eta) = 
     \begin{cases}
     2\Psi^* (\xi,\eta), &\hspace{-0.3cm}\text{if }\det (\gamma(\xi),\gamma(\eta))\geq 0,\\
    \Psi(\xi,\eta), &\hspace{-0.3cm}\text{if }\det (\gamma(\xi),\gamma(\eta))< 0.
    \end{cases}
\end{equation}     
     Hence, $(\hat{\xi},\hat{\eta})$ is a critival point of $\Psi_a$ and $\Psi_c$, and it is isolated. 
     Recall that for a piecewise $C^1$ map the index of an isolated critical point can be computed as the winding number of the image through its gradient of a small circle enclosing it and free of other critical points. Let us then pick a small circle $c:[0,2 \pi] \to \mathbb{S}^1\times \mathbb{S}^1$ enclosing $(\hat{\xi},\hat{\eta})$ and no other critical point of neither $\Psi$ or $\Psi^*$. We write $c$ as the concatenation of two curves  $c_1$ and $c_2$ joining antipodal points. It follows that
            \[
            \mathrm{index}(\nabla \Psi,(\hat{\xi},\hat{\eta})) = 
            \frac{1}{2\pi}\int_{\nabla \Psi\circ c} d \theta  = \frac{1}{2 \pi} \left(\int_{\nabla \Psi\circ c_1} d \theta  +\int_{\nabla \Psi\circ c_2} d \theta \right), 
            \]
            where $d \theta \in \Omega^1(\mathbb{R}^2\setminus \{0\})$ is the closed one form 
            \[
              d \theta  = \frac{ -y dx+ x dy}{x^2+y^2}.
             \]
\begin{figure}
\centering
\begin{tikzpicture}[x=0.75pt,y=0.75pt,yscale=-1,xscale=1]

\draw [color={rgb, 255:red, 74; green, 74; blue, 74 }  ,draw opacity=1 ] [dash pattern={on 4.5pt off 4.5pt}]  (181.33,45) .. controls (165.33,78) and (176.33,73) .. (161.33,116) .. controls (146.33,159) and (126.33,154) .. (102.33,202) ;
\draw  [dash pattern={on 0.84pt off 2.51pt}]  (102.33,202) -- (93.33,217) ;
\draw  [dash pattern={on 0.84pt off 2.51pt}]  (181.33,45) -- (194.33,34) ;
\draw  [fill={rgb, 255:red, 0; green, 0; blue, 0 }  ,fill opacity=1 ] (151,132) .. controls (151,133.66) and (152.34,135) .. (154,135) .. controls (155.66,135) and (157,133.66) .. (157,132) .. controls (157,130.34) and (155.66,129) .. (154,129) .. controls (152.34,129) and (151,130.34) .. (151,132) -- cycle ;
\draw  [line width=1.5]  (123.33,162) .. controls (108.33,150) and (85.33,147) .. (103.33,125) .. controls (121.33,103) and (144.33,99) .. (163,86) .. controls (181.67,73) and (225.33,89) .. (203.33,141) .. controls (181.33,193) and (138.33,174) .. (123.33,162) -- cycle ;
\draw  [fill={rgb, 255:red, 0; green, 0; blue, 0 }  ,fill opacity=1 ] (135.57,108.82) -- (115.47,113.12) -- (128.16,96.94) -- (123.67,108) -- cycle ;
\draw  [fill={rgb, 255:red, 0; green, 0; blue, 0 }  ,fill opacity=1 ] (201.9,125.83) -- (208.44,106.34) -- (215.89,125.5) -- (208.67,116) -- cycle ;
\draw  [line width=1.5]  (372.33,208.67) .. controls (357.33,196.67) and (365.33,165) .. (383.33,143) .. controls (401.33,121) and (384.67,85.67) .. (403.33,72.67) .. controls (422,59.67) and (518.33,124.67) .. (496.33,176.67) .. controls (474.33,228.67) and (387.33,220.67) .. (372.33,208.67) -- cycle ;
\draw [color={rgb, 255:red, 74; green, 74; blue, 74 }  ,draw opacity=1 ] [dash pattern={on 4.5pt off 4.5pt}]  (490.33,49.67) .. controls (458.33,91.67) and (441.33,96.67) .. (433.33,150.67) .. controls (425.33,204.67) and (436.33,209.67) .. (476.33,221.67) .. controls (516.33,233.67) and (570.33,213.67) .. (581.33,199.67) ;
\draw  [dash pattern={on 0.84pt off 2.51pt}]  (596.33,187.67) -- (581.33,199.67) ;
\draw  [dash pattern={on 0.84pt off 2.51pt}]  (490.33,49.67) -- (504.33,33.67) ;
\draw    (334.33,150.67) -- (585.33,150.67) ;
\draw [shift={(587.33,150.67)}, rotate = 180] [color={rgb, 255:red, 0; green, 0; blue, 0 }  ][line width=0.75]    (10.93,-3.29) .. controls (6.95,-1.4) and (3.31,-0.3) .. (0,0) .. controls (3.31,0.3) and (6.95,1.4) .. (10.93,3.29)   ;
\draw    (433.33,273.67) -- (433.33,29.67) ;
\draw [shift={(433.33,27.67)}, rotate = 90] [color={rgb, 255:red, 0; green, 0; blue, 0 }  ][line width=0.75]    (10.93,-3.29) .. controls (6.95,-1.4) and (3.31,-0.3) .. (0,0) .. controls (3.31,0.3) and (6.95,1.4) .. (10.93,3.29)   ;
\draw  [line width=1.5]  (456.33,118.67) .. controls (441.33,106.67) and (414.33,150.67) .. (406.33,136.67) .. controls (398.33,122.67) and (416.33,105.67) .. (456.33,91.67) .. controls (496.33,77.67) and (486.33,149.67) .. (451.33,211.67) .. controls (416.33,273.67) and (471.33,130.67) .. (456.33,118.67) -- cycle ;
\draw  [fill={rgb, 255:red, 0; green, 0; blue, 0 }  ,fill opacity=1 ] (382,136.2) -- (392.46,118.5) -- (395.75,138.8) -- (390.67,128) -- cycle ;
\draw  [fill={rgb, 255:red, 0; green, 0; blue, 0 }  ,fill opacity=1 ] (488.5,123) -- (479.8,141.63) -- (474.56,121.75) -- (480.67,132) -- cycle ;
\draw  [fill={rgb, 255:red, 155; green, 155; blue, 155 }  ,fill opacity=1 ] (451,92) .. controls (451,93.66) and (452.34,95) .. (454,95) .. controls (455.66,95) and (457,93.66) .. (457,92) .. controls (457,90.34) and (455.66,89) .. (454,89) .. controls (452.34,89) and (451,90.34) .. (451,92) -- cycle ;
\draw  [fill={rgb, 255:red, 155; green, 155; blue, 155 }  ,fill opacity=1 ] (448,213) .. controls (448,214.66) and (449.34,216) .. (451,216) .. controls (452.66,216) and (454,214.66) .. (454,213) .. controls (454,211.34) and (452.66,210) .. (451,210) .. controls (449.34,210) and (448,211.34) .. (448,213) -- cycle ;
\draw  [fill={rgb, 255:red, 155; green, 155; blue, 155 }  ,fill opacity=1 ] (167,83) .. controls (167,84.66) and (168.34,86) .. (170,86) .. controls (171.66,86) and (173,84.66) .. (173,83) .. controls (173,81.34) and (171.66,80) .. (170,80) .. controls (168.34,80) and (167,81.34) .. (167,83) -- cycle ;
\draw  [fill={rgb, 255:red, 155; green, 155; blue, 155 }  ,fill opacity=1 ] (123,164) .. controls (123,165.66) and (124.34,167) .. (126,167) .. controls (127.66,167) and (129,165.66) .. (129,164) .. controls (129,162.34) and (127.66,161) .. (126,161) .. controls (124.34,161) and (123,162.34) .. (123,164) -- cycle ;

\draw (145,191) node [anchor=north west][inner sep=0.75pt]   [align=left] {$\displaystyle \{\theta ( \xi ,\eta ) \ =\ \pi $$\displaystyle \}$};
\draw (143,140) node [anchor=north west][inner sep=0.75pt]   [align=left] {$\displaystyle (\hat{\xi } ,\hat{\eta })$};
\draw (84,104.4) node [anchor=north west][inner sep=0.75pt]    {$c_{1}$};
\draw (205.33,144.4) node [anchor=north west][inner sep=0.75pt]    {$c_{2}$};
\draw (392,156) node [anchor=north west][inner sep=0.75pt]   [align=left] {$\displaystyle ( 0,0)$};
\draw (440,233) node [anchor=north west][inner sep=0.75pt]   [align=left] {$\displaystyle \nabla \Psi (\{\theta ( \xi ,\eta ) \ =\ \pi \})$};
\draw (308,163.4) node [anchor=north west][inner sep=0.75pt]    {$\nabla \Psi \circ c$};
\draw (477,77.4) node [anchor=north west][inner sep=0.75pt]    {$\nabla \Psi ^{*} \circ c$};
\end{tikzpicture}
\caption{An illustration of the proof of Lemma  
\ref{lemma:grado_psi_A}. The dotted curve in the left picture is the graph of $\eta=\eta(\xi)$ such that $\gamma(\xi)$ and $\gamma(\eta(\xi))$ are antipodal on $\partial \Omega$. Its image trough $\nabla \Psi$ and $\nabla \Psi^*$ is the same. For this reason its intersections with both $\nabla \Psi \circ c$ and $\nabla \Psi ^{*} \circ c$ coincide.}
\end{figure}
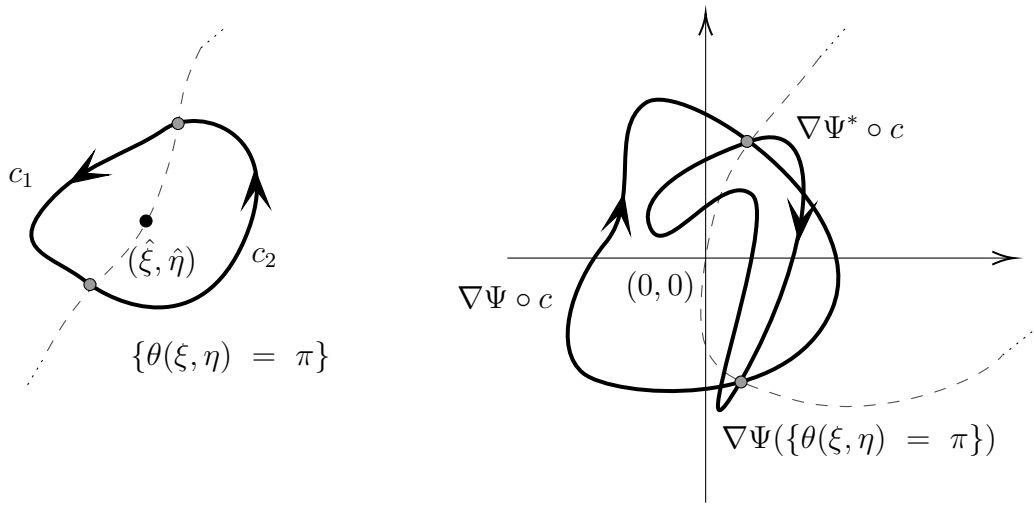
             
Since, by assumption, $\mathrm{index}(\nabla \Psi^*,(\hat{\xi},\hat{\eta})) =0$, we have
\begin{align*}
\mathrm{index}(\nabla \Psi,(\hat{\xi},\hat{\eta})) &=  \mathrm{index}(\nabla \Psi,(\hat{\xi},\hat{\eta}))+ \mathrm{index}(\nabla \Psi^*,(\hat{\xi},\hat{\eta})) \\
                 &=\frac{1}{2 \pi}\left(\int_{\nabla \Psi\circ c_1} d \theta  +\int_{\nabla \Psi\circ c_2} d \theta +\int_{\nabla \Psi^*\circ c_1} d \theta  +\int_{\nabla \Psi^*\circ c_2} d \theta\right)\\
                 &=\frac{1}{2 \pi} \left(\int_{\nabla \Psi_a\circ c} d \theta  +\int_{\nabla \Psi_c\circ c} d \theta\right)\\
                 &= \mathrm{index}(\nabla \Psi_a,(\hat{\xi},\hat{\eta}))+ \mathrm{index}(\nabla \Psi_c,(\hat{\xi},\hat{\eta})).
            \end{align*}
            It follows that at least one of the two last indices is non-zero and the proof is concluded.
        \end{proof}